\documentclass[a4paper,12pt]{article}
\usepackage[letter paper, margin=1 in]{geometry}
\usepackage{hyperref}
\usepackage{xcolor}
\usepackage{subcaption}
\usepackage{authblk}
\usepackage[numbers]{natbib}
\usepackage[utf8]{inputenc}
\usepackage{graphicx}
\usepackage{amsmath}
\usepackage{amssymb}
\usepackage[version=4]{mhchem}
\usepackage{siunitx}
\usepackage{longtable,tabularx}
\usepackage{microtype}
\usepackage{subcaption}
\usepackage{booktabs} 
\usepackage{multirow}
\usepackage{threeparttable}
\usepackage{siunitx}
\usepackage{float}
\usepackage[algo2e,ruled,vlined,linesnumbered]{algorithm2e}
\usepackage{amsfonts}
\usepackage{bm,bbm}
\usepackage{cleveref}
\usepackage{algorithm}
\usepackage{algorithmic}
\usepackage[algo2e,ruled,vlined,linesnumbered]{algorithm2e}
\usepackage{subcaption}
\usepackage{url}
\newcommand{\mbb}{\mathbb}
\newcommand{\mbf}{\mathbf}
\newcommand{\mcl}{\mathcal}
\newcommand{\bs}{\boldsymbol}

\newcommand{\T}{\textnormal}
\newcommand{\x}{\mathbf{x}}

\newcommand{\X}{\bm{\mathcal{X}}}
\newcommand{\D}{\bm{\mathcal{D}}}

\usepackage[normalem]{ulem} 
\newcommand{\argmax}{\operatornamewithlimits{arg\,max}}

\newcommand{\bo}{\texttt{BO}}
\newcommand{\bfgs}{\texttt{L-BFGS-B}}
\newcommand{\nm}{\texttt{NM}}
\newcommand{\cobyla}{\texttt{COBYLA}}
\newcommand{\tc}{\texttt{Trust-constrained}}
\newcommand{\tnc}{\texttt{TNC}}
\newcommand{\slsqp}{\texttt{SLSQP}}
\newcommand{\change}{\textcolor{black}}
\newcommand{\changetwo}{\textcolor{black}}
\newcommand{\changethree}{\textcolor{black}}

\title{Derivative-free optimization is competitive for aerodynamic design optimization in moderate dimensions}
\author[1]{Punya Plaban\footnote{Graduate assistant, Aerospace Engineering, 556 White Course drive.}}
\author[1]{Peter Bachman\footnote{Undergraduate research assistant, Aerospace Engineering, 556 White Course drive.}}
\author[1,2,*]{Ashwin Renganathan\footnote{Assistant professor, Aerospace Engineering, 556 White Course drive.
}}
\affil[1]{The Pennsylvania State University, University Park, PA, 16802}
\affil[2]{Penn State Institute of Computational and Data Sciences, University Park, PA, 16802}
\affil[*]{Corresponding author. \href{mailto:ashwin.renganathan@psu.edu}{ashwin.renganathan@psu.edu}}
\date{}

\begin{document}

\maketitle

\begin{abstract}
Aerodynamic design optimization is an important problem in aircraft design that depends on the interplay between a numerical optimizer and a \change{flow physics} solver. Derivative-based, first and (quasi) second order, optimization techniques are the de facto choice, particularly given the availability of the adjoint method and its ability to efficiently compute gradients at the cost of just one solution of the forward problem. However, the implementation of the adjoint method requires careful mathematical treatment, and its sensitivity to changes in mesh quality limits its widespread applicability. Derivative-free approaches are often overlooked for large scale optimization, citing their lack of scalability in higher dimensions and/or the lack of practical interest in globally optimal solutions that they often target. However, breaking free from an adjoint solver can be paradigm-shifting in broadening the applicability of aerodynamic design optimization. \change{We compare derivative-free and derivative-based optimizers on a suite of two- and three-dimensional aerodynamic shape optimization problems with \textbf{design-space dimensionalities ranging from $\mathcal{O}(10)$ to $\mathcal{O}(50)$ shape parameters}, including unconstrained and constrained variants of the NACA0012, RAE2822, and ONERAM6 configurations. Our results demonstrate that, for these problems, derivative-free methods are competitive with derivative-based methods and can outperform them in the higher-dimensional cases.}
Our results demonstrate that derivative-free methods are competitive with derivative-based methods, while outperforming them consistently in the high-dimensional setting. These findings highlight the practical \changethree{viability} of modern derivative-free strategies, offering a \change{promising} alternative for aerodynamic design optimization when adjoint-based gradients are unavailable or unreliable.  
\end{abstract}


\section{Introduction}
\label{sec:intro}
The environmental impact of aircraft, via reduced fuel burn, is achieved primarily by optimizing the shapes of the aircraft wing and control surfaces for better aerodynamic efficiency (e.g., higher lift-drag ratios). Aerodynamic shape optimization (ASO), which refers to optimizing parameterized shapes of aerodynamic bodies for aerodynamic performance under various constraints, is a popular and active research area~\cite{marchildon2024}. Solving this optimization problem with a high-fidelity computational fluid dynamics (CFD) model of the aircraft is computationally demanding; this is further exacerbated by the dimensionality of the parameterization. Derivative-based optimization, with sensitivities computed via the adjoint method~\cite{jameson1988control}, makes high-fidelity gradient computation feasible at only the additional cost of one primal (forward) solution. Therefore, the derivative-based approach to ASO is considered the most scalable approach and, hence, is the gold standard. However, this consideration is underpinned by the assumption that a robust adjoint formulation and an associated adjoint solver are available. 

First, implementation of the adjoint method for high-fidelity CFD solvers requires careful mathematical treatment and an involved numerical implementation. Second, convergence of adjoint solvers can be very sensitive to mesh quality and changes in boundary conditions~\cite{marchildon2024, meo2024}. Third, the adjoint solve still costs as much as the forward solve—practically speaking, derivative-based algorithms require multiple gradient evaluations per step to compute appropriate step lengths or for Hessian approximation, which can get prohibitively expensive. \citet{renganathan2021enhanced} showed that gradient-free approaches can be $10\times$ computationally cheaper than adjoint-based approaches for ASO. Finally, in the presence of non-computable constraints, such as those arising from numerical divergence, derivative-based methods often lack robustness~\cite{jim2021}. Given these challenges, despite the benefits of adjoint-based approaches, we are interested in exploring derivative-free approaches to ASO that could potentially widen its applicability. 

One case in point is ASO with large eddy simulations (LES) which better captures high angle-of-attack flows (e.g., takeoff and landing) where adjoints are inherently stochastic -- approximations may be required such as using time-averaged states to compute the adjoints. Another example is chemically reacting flows in the hypersonic regime where the time scales for flow, chemistry, and heat transfer can be highly disparate, making adjoint computation challenging. While the need for derivative-free ASO is compelling, the question that remains is: \emph{are derivative-free approaches competitive for ASO?} We intend to address this question systematically in this work. \changethree{Throughout this paper, we use ``competitive'' to mean competitive in terms of the best
feasible objective value achieved for a fixed CFD evaluation budget (and its variability
across runs/initializations), rather than asymptotic stationarity or first-order optimality
guarantees.}

Generally speaking, derivative-free methods, such as Nelder-Mead~\cite{nelder1965simplex}, COBYLA~\cite{powell1994direct}, and DIRECT~\cite{jones1993direct}, are at a disadvantage due to the additional work they have to do to identify descent directions in the absence of derivative information. However, surrogate-based ``global optimization'' approaches, such as Bayesian optimization~\cite{rasmussen2006gpml, frazier2018tutorial, garnett2023bo}, exploit knowledge gained on the design space through a data-driven global surrogate model. Our focus here is primarily on pitting a surrogate-based approach, such as \texttt{BO}, against the classic adjoint-based method for ASO; however, we also include classic derivative-free methods, Nelder-Mead and COBYLA, in this study.

\subsection{\change{Surrogate based approaches}}
Surrogate-based approaches to aerodynamic optimization are not new. Arguably, Gaussian process (GP) regression~\cite{rasmussen2006gpml} is the most popular choice of surrogates, which, along with a Bayesian decision-theoretic framework, leads to what is called ``Bayesian optimization''--a popular derivative-free optimization method. \citet{liu2022} showed that multiobjective Bayesian optimization (MOBO) can be effectively applied to airfoil shape optimization, achieving high-quality Pareto fronts with minimal CFD evaluations. Gradient-enhanced Bayesian optimization has been shown to be competitive with quasi-Newton optimizers such as SNOPT~\cite{marchildon2024}. However, it is not clear if, in the presence of derivative information, one is better off sticking to classic derivative-based line search or trust-region methods. \citet{reist2020} highlights the potential benefit of combining global and local search strategies for better results. \texttt{BO} has also been successfully applied to conceptual aircraft design problems with hidden constraints, such as simulation failures, by leveraging machine-learning classifiers~\cite{tfaily2024}, \change{and in system} reliability analysis~\cite{renganathan2024efficient}. Other applications include active flow control, multifidelity optimization~\cite{xu2024, erhard2024, renganathan2023camera}, and chance-constrained optimization~\cite{kundu2018}.

Recent advances have integrated deep learning techniques into surrogate-based aerodynamic optimization to address challenges related to high-dimensional scaling. \citet{li2020efficient} proposed a framework combining deep convolutional generative adversarial networks (DCGAN) to generate realistic airfoil and wing shapes with a convolutional neural network (CNN)-based discriminator to filter out geometrically invalid samples -- this led to an efficient low-dimensional parametrization. Building on the potential of generative models, \citet{chen2019aerodynamic} introduced a Bézier-GAN framework that learns a low-dimensional latent space of airfoil geometries, which is then wrapped into a Bayesian optimizer. In a related effort, \citet{sheikh2022optimization} a Design-by-Morphing (DbM) strategy with a mixed-variable multi-objective Bayesian optimization (MixMOBO) algorithm. The DbM approach enables exploration across a rich design space by morphing between baseline geometries, while MixMOBO efficiently navigates this space using a portfolio of acquisition functions and batch evaluations. Another work uses partial least squares with Kriging surrogates to tackle high-dimensional scaling ~\cite{bouhlel2019python}.
\citet{queipo2005surrogate} and \citet{li2022machine} provide comprehensive reviews of surrogate-based approaches.

 Overall, we observe several gaps in the existing work. First, to our knowledge, a systematic comparison between classic derivative-based methods and derivative-free methods doesn't exist. Specifically, comparisons in terms of convergence rate, achieving optimality, and sample efficiency are lacking. Second, existing work generally applies surrogate-based methods to specific application problems. However, evaluating them on a suite of problems of varying characteristics and a consistent comparison against derivative-based methods is likely to offer better insights into the viability of derivative-free methods for aerodynamic optimization. In this work, we leverage canonical aerodynamic shape optimization problems which standardizes the benchmarking exercise.
Our objective is to address the aforementioned gaps; however, our study has limited scope that we clarify below:

\noindent {\bf Remark 1.} Our choice of derivative-based and derivative-free algorithms for benchmarking is by no means exhaustive. However, we select algorithms that are widely found in well-established open-source optimization libraries, and we ensure there is sufficient diversity in their respective methodologies.

\noindent {\bf Remark 2.} Our study is not an evaluation of global versus local optimization approaches to ASO. Rather, we are interested in the practical utility of derivative-free approaches to ASO. In this regard, we entertain both local and global derivative-free approaches in this study.

\noindent {\bf Remark 3.} As a first step, we restrict ourselves to moderate dimensions ($< 50$); higher dimensional investigation is reserved for future work.

The rest of the article is organized as follows. In \Cref{sec:aso}, we provide details of the optimization problem, aerodynamic model, and the experiment design. In \Cref{sec:results}, we present the results and the associated discussion. We conclude in \Cref{sec:conclusions} with an outlook for future work.

\section{Optimization problem}
\label{sec:aso}
We are interested in solving aerodynamic design optimization problems that fall under the following generic class of nonlinear programming
\begin{equation}
\begin{split}
\min_{\mathbf{x} \in \mcl{X} \subset \mathbb{R}^n} &f(\mathbf{x}) \\
\text{s.t.}~ & c_i(\textbf{x}) = 0,~ \forall ~i \in \mathcal{E} \\
& c_i(\textbf{x}) \geq 0, ~ \forall i \in \mathcal{I},
\end{split}
\label{eqn:nonlin_prog}
\end{equation}
where $f: \mbb{R}^n \rightarrow \mbb{R}$ is the objective function, $c_i: \mbb{R}^n \rightarrow \mbb{R}$ are constraint functions, and $\mathcal{E}$ and $\mathcal{I}$ are the sets of indices representing equality and inequality constraints, respectively, $\x \in \mcl{X}$ are the design variables and $\mcl{X} \subset \mbb{R}^n$ is a compact domain. In aerodynamic design, the objective function is the drag coefficient; the constraints involve the aerodynamic lift and moment coefficients in addition to the thickness of the geometry. The design variables $\x$ are the free-form deformation (FFD)~\cite{hsu1992direct, sederberg1986ffd} control point displacements as described in Section \ref{model_overview}. 
The aerodynamic coefficients that are part of the objectives and constraints of \eqref{eqn:nonlin_prog} are computed via computational fluid dynamics, which is explained in the following section. 

Let $\Omega \subset \mcl{X}$ denote the feasible set. Then, derivative-based methods typically generate a sequence of iterates according to 
\begin{equation}
    \mbf{x}_{k+1} = \mbf{x}_k + \alpha_k \mbf{p}_k,~k=0,1,2,\ldots
\end{equation}
where $\mbf{x}_k$ denotes the $k$th iterate, $\alpha_k$ and $\mbf{p}_k$ are the step length and search direction at the $k$th iterate, respectively, and $\mbf{x}_0$ is typically a user-provided ``guess''.  The search direction $\mbf{p}_k$ is some function of the gradient of the objective and constraints, and if $\alpha_k$ is chosen to satisfy certain ``sufficient decrease'' and ``curvature'' conditions, then this sequence of iterates is guaranteed to converge to a stationary point $\mbf{x}_* \in \Omega$, provided the objective and constraints are smooth and bounded from below~\cite{nocedal1999numerical}. However, much work goes into finding a good $\alpha_k$ (involving several evaluations of the objective, constraints and their gradients) which also impacts their convergence rate.

Derivative-free methods also generate a sequence of iterates, except the iterates may depend on computing a search direction and step length (e.g., coordinate/pattern search~\cite{lewis1999pattern}), or simplexes (e.g., Nelder-Mead~\cite{olsson1975nelder}), or surrogate models (e.g., Bayesian optimization~\cite{frazier2018tutorial}), among others. Of course, these computations don't involve any higher-order information and, as a result, have to make more queries to the objective and constraints to compensate for it. Among the derivative-free methods, Bayesian optimization is particularly unique. It first involves the development of a probabilistic surrogate model of the objective (and constraints). Then, new iterates can be \emph{optimally} chosen in the sense of maximizing the predicted utility of the new iterate. We provide a quick overview in the following section, and the reader may refer to other sources for more details~\cite{renganathan2025qpots,frazier2018tutorial,shahriari2015taking,greenhill2020bayesian}.

\subsection{Bayesian optimization with Gaussian process priors}
Amongst DFO methods, our main emphasis is on Bayesian optimization (BO) -- therefore, we provide a brief overview here.
We begin by placing a Gaussian process (GP) prior distribution on the objectives $f(\x) \sim \mcl{GP}(0, k(\x, \cdot))$, where $k(\cdot, \cdot): \X \times \X \rightarrow \mbb{R}_+$ is a covariance function (a.k.a., \emph{kernel}). We denote the observations of $f$ as $y_i = f(\x_i),\ldots,n$. We begin by fitting a posterior GP (parametrized by hyperparameters $\bs{\Omega}$) for the observations $\D_n = \{\x_i, y_i\},~i=1,\ldots,n$, which gives the conditional aposteriori distribution~\cite{rasmussen2006gpml}:
\begin{equation}
    \begin{split}
        Y(\x) |
      \D_n, \bs{\Omega} &\sim \mcl{GP}(\mu_n(\x), {\sigma}^2_n (\x)), \\
\mu_n(\x) &= \mbf{k}_n^\top \mbf{K}_n ^{-1} \mbf{y}_n\\
\sigma^2_n(\x) &=  k(\x, \x) - \mbf{k}_n^\top \mbf{K}_n ^{-1} \mbf{k}_n,
    \end{split}
    \label{e:GP}
\end{equation}
where $\mbf{k}_n \equiv k(\x, X_n)$ is a vector of covariances between $\x$ and all observed points in $\D_n$, $\mbf{K}_n \equiv k(X_n, X_n)$ is a sample covariance matrix of observed points in $\D_n$, $\mbf{I}$ is the identity matrix, and $\mbf{y}_n$ is the vector of all observations in $\D_n$. \change{ Note that we assume the aerodynamics are computed from deterministic PDE-based computer models and hence don't account for any observation noise in fitting the GP.
The GP hyperparameters (which typically include parameters of the covariance function) are estimated via maximizing the marginal log likelihood of the GP~\cite{rasmussen2006gpml}.} We also denote by $X_n = [\x_1,\ldots,\x_n]^\top \in \mbb{R}^{n\times d}$ the observation sites.  BO seeks to make sequential decisions based on a probabilistic utility function constructed from the posterior GP $u(\x) = u(Y(\x)|\D_n)$. 
Sequential decisions are the result of an ``inner'' optimization subproblem that optimizes an acquisition function, typically of the form $\alpha(\x) = \mbb{E}_{Y|\D_n}[u(Y(\x))]$. 
For example, the expected improvement (EI)~\cite{jones1998efficient} acquisition function can be expressed as
\[
\begin{split}
    \alpha_\T{EI}(\x) =& \mbb{E}_{Y|\D_n}[\max(0, \left(Y(\x) - \xi \right) )],
\end{split}
\]
where $\xi = \max_i \{f(\x_i), \ldots, f(\x_n)\}$ is the current best observed value. Although $\alpha_{\T{EI}}$ has an analytically closed and continuously differentiable form, this depends on the choice of $u$. Other alternatives to EI include, notably, the information-theoretic approaches~\cite{hernandez2016predictive,hernandez2014predictive,wang2017max}. We provide a generic \bo~algorithm in \Cref{a:BO}. \change{Extension to constrained problems abound, e.g., \cite{gardner2014bayesian,letham2019constrained}; in this work, we use a trust-region based approach to constrained BO~\cite{eriksson2019scalable}.}  \change{It is worth noting that the GP hyperparameter estimation and acquisition function optimization are additional nonconvex optimization required at every step in BO. However, this optimization depends only on the GP surrogate (and hence is cheap) and can exploit high quality derivative information. As a result, we assume that these ``inner'' optimizations within BO are carried out with sufficient accuracy and their cost is trivial compared to that of objective and constraint evaluations.}

\begin{algorithm}[t]
\textbf{Given:} $\D_n = \lbrace
\x_i, y_i \rbrace _{i=1}^n$, 
 total budget $B$, 
 and GP hyperparameters $\bs{\Omega}$ \\
\KwResult{$\{\max_{i=1,\ldots,B} f(\x_i),~\argmax_{i=1,\ldots,B} f(\x_i),\}$}
  \For{$i=n+1, \ldots, B$, }{
  Find $\x_i \in \underset{\x \in \X }{\argmax}~ \alpha(\x)$ \qquad (acquisition function maximization)\\
    Observe $y_i$ = $f(\x_i) + \epsilon_i$\\ 
    Append $\D_i = \D_{i-1} \cup \lbrace \x_i, y_i \rbrace$\\ 
    Update GP hyperparameters $\bs{\Omega}$ \\
 }
 \caption{Generic Bayesian optimization}
 \label{a:BO}
\end{algorithm}

\subsection{Overview of the aerodynamic model}\label{model_overview}
  \par In this work, we leverage SU2 for high-fidelity aerodynamic calculations.
  SU2 is a finite-volume based open-source CFD and multi-physics simulation suite developed for solving partial differential equations (PDEs) on unstructured meshes using numerical methods with adjoint capability \cite{palacios_stanford_2013}. We specifically leverage the inviscid (Euler) and Reynolds averaged Navier-Stokes (RANS) models within SU2.
  As all the test cases are in the transonic regime, the Jameson-Schmidt-Turkel (JST) scheme \cite{jameson2017origins} is used for spatial discretization, which is a central scheme with 2nd order accuracy in time and space. Furthermore, an implicit Euler scheme is used for the time marching of the solution, as it is stable for any time-step size unlike the explicit scheme until a steady state is reached. For the RANS simulations, the popular one-equation Spalart-Allmaras (SA) turbulence model \cite{Spalart1992} is used as a closure for the system of equations. 
  Slip wall boundary conditions and no-slip wall conditions are used to specify the aerodynamic surfaces in Euler simulations and RANS simulations respectively. In addition to that, a freestream boundary condition is applied to the farfield boundaries in each case and a symmetry boundary condition is used for the three-dimensional test cases to reduce the computational complexity. A constant CFL number of $50$ is used in each case without any adaptation parameter. 
  Further, W-cycle multigrid solver is used in the case of Euler simulations to accelerate convergence; however, this is not done for RANS simulations as it worsens the convergence of adjoint simulations, thereby negatively impacting the accurate computation of gradients of the objective function. Finally, we employ the minimum value of the $\log_{10}$ of the residual as the convergence criterion. This criterion is evaluated over the last $200$ iterations of each simulation to ensure stability and consistency in the convergence assessment. \Cref{Table:1} summarizes the three canonical
  geometries and their operating conditions considered in this study.
  
  \par We use free form deformation (FFD) parametrization~\cite{sederberg1986ffd} in this work. 
  The control point displacements on the FFD box are the design variables of the optimization problem; see \Cref{Mesh-Geometries} for an illustration of the FFD parameterization and we provide a few more in the Appendix. The arrangement of control points is carefully designed to avoid influencing the deformation near critical regions such as the leading edge, trailing edge, wing root, and wing tip, thereby preserving the structural fidelity of the volume mesh in these sensitive areas. Additionally, as the control point displacements deform the geometry and the mesh, we apply necessary smoothing during deformation to ensure minimal mesh quality distortion; see the bottom of \Cref{Mesh-Geometries}. 



  The three test cases used in this work are summarized in \Cref{Table:1}.  The NACA0012 mesh comprises $10216$ unstructured triangular cells without any boundary layers. The RAE2822 mesh includes $22800$ unstructured cells comprising of triangles and quadrilaterals. To capture the viscous phenomena, $26$ additional layers of quad elements are added. The first layer of the boundary layer is placed approximately $1.008 \times 10^{-5}$ m away from the wall such that the wall $y^+ \approx 1$. The 3D ONERA test case involves an unstructured mesh as well with $582752$ elements. \change{The Mach and Reynolds numbers and the angle of attack are chosen to roughly match the AIAA Aerodynamic Design Optimization Discussion Group ADODG~\footnote{\url{https://sites.google.com/view/mcgill-computational-aerogroup/adodg}} benchmark problem settings.}
\begin{table}[ht]
\centering
\caption{Summary of test cases.}
\label{Table:1}
\begin{tabular}{c c c c c c c c} 
\toprule 
{\bf Test case} & $d$ & {\bf Obj.} & {\bf Cons.} & {\bf Mach} & {\bf Re.} ($\times 10^6$) & {\bf AOA} & {\bf Physics} \\
\midrule 
NACA & $4,8,16$ & $C_d$ & - & 0.85 & $\infty$ & 0.0 & Euler \\
RAE & $4,8,16,32$ & $C_d$ & $C_m,~C_l$, thick. & 0.729 & 6.50 & 2.31 & RANS \\
ONERA & $12$ & $C_D$ & $C_L$, thick. & 0.8395 & $\infty$ & 3.06 & Euler \\
\bottomrule
\end{tabular}
\end{table}

\begin{figure}[htb!]
    \centering
    \begin{subfigure}{0.33\textwidth}
\includegraphics[width=\linewidth]{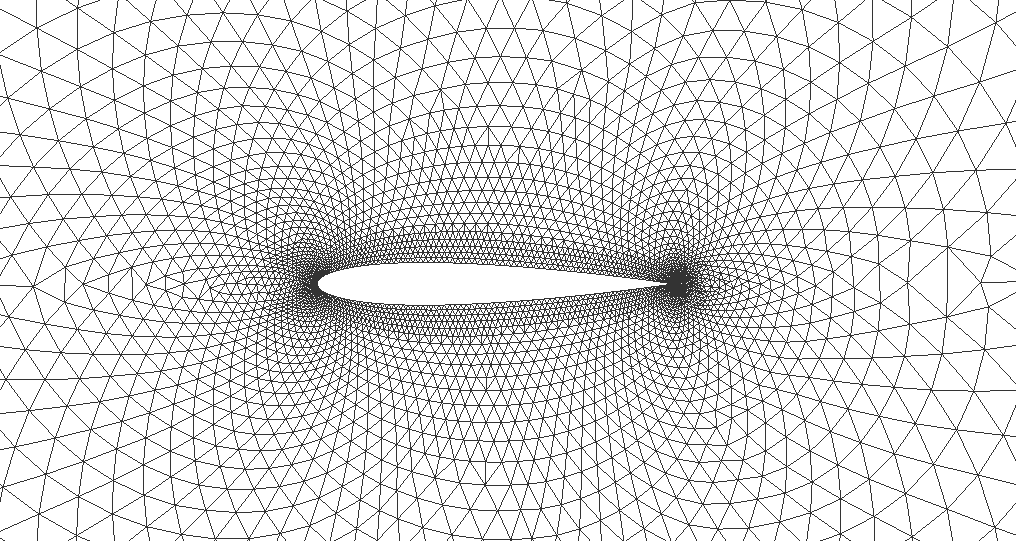}
\caption{NACA0012}\label{Mesh-NACA0012}
    \end{subfigure}%
\begin{subfigure}{0.33\textwidth}
\includegraphics[width=\linewidth]{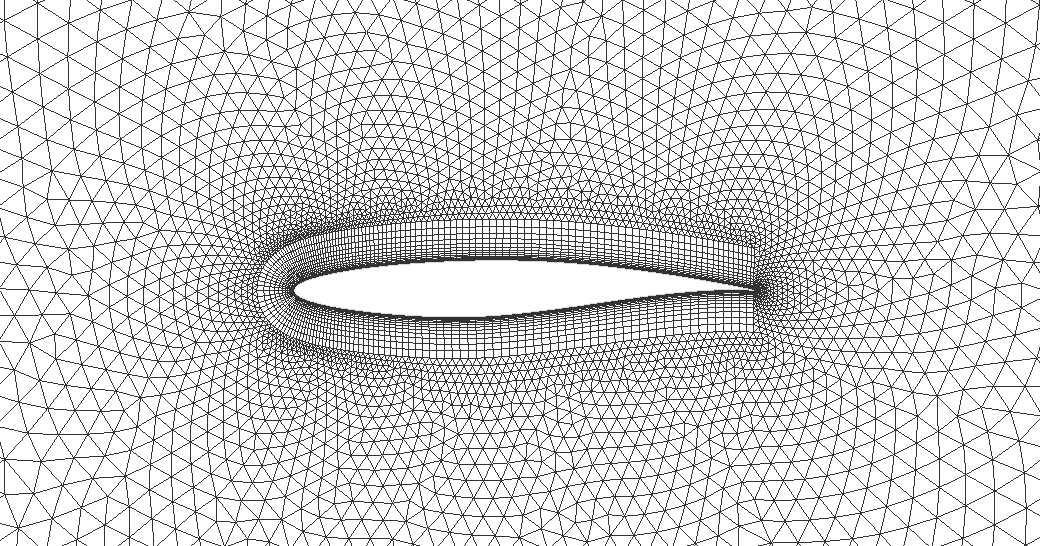} \caption{RAE2822}\label{Mesh-RAE2822}
    \end{subfigure}%
    \begin{subfigure}{0.3\textwidth}
\includegraphics[width=\linewidth]{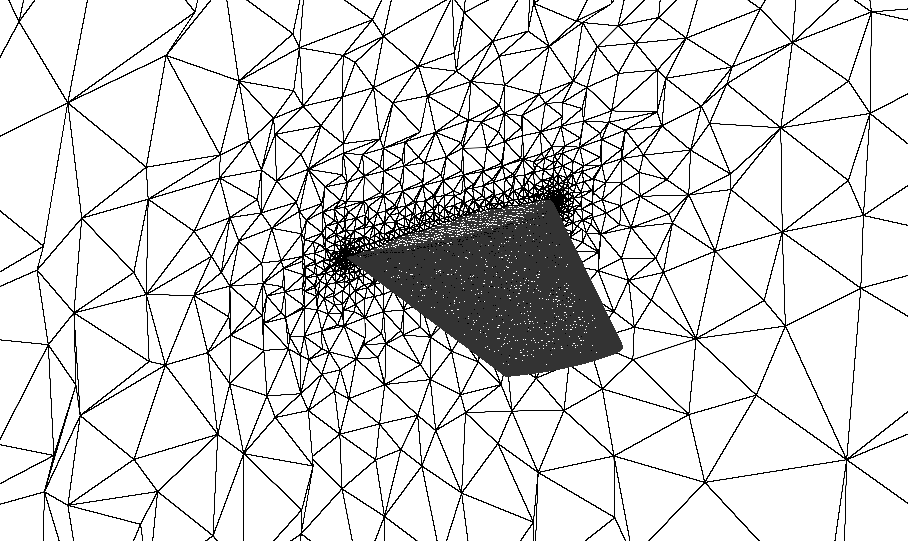} 
\caption{ONERAM6}\label{Mesh-ONERAM6}
    \end{subfigure} \\
    \begin{subfigure}{0.5\textwidth}
\includegraphics[width=\linewidth]{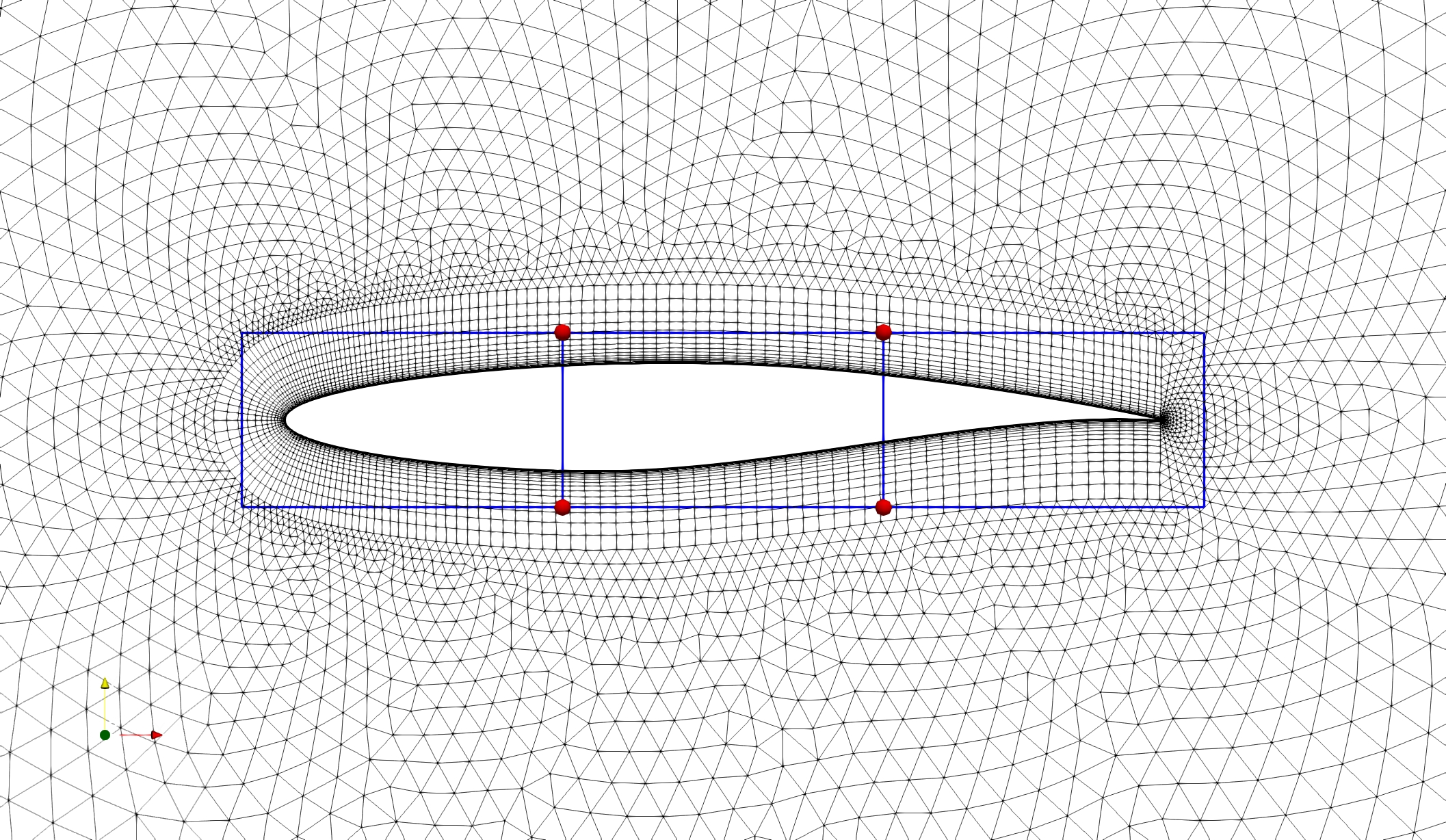}
\caption{Undeformed Mesh}\label{Undeformed-RAE2822}
    \end{subfigure}%
\begin{subfigure}{0.5\textwidth}
\includegraphics[width=\linewidth]{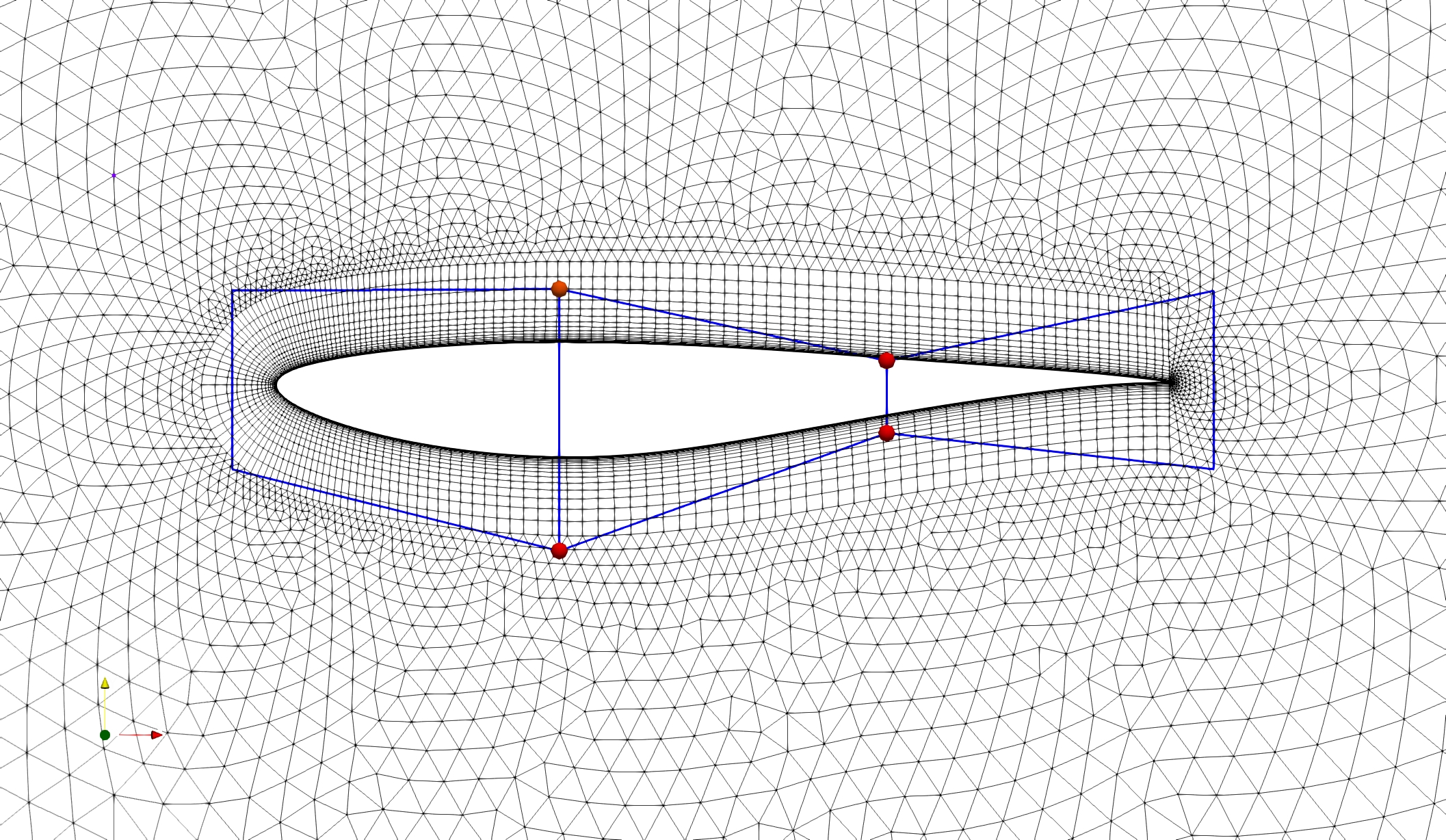} \caption{Deformed Mesh}\label{Deformed-RAE2822}
    \end{subfigure}\\
\begin{subfigure}{.5\textwidth}
        \includegraphics[width=1\linewidth]{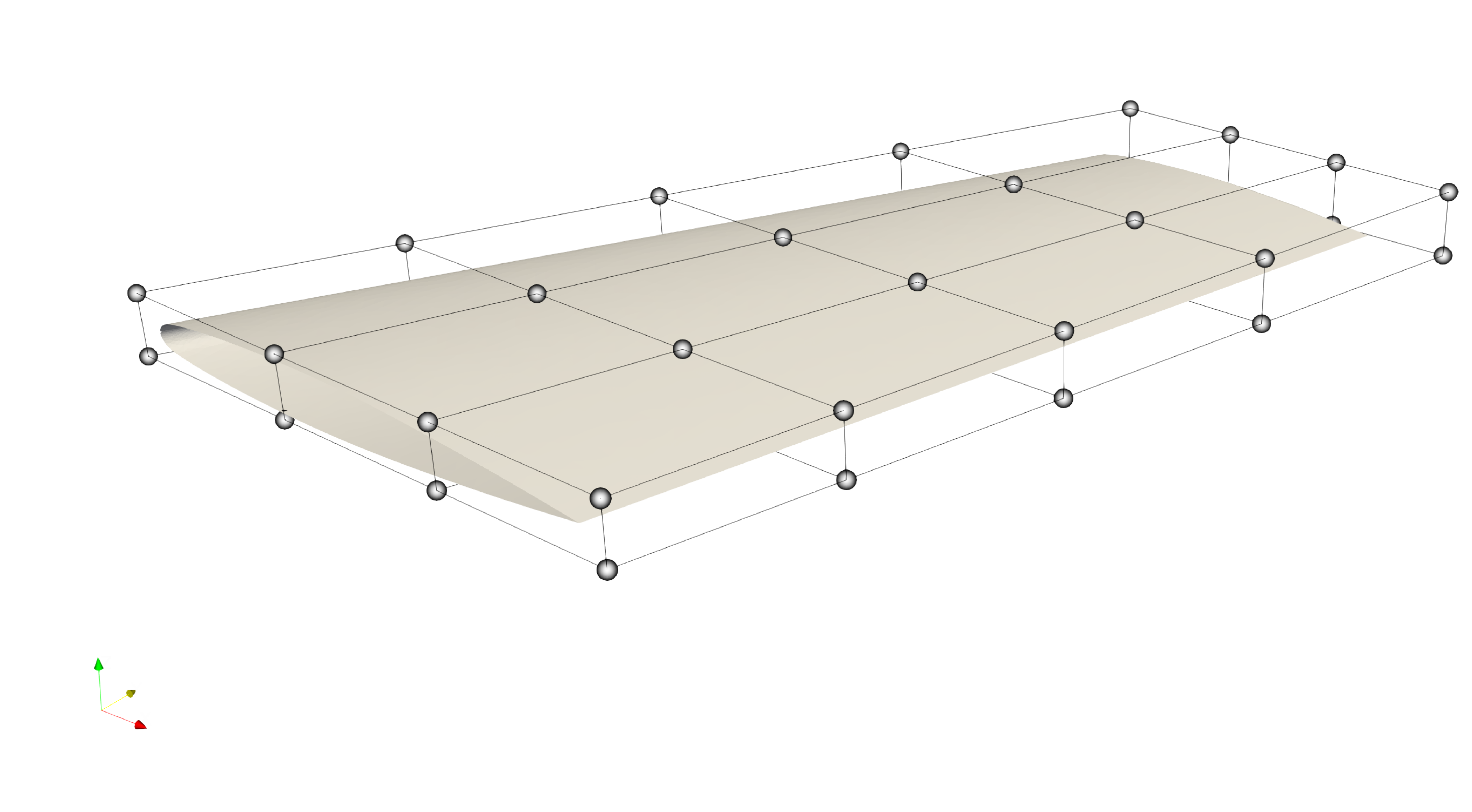}    
    \end{subfigure}%
    \begin{subfigure}{.5\textwidth}
        \includegraphics[width=1\linewidth]{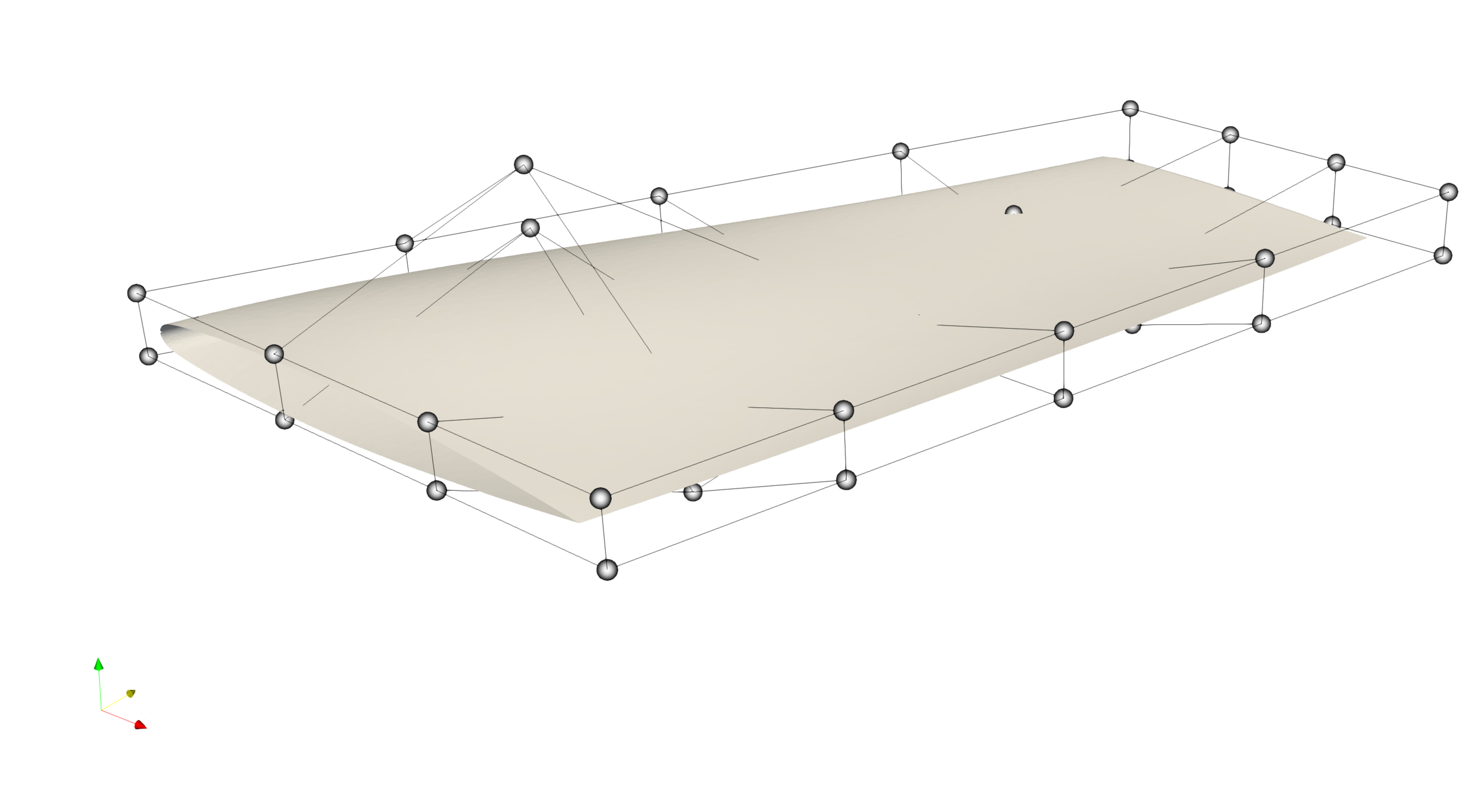}  
    \end{subfigure}    
        
    \caption{Computational mesh and free-form deformation (FFD) parametrization for shape deformation.}
    \label{Mesh-Geometries}
\end{figure}

\subsection{Experiment design}
\label{Experiments}
To answer the central question of this work, we consider a combination of widely used derivative-based and derivative-free optimization algorithms. This includes \change{a} quasi-Newton method (\texttt{L-BFGS-B}) \cite{Zhu1997}, \change{a} first-order trust-region method (\tc) \cite{Conn2000}, the truncated Newton conjugate gradient (\texttt{TNC}) \cite{Nash1984}, and sequential least squares quadratic programming (\texttt{SLSQP}) \cite{Gill1981} for derivative-based methods. On the derivative-free side, we entertain Nelder-Mead (\texttt{NM}) \cite{nelder1965simplex}, constrained optimization by linear approximation (\texttt{COBYLA}) \cite{powell1994direct}, and Bayesian optimization (\texttt{BO}) with Gaussian process priors \cite{rasmussen2006gpml, frazier2018tutorial, garnett2023bo}. \change{The choice of our algorithms, therefore, span zeroth-order (both model based and model free), first-order, and quasi second-order methods -- these choices, while not exhaustive, are representative choices of what is widely leveraged in the optimization community.}
\begin{table}[ht]
\centering
\caption{Optimization stopping criteria tolerance settings.}
\label{tab:tolerances}
\small
\begin{tabularx}{\linewidth}{l l l c}
\toprule
\textbf{Algorithm} & \textbf{Stopping crit.} & \textbf{Value} & \textbf{Description} \\
\midrule
L\mbox{-}BFGS\mbox{-}B   &  \texttt{gtol}, \texttt{maxiter} & $10^{-4},~500$  &  Projected grad. inf norm\\
Trust-constrained        & \texttt{gtol}  & $10^{-4}$ &  Lagrangian grad. norm and MCV.\\
TNC                      & \texttt{ftol}, \texttt{gtol} & $10^{-4},~10^{-3}$ & Func. precision and projected grad. inf norm \\
SLSQP                    & \texttt{ftol} & $10^{-9}$  &  Lagrangian grad. norm and MCV\\
Nelder\mbox{--}Mead      & \texttt{fatol}  & $10^{-4}$  &  Abs. func. diff. between iterations\\
COBYLA                   & \texttt{tol}, \texttt{catol}  & $10^{-9},~10^{-4}$ &  Trust-region lower bound\\
BO                       & - & - & Budget of function evaluations  \\
\bottomrule
\end{tabularx}
\end{table}

\par In all optimization runs, we use the following solver settings for each algorithm, that were held fixed across the experiments. \change{For the \texttt{L-BFGS-B} we set a projected gradient infinity norm tolerance (\texttt{gtol}) while limiting the maximum number of search steps (\texttt{maxiter}) to $500$. For \tnc, in addition to \texttt{gtol}, we also set a tolerance for relative objective function change \texttt{ftol}. The tolerance for \slsqp~(\texttt{ftol}) and \tc~(\texttt{gtol}) collectively controls Lagrangian gradient infinity norm and maximum constraint violation (MCV). Whereas for the unconstrained DFO method \nm, the tolerance (\texttt{fatol}) controls the absolute difference in the objective value between iterations, for \cobyla, \texttt{tol} and \texttt{catol} are the lower bound on the trust-region size and MCV, respectively. Finally, for \bo, the algorithm is run for a specified budget of function evaluations. The specific tolerances and their respective values for each algorithm are summarized in \Cref{tab:tolerances}.}
Additionally, we normalize the design variables to be bounded by $[0,1]$ to remove any undesirable influence of the scales of the design variables. We leverage the optimization routines within the open-source \texttt{SciPy}~\cite{virtanen2020scipy} library.

\noindent {\bf Remark 4.} Note that we don't restart the local optimization algorithms with a randomized choice of $\x_0$. Instead, $\x_0$ is chosen from the initial pool of seed points provided to \bo~that corresponds to the best feasible objective value. \changethree{This way, we ensure we don't disadvantage local optimization methods with a starting point worse than that of BO. On the other hand, BO is repeated $3$ times particularly to address the stochasticity of the seed point selection.}

\noindent \change{{\bf Remark 5.}
We choose disparate tolerances for algorithms due to the following reasons. First, \texttt{SciPy} exposes different termination metrics and semantics for each solver. For example, \slsqp~and \tc~uses a single tolerance that controls several tests simultaneously (objective change, Lagrangian gradient, and maximum constraint violation). On the other hand, \bfgs~exposes both a tolerance for relative decrease in the objective and  the projected gradient infinity norm; the two control fundamentally different notions -- function decrease vs. stationarity. Second, we had to make tweaks to algorithm tolerances independently to drive them toward convergence.
}

\noindent {\bf Remark 6.} \changetwo{We do have strong reasons for restricting our study to SciPy, as opposed to other libraries used for aerodynamic design optimization such as SNOPT~\cite{gill2005snopt}. First, SciPy is a widely adopted research and industry standard software (see, e.g., \cite{virtanen2020scipy}). Industry standard machine learning libraries still use SciPy routines under the hood (see, e.g., \cite{balandat2020botorch}). Second, SNOPT has restrictions in its usage -- for academic research and student learning, an unrestricted version is available but with a time limit~\footnote{\url{https://ccom.ucsd.edu/~optimizers/downloads/}}. SciPy, on the other hand, is openly available and is actively under maintenance and development. Finally, SNOPT presents a sparse version of an SQP algorithm with a limited-memory approximation to the Hessian of the Lagrangian -- the software is suitable for smooth nonconvex programs. The same problem class is also targeted by SciPy’s SLSQP and by SciPy’s \tc~(a trust-region constrained method that, for equality-constrained problems, implements a Byrd–Omojokun trust-region SQP~\cite{omojokun1989trust} variant and, with inequalities, switches to a trust-region interior-point approach). In other words, our SciPy baseline already exercises the same family of methods as SNOPT. For these reasons, we choose SciPy for this work that also promotes replicability of results.}

For \bo, we use the expected improvement (EI) acquisition function~\cite{jones1998efficient} for unconstrained problems and the scalable constrained Bayesian optimization (SCBO)~\cite{eriksson2019scalable} implementation for constrained problems, which combines Thompson sampling~\cite{thompson1933likelihood} with a trust-region method to select the iterates. \change{It is worth noting that all reported optima are the best observed feasible designs under SCBO; infeasible points are not candidates for convergence, but are still included in updating the GP surrogate.} All \bo~related code is implemented in the open-source library \texttt{BoTorch}~\cite{balandat2020botorch}. All the code, including the SU2 mesh and configuration files, is available publicly~\footnote{\url{https://github.com/csdlpsu/aso}} which will facilitate reproducing our results as well as further benchmarking.

For each unconstrained test case, the primary performance metrics include (i) the drag coefficient and (ii) the gradient infinity norm of the drag coefficient (for first-order optimality) versus optimization iteration and number of objective and gradient evaluations; we assume the cost of objective and gradient evaluations is roughly equal.  For the constrained test cases, we additionally monitor the maximum constraint violation--this is the infinity norm of an aggregated vector of all constraint violations evaluated at an iterate. For Bayesian optimization, we compare the ``best observed'' (that is, the smallest feasible objective observed until the current iteration) drag coefficient at every iteration. Note that, unlike all other algorithms, \bo~requires an ensemble of ``seed'' points to start the algorithm; to account for this, \bo~convergence histories are offset to the right.
We will provide further details in the following sections.

\section{Results and discussion}
\label{sec:results}

\subsection{Unconstrained NACA0012}
\begin{figure}[htb!]
    \centering
    \begin{subfigure}{0.33\textwidth}
\includegraphics[width=\linewidth]{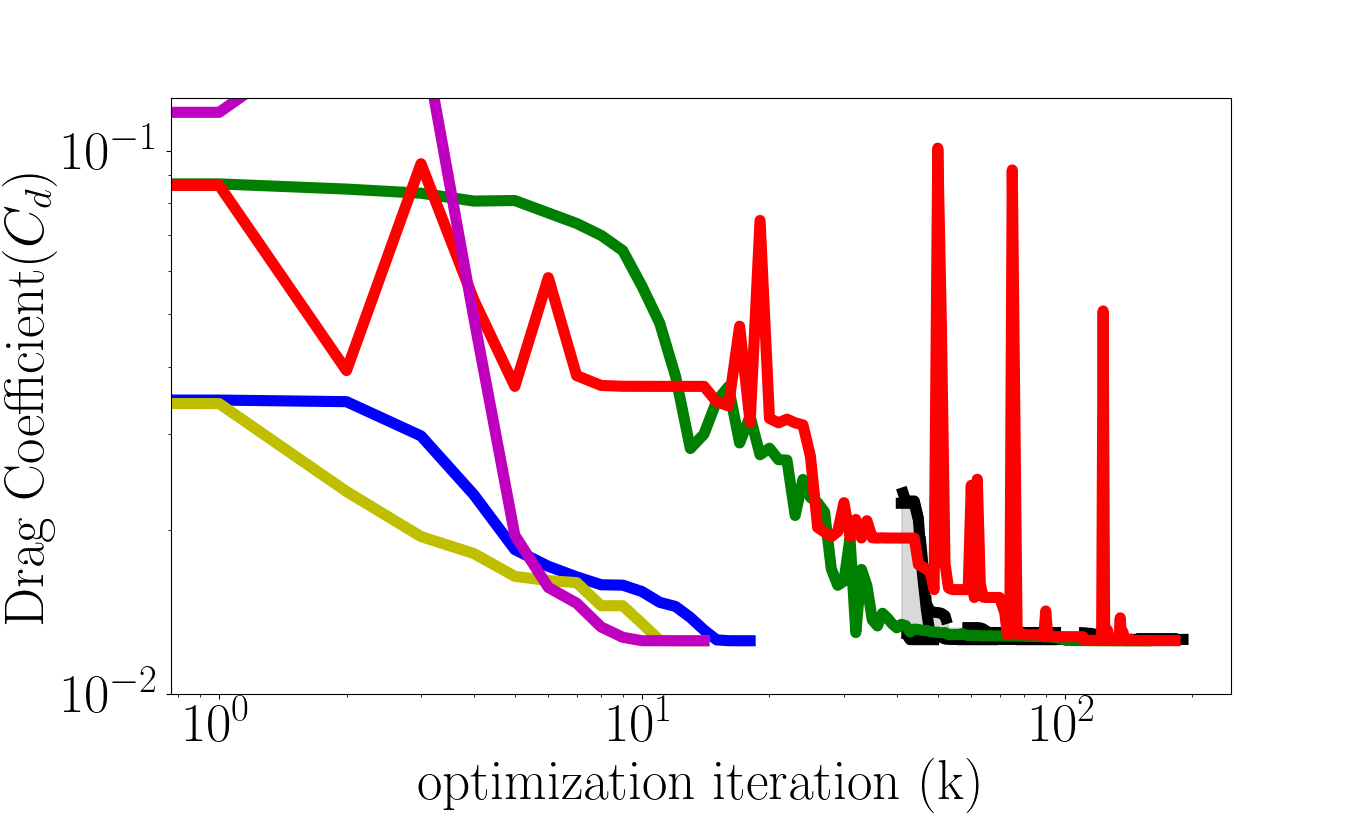}
\caption{4 control points}\label{Cd-4pts-NACA0012}
    \end{subfigure}%
\begin{subfigure}{0.33\textwidth}
\includegraphics[width=\linewidth]{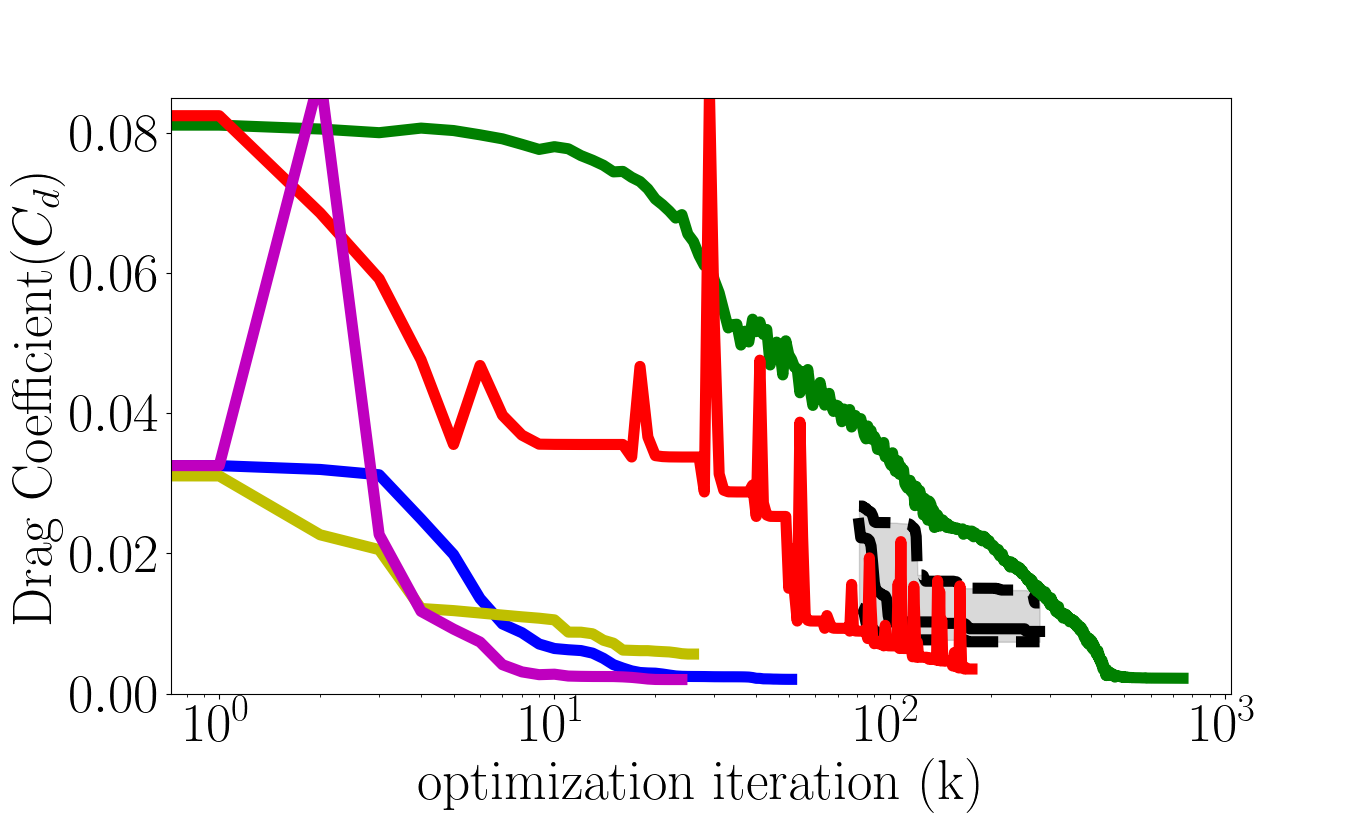} \caption{8 control points}\label{Cd-8pts-NACA0012}
    \end{subfigure}%
    \begin{subfigure}{0.33\textwidth}
\includegraphics[width=\linewidth]{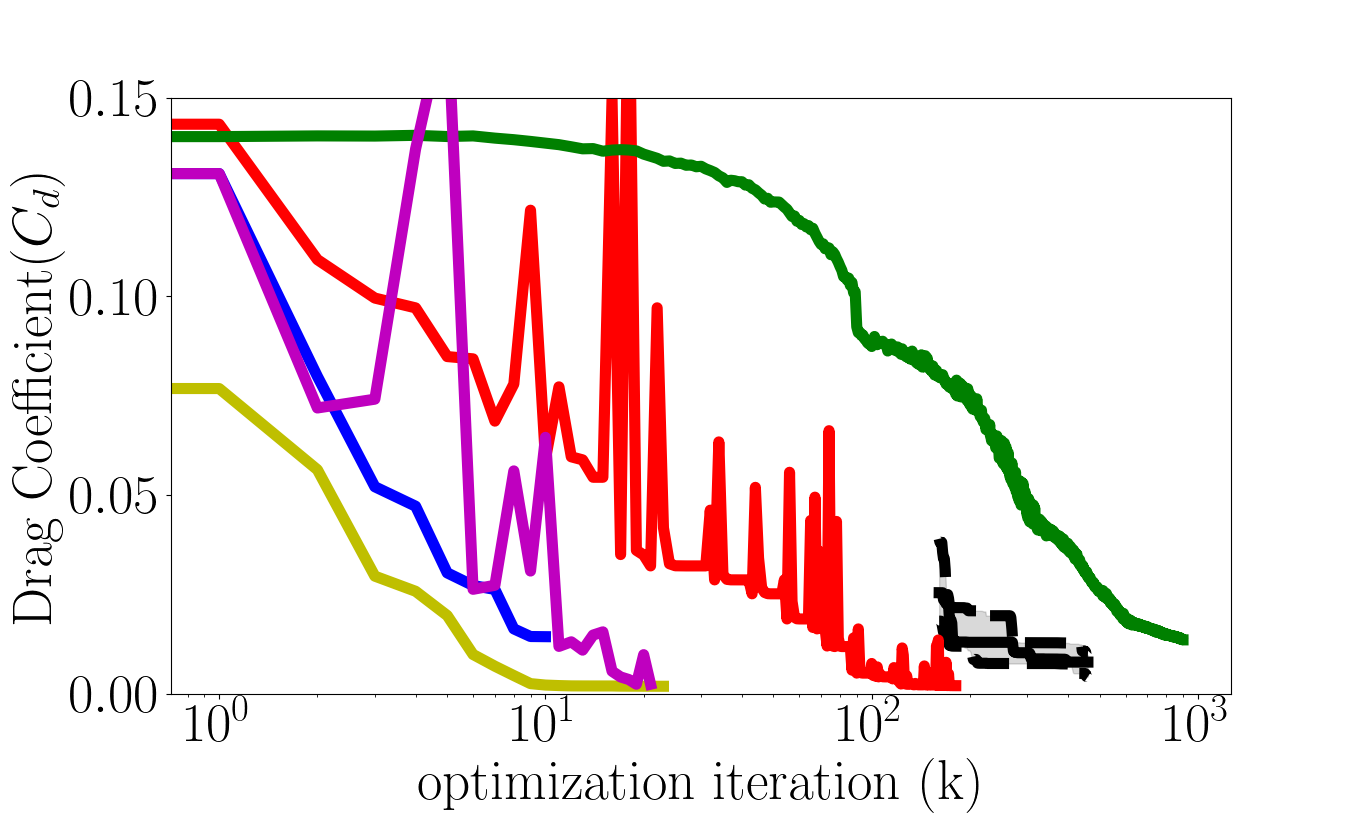} 
\caption{16 control points}\label{Cd-16pts-NACA0012}
    \end{subfigure}
    \begin{subfigure}{\textwidth}
        \centering
        \includegraphics[width=0.9\linewidth]{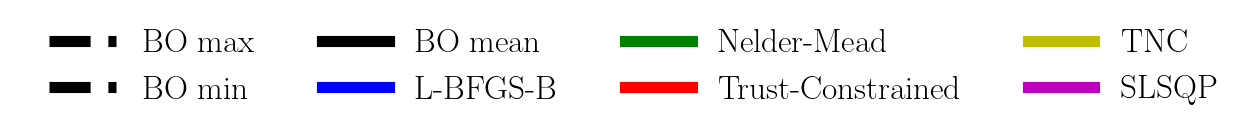} 
    \end{subfigure}
    \caption{{\bf Unconstrained NACA0012.} Convergence histories versus optimization iterations. \change{\bo~is offset to indicate that seed samples were used to start the algorithm.}}
    \label{Cd-NACA0012}
\end{figure}

{\label{NACA0012}}
We begin with the unconstrained optimization of NACA0012 under inviscid conditions. The control point displacements are contained to $[-0.1, 0.1]$; for example, with $4$ control points, the domain $\mcl{X} = [-0.1, 0.1]^4$. 
\Cref{Cd-NACA0012} shows the drag coefficient (objective) history for $4$, $8$ and $16$ control points; notice that derivative-based methods have highly varying performance -- \bfgs, \slsqp, and \tnc~rapidly converge to the minimum while \tc~are relatively slow. This exemplifies the fact that derivative-based methods can vary drastically in performance depending on the search direction $(\mathbf{p}_k)$ and step length $(\alpha_k)$ calculation. Derivative-free methods (\nm~and \bo)  take longer to converge, while still being quite comparable to \change{\tc}~in the $4$ control points case. 

\begin{figure}[htb!]
    \centering
    \begin{subfigure}{0.33\textwidth}
\includegraphics[width=\linewidth]{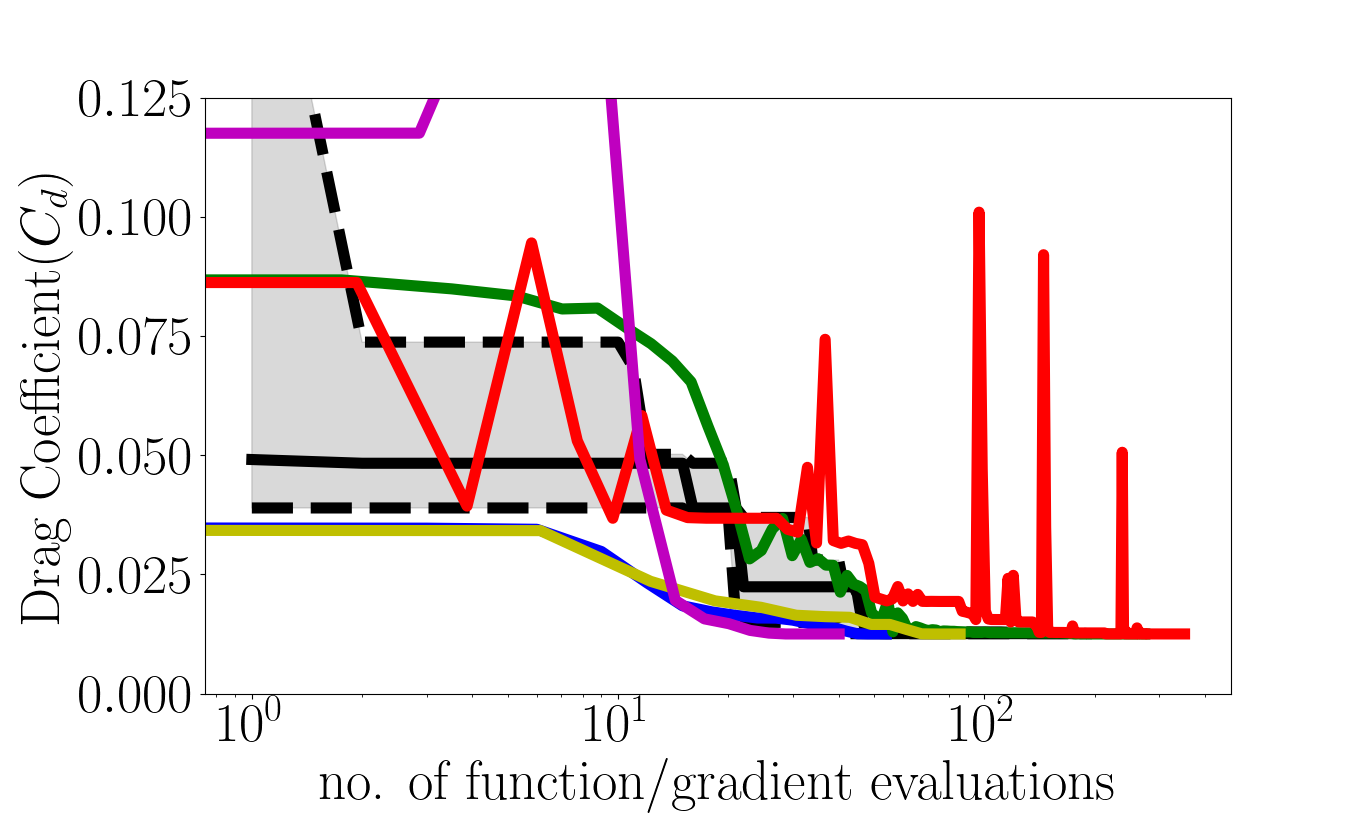}
\caption{4 control points}\label{feval-4pts-NACA0012}
    \end{subfigure}%
\begin{subfigure}{0.33\textwidth}
\includegraphics[width=\linewidth]{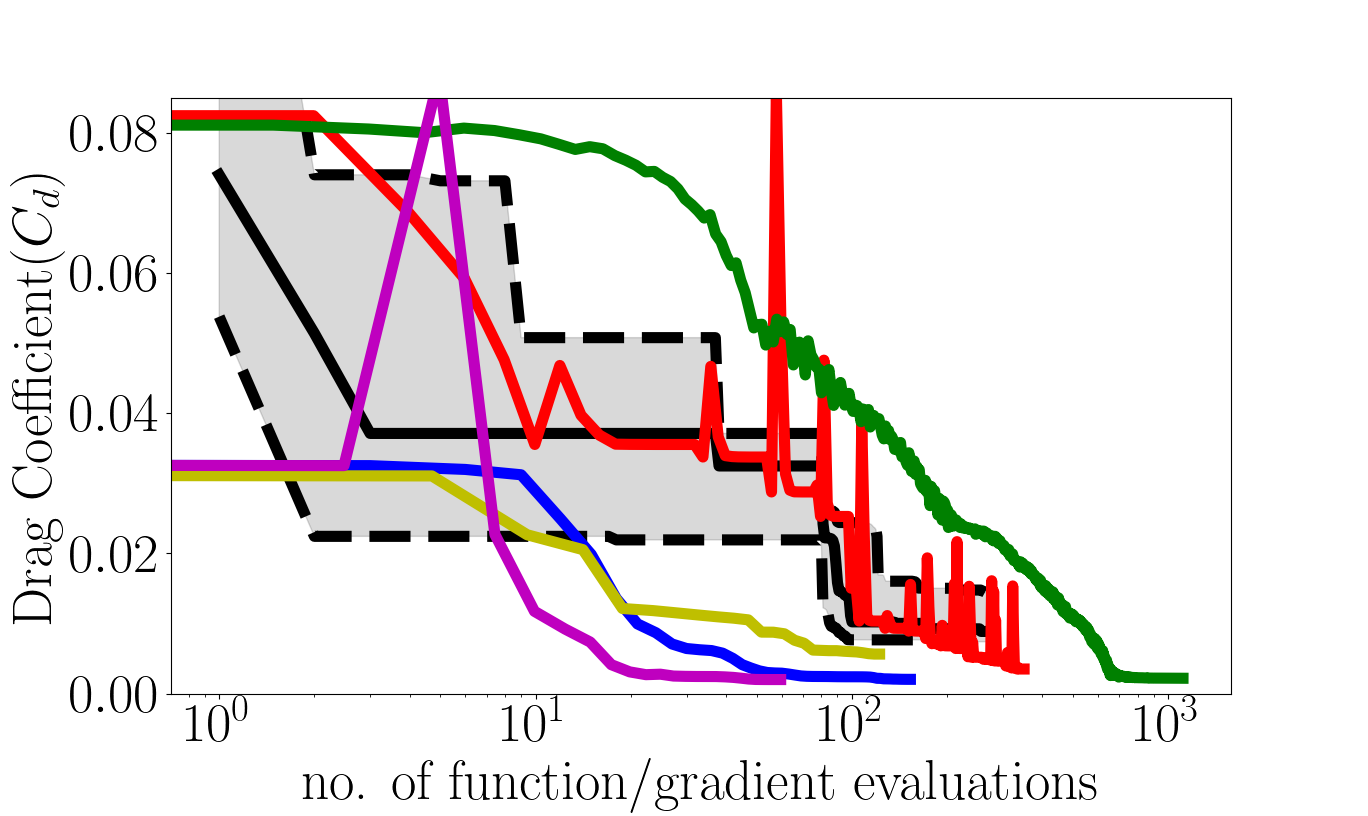} \caption{8 control points}\label{feval-8pts-NACA0012}
    \end{subfigure}%
    \begin{subfigure}{0.33\textwidth}
\includegraphics[width=\linewidth]{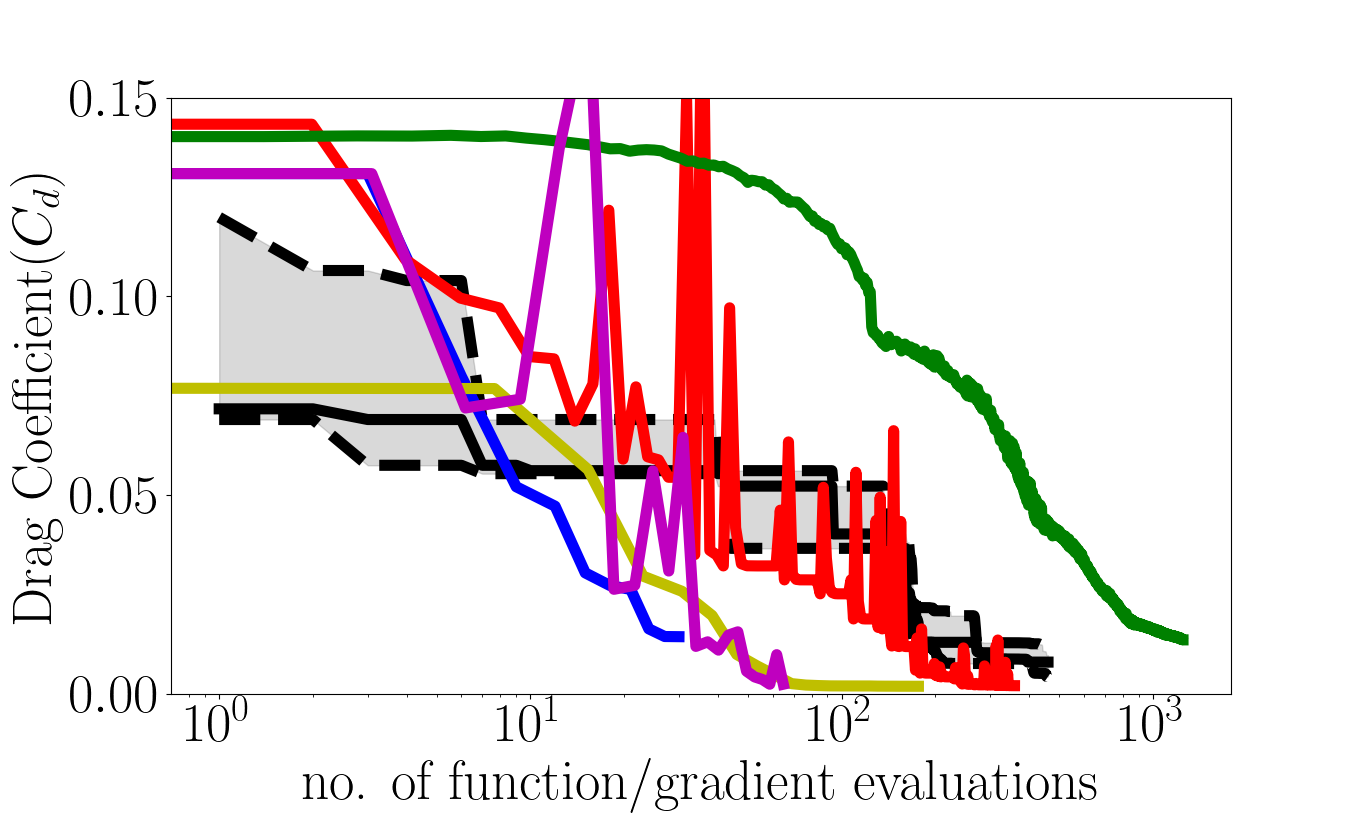} 
\caption{16 control points}\label{feval-16pts-NACA0012}
    \end{subfigure}
    \begin{subfigure}{\textwidth}
        \centering
        \includegraphics[width=0.9\linewidth]{Legends_plot.png} 
    \end{subfigure}
    \caption{{\bf Unconstrained NACA0012.} Convergence histories versus number of function/gradient evaluations.}
    \label{feval-NACA0012}
\end{figure}
Note that for \bo~, we show the $\max,~\min$, and the average of $3$ repetitions with randomized seed points -- \change{the gray shaded regions emphasize the variabilty in the \bo~algorithm realizations}. \change{Recall that the unconstrained \bo~uses the EI~\cite{jones1998efficient} acquisition function.} It is worth noting that \bo~costs one objective evaluation per iteration, in contrast to other methods, which can consume additional evaluations to compute the search step and gradient (when applicable). So, we compare convergence histories against the number of function (objective + gradient) evaluations in \Cref{feval-NACA0012}. Now notice that \bo~performs comparably to derivative-based methods \change{(particularly \tc)}, while outperforming \nm~in the cases with $8$ and $16$ control points. It should be noted that we assign equal costs for objective function and adjoint evaluation; however, in practice, the adjoint solver in SU2 tends to take significantly longer to converge, while occasionally failing to do so. Given these factors, we argue that \bo~is quite competitive for this 2D inviscid test case, despite not depending on gradient information.
\begin{figure}[htb!]
    \centering
    \begin{subfigure}{0.33\textwidth}
\includegraphics[width=\linewidth]{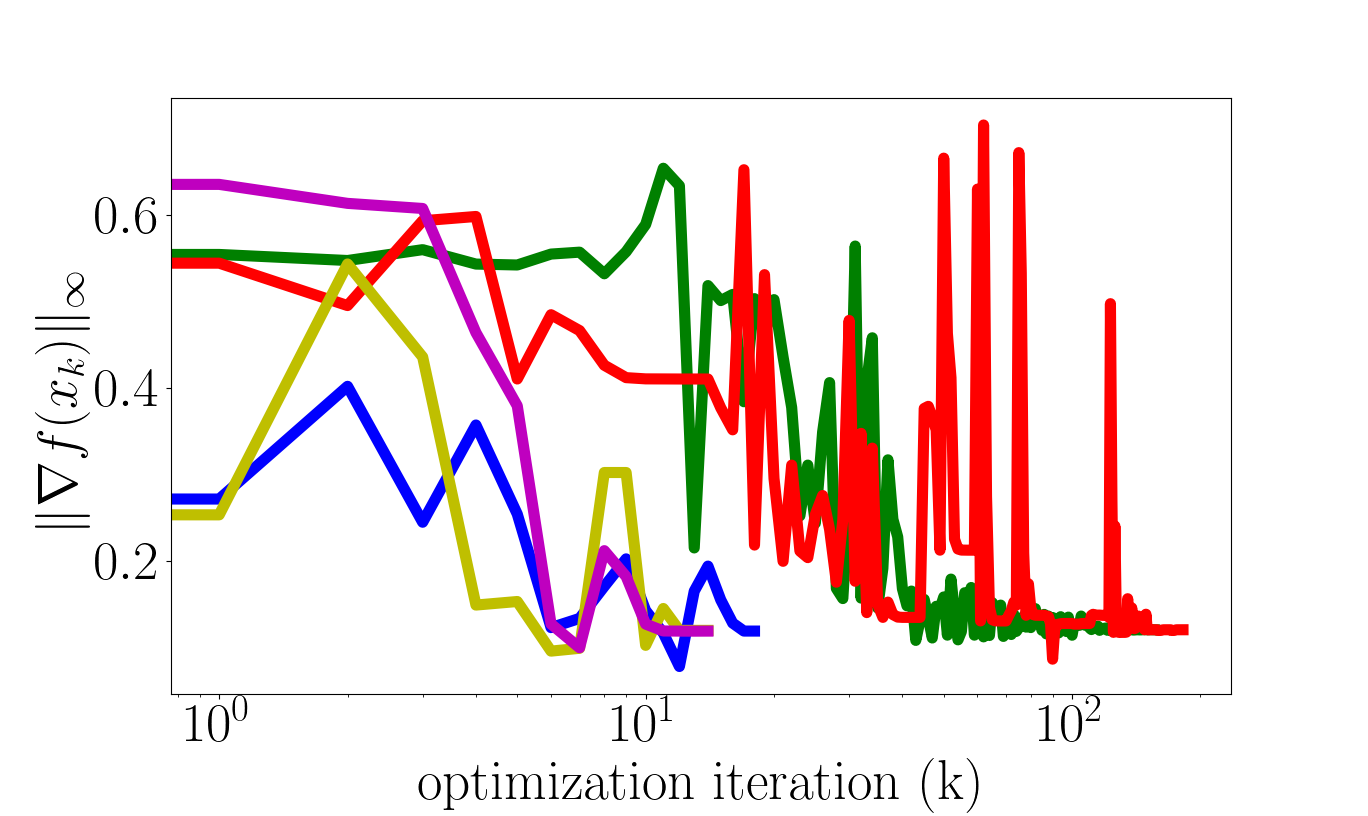}
\caption{4 control points}\label{norm-4pts-NACA0012}
    \end{subfigure}%
\begin{subfigure}{0.33\textwidth}
\includegraphics[width=\linewidth]{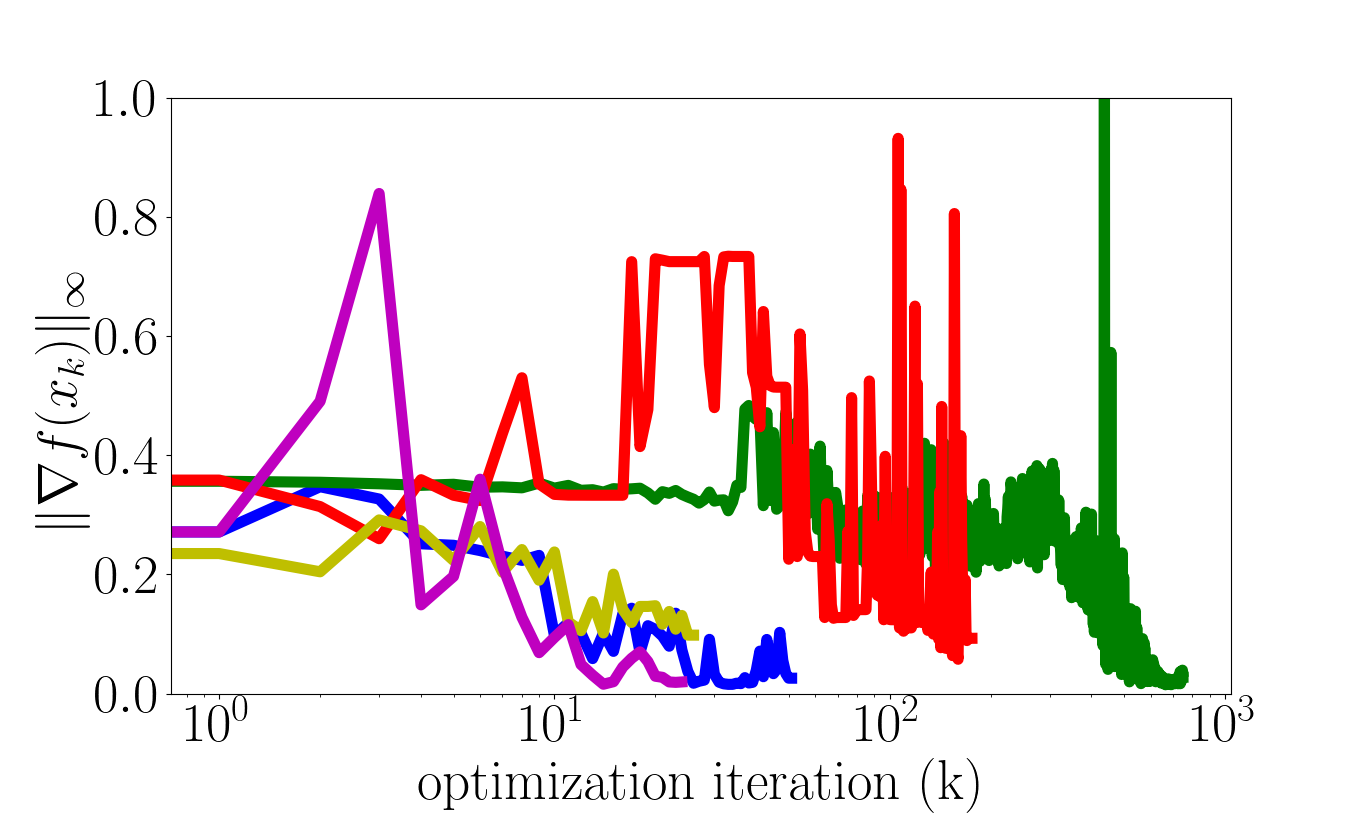} \caption{8 control points}\label{norm-8pts-NACA0012}
    \end{subfigure}%
    \begin{subfigure}{0.33\textwidth}
\includegraphics[width=\linewidth]{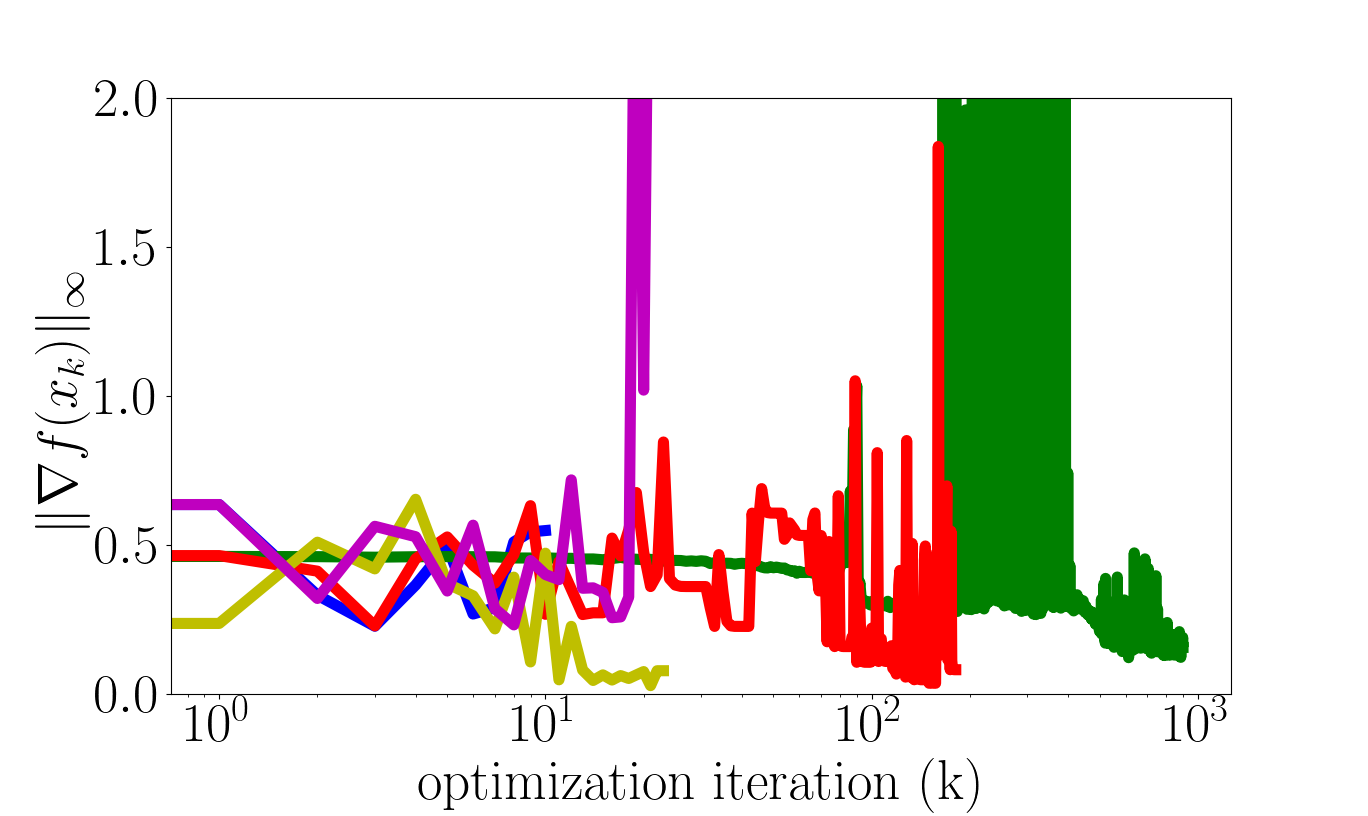} 
\caption{16 control points}\label{norm-16pts-NACA0012}
    \end{subfigure}
    \begin{subfigure}{\textwidth}
        \centering
        \includegraphics[width=0.9\linewidth]{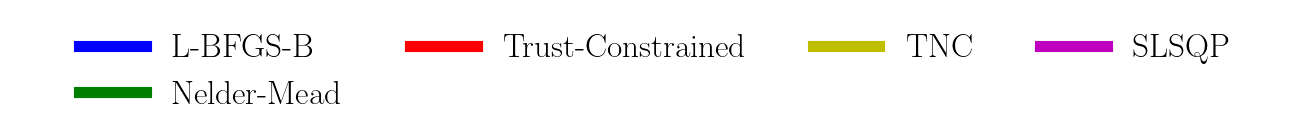} 
    \end{subfigure}
    \caption{{\bf Unconstrained NACA0012.} Gradient infinity norm history.}
    \label{norm-NACA0012}
\end{figure}


Finally, we also show the gradient infinity norm of the objective function in \Cref{norm-NACA0012}. Notice that we don't achieve strict optimality ($\|\nabla f(\x)\|_\infty = 0$) in any of the algorithms; instead, the optimizer terminates after a few attempts to find a stationary point. Typically, the off-the-shelf optimizers are designed to terminate upon encountering diminishing returns in the improvement of the objective function. Furthermore, note that the gradient norm can drastically vary with iterations (particularly in \tc), which is likely exacerbated due to incomplete convergence of the adjoint solver.  Therefore, we remark that in practice for aerodynamic design optimization problems, derivative-based optimizers don't always necessarily find a stationary point.
Note that in \Cref{norm-NACA0012}, we also include \nm, where we manually compute the gradient (by calling the adjoint solver) at every step for comparison. This reveals that the variability in gradient norm, for derivative-based approaches, can be on par with that of derivative-free optimizers which do not seek a stationary point.

\subsection{Unconstrained RAE2822}
\begin{figure}[htb!]
\begin{minipage}{0.47\textwidth}
    \includegraphics[width=\linewidth] {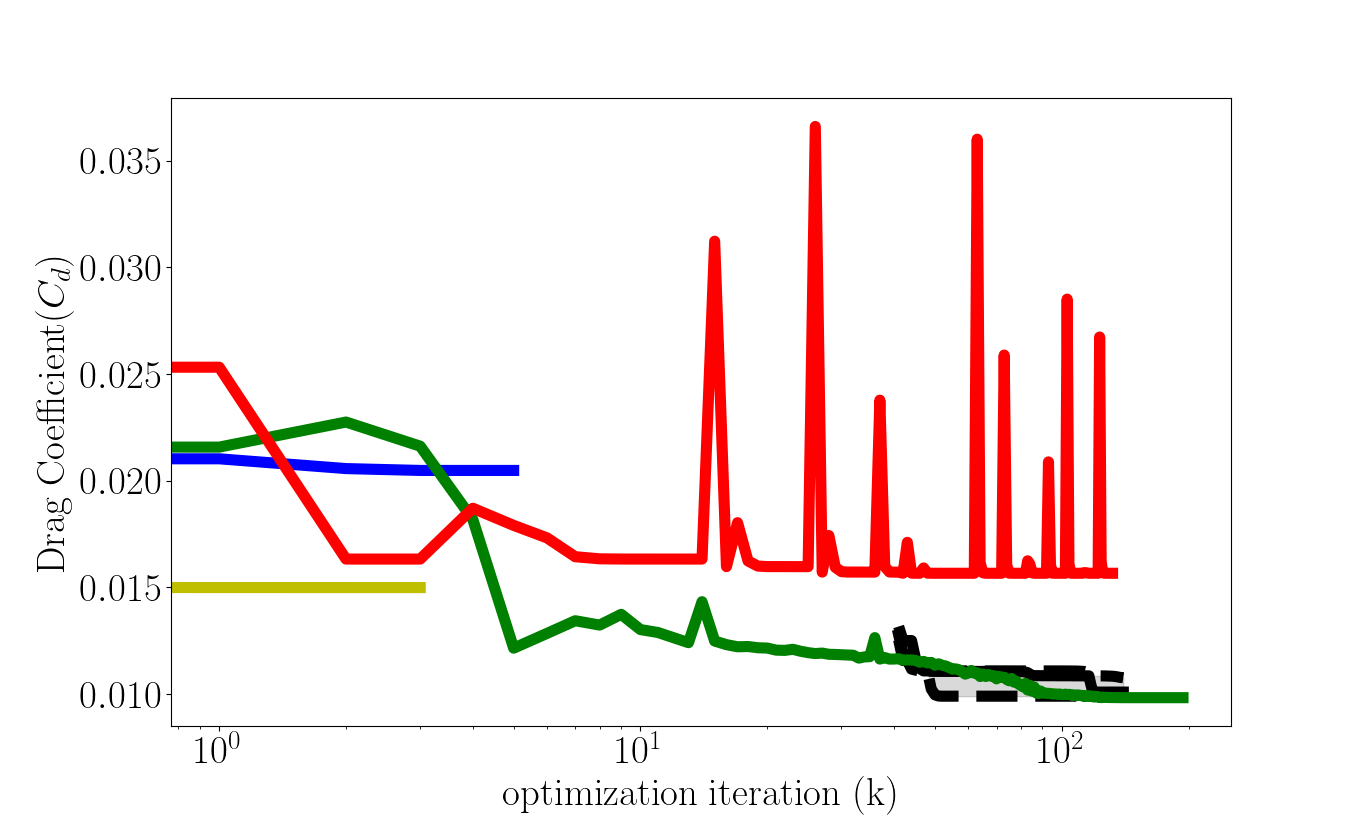} 
    \subcaption{4 control points}\label{Cd-4pts-RAE2822}
    \end{minipage}
    \hspace{\fill} 
    \begin{minipage}{0.47\textwidth}
    \includegraphics[width=\linewidth]{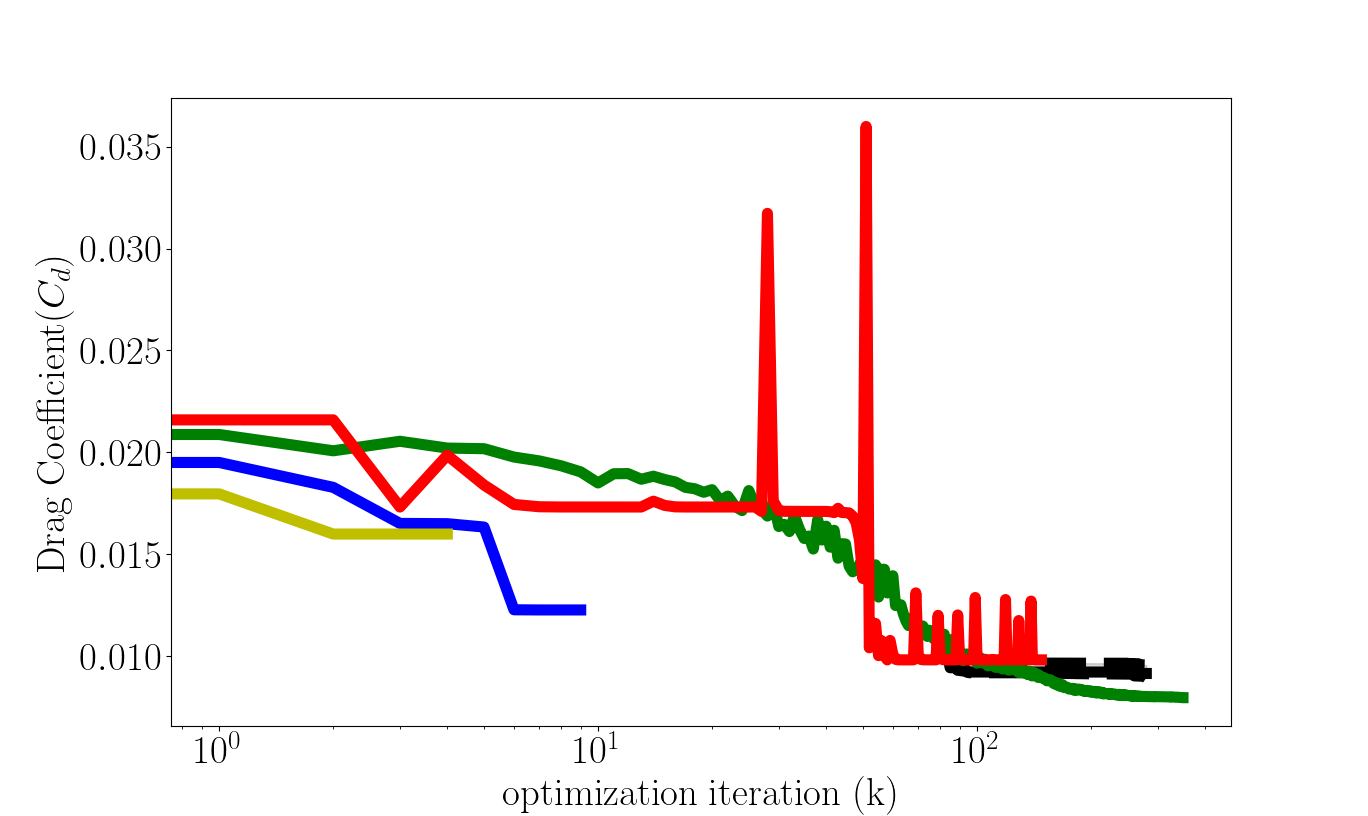} 
    \subcaption{8 control points}\label{Cd-8pts-RAE2822}
    \end{minipage}


    \begin{minipage}{0.47\textwidth}
    \includegraphics[width=\linewidth]{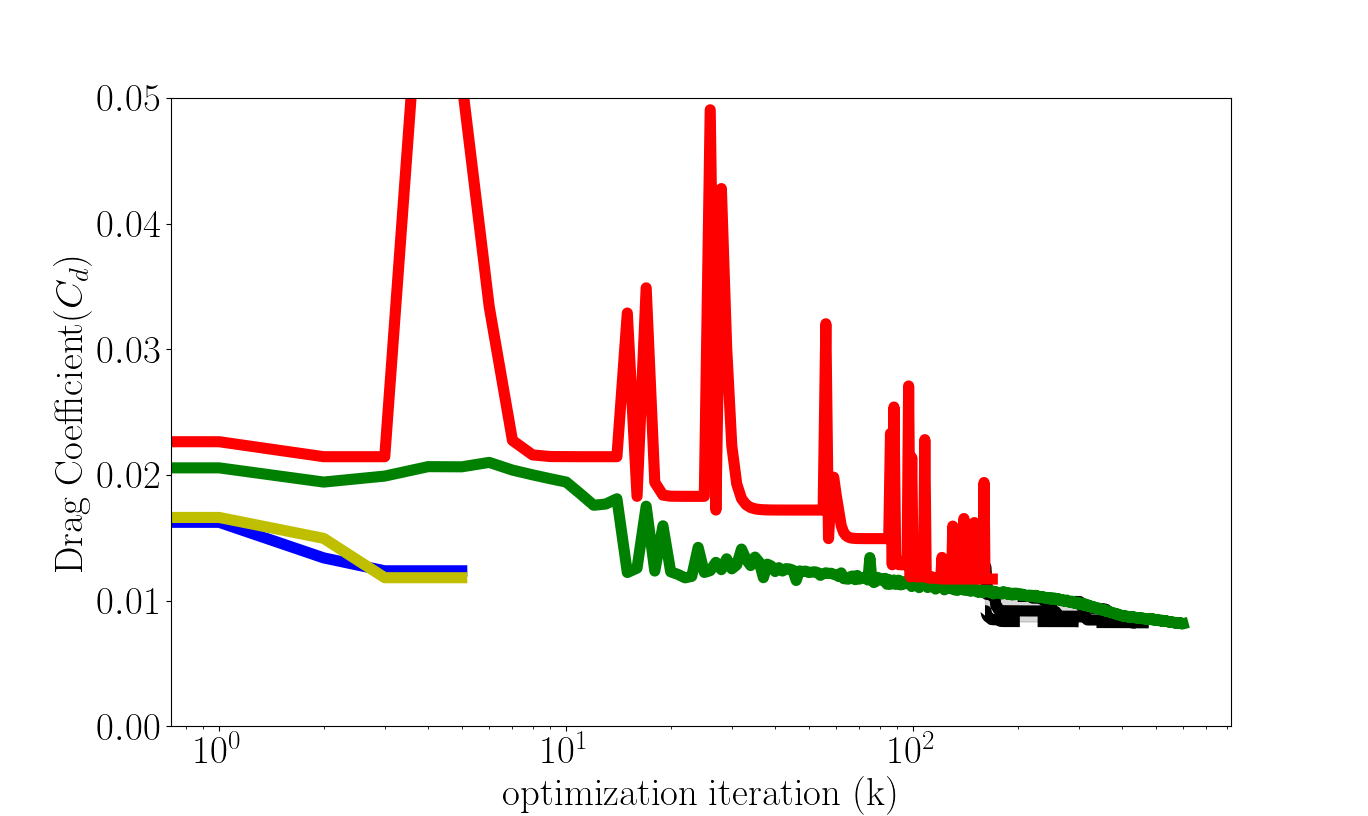} 
    \subcaption{16 control points}\label{Cd-16pts-RAE2822}
    \end{minipage}
    \hspace{\fill} 
    \begin{minipage}{0.47\textwidth}
    \includegraphics[width=\linewidth]{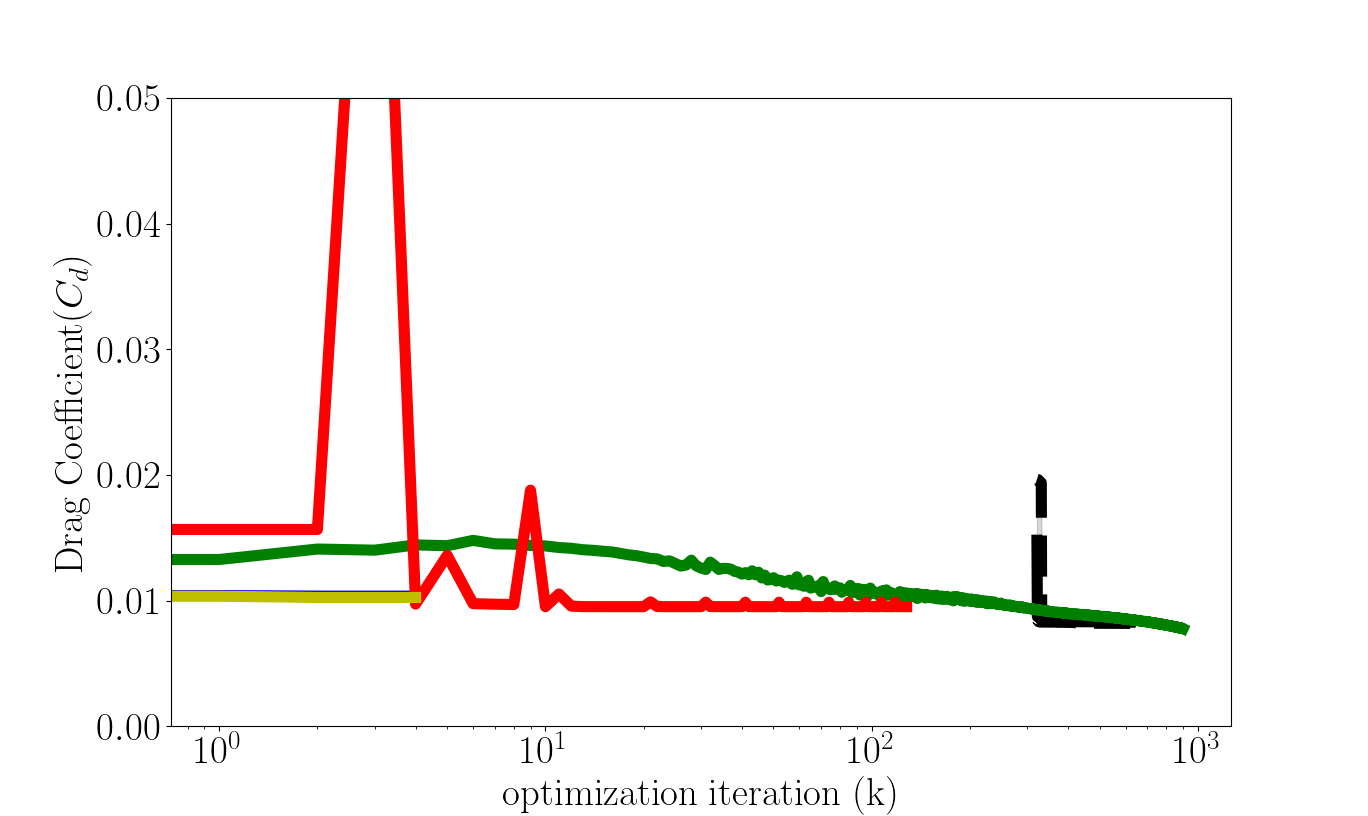} 
    \subcaption{32 control points}\label{Cd-32pts-RAE2822}
    \end{minipage}
    \begin{center}
     \begin{minipage}{\textwidth}
     \includegraphics[width=0.9\linewidth]{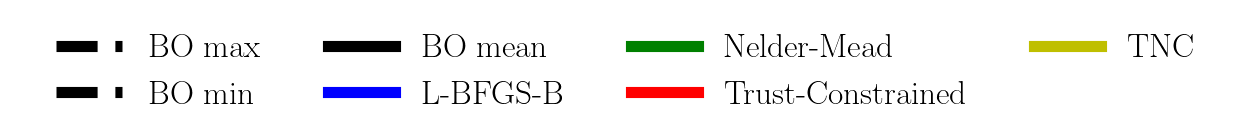} 
     \end{minipage}
     \end{center}
\caption{{\bf Unconstrained RAE2822.} Convergence histories versus optimization iterations.} \label{Cd-RAE2822}
\end{figure}

{\label{RAE2822}}
We now present the optimization results for RAE2822 under viscous conditions -- this leads to more expensive forward and adjoint solver costs. The control point displacements are contained to $[-0.002, 0.002]$. \change{Note that these bounds include an internal scaling factor ($\texttt{OPT\_RELAX\_FACTOR} = 10^2$) in SU2 -- thus $\pm 0.0002$ is equivalent to a $20\%$ chord length displacement. The scaling factor with a smaller bound helps improve the primal and adjoint solve convergence in SU2.} \Cref{Cd-RAE2822} shows the drag coefficient history for $4,8,16,$ and $32$ control points. Notice that derivative-free approaches (\nm~and~\bo) find the lowest objective value compared to all derivative-based approaches. \bfgs~and \tnc, terminated early without achieving much reduction in the objective function. And, as in the NACA0012 experiment, the \tc~objective history shows significantly more volatility than other methods. \changetwo{For this experiment, \slsqp~repeatedly terminated with linesearch failures or without satisfying our KKT-based convergence criterion despite extensive tuning of step sizes, tolerances, and basic scaling, and thus we classify it as having failed to converge within the allotted evaluation budget and omit it from the main performance plots.
failed to converge for this experiment and hence was not included in the results. We attribute this behavior to the
adjoint-based gradient
information being contaminated by solver tolerances and turbulence-model non-smoothness,
which can cause the Hessian approximation in SLSQP to become ill-conditioned, leading to
step-rejection and stalling.} \change{Further, the choice of the initial point for derivative-based algorithms plays a major role in their convergence.} On the other hand, despite the randomly chosen seed points for \bo, their variability in performance (visualized as shaded gray in the figure) is quite acceptable for practical applications. 
\begin{figure}[htb!]
\begin{minipage}{0.47\textwidth}
    \includegraphics[width=\linewidth] {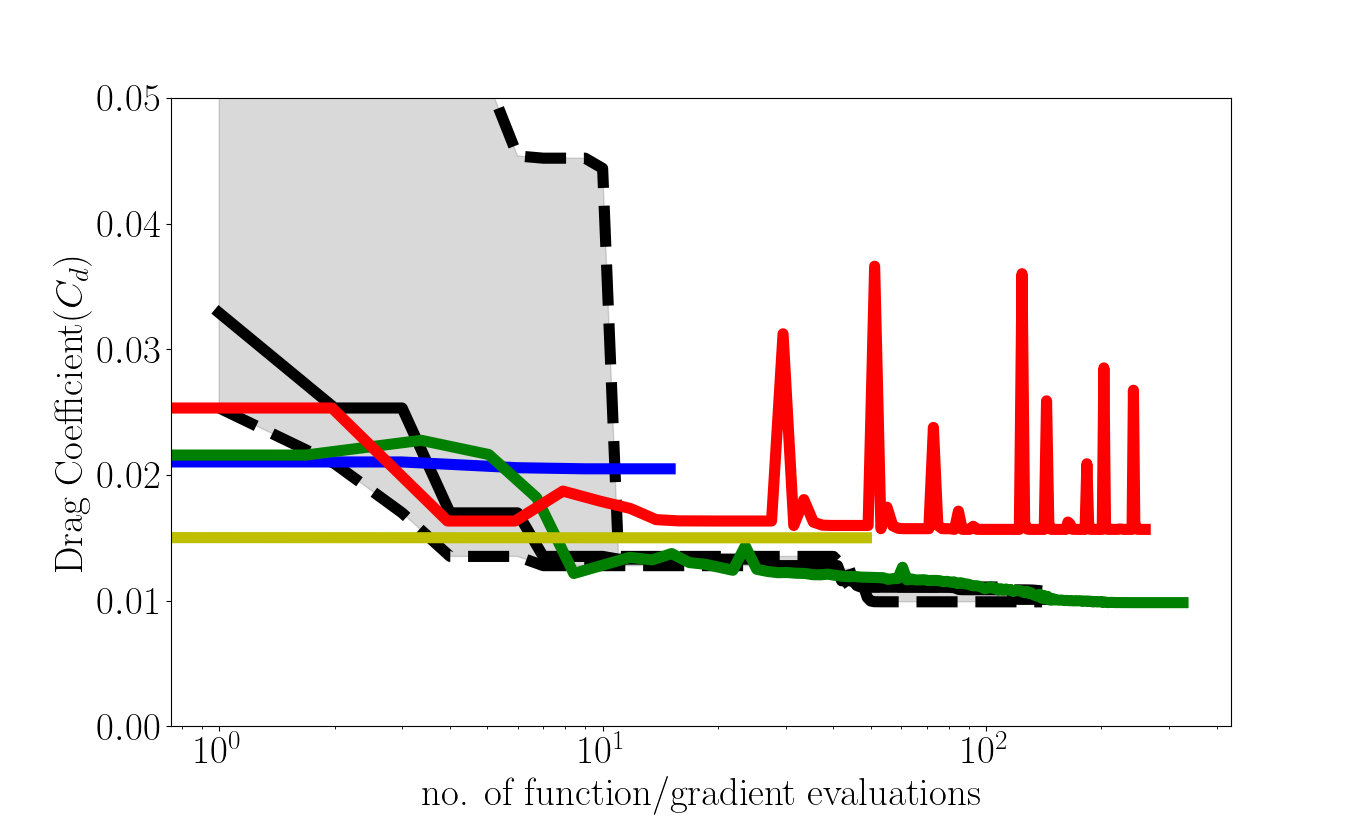}
    \subcaption{4 control points}\label{feval-4pts-RAE2822}
    \end{minipage}
    \hspace{\fill} 
    \begin{minipage}{0.47\textwidth}
    \includegraphics[width=\linewidth]{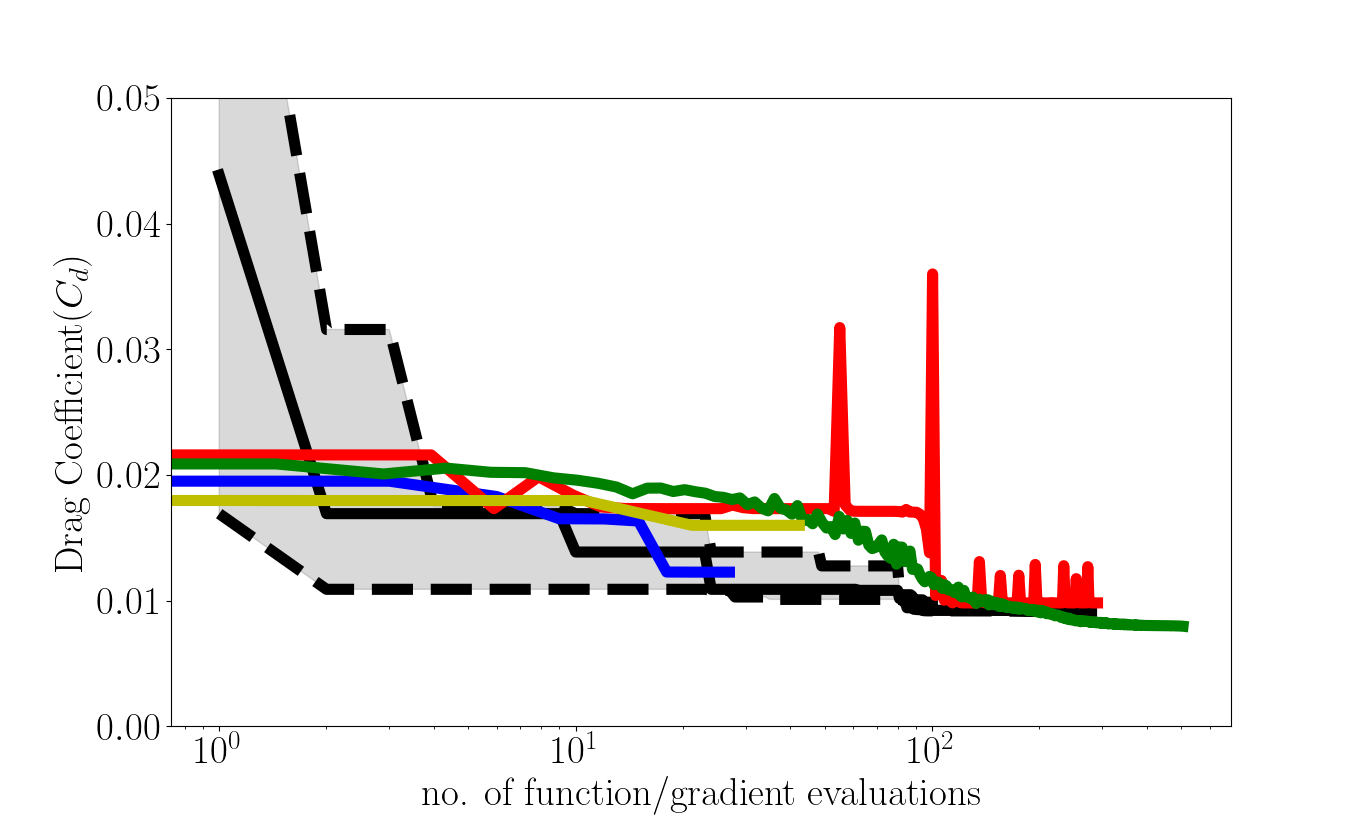} 
    \subcaption{8 control points}\label{feval-8pts-RAE2822}
    \end{minipage}


    \begin{minipage}{0.47\textwidth}
    \includegraphics[width=\linewidth]{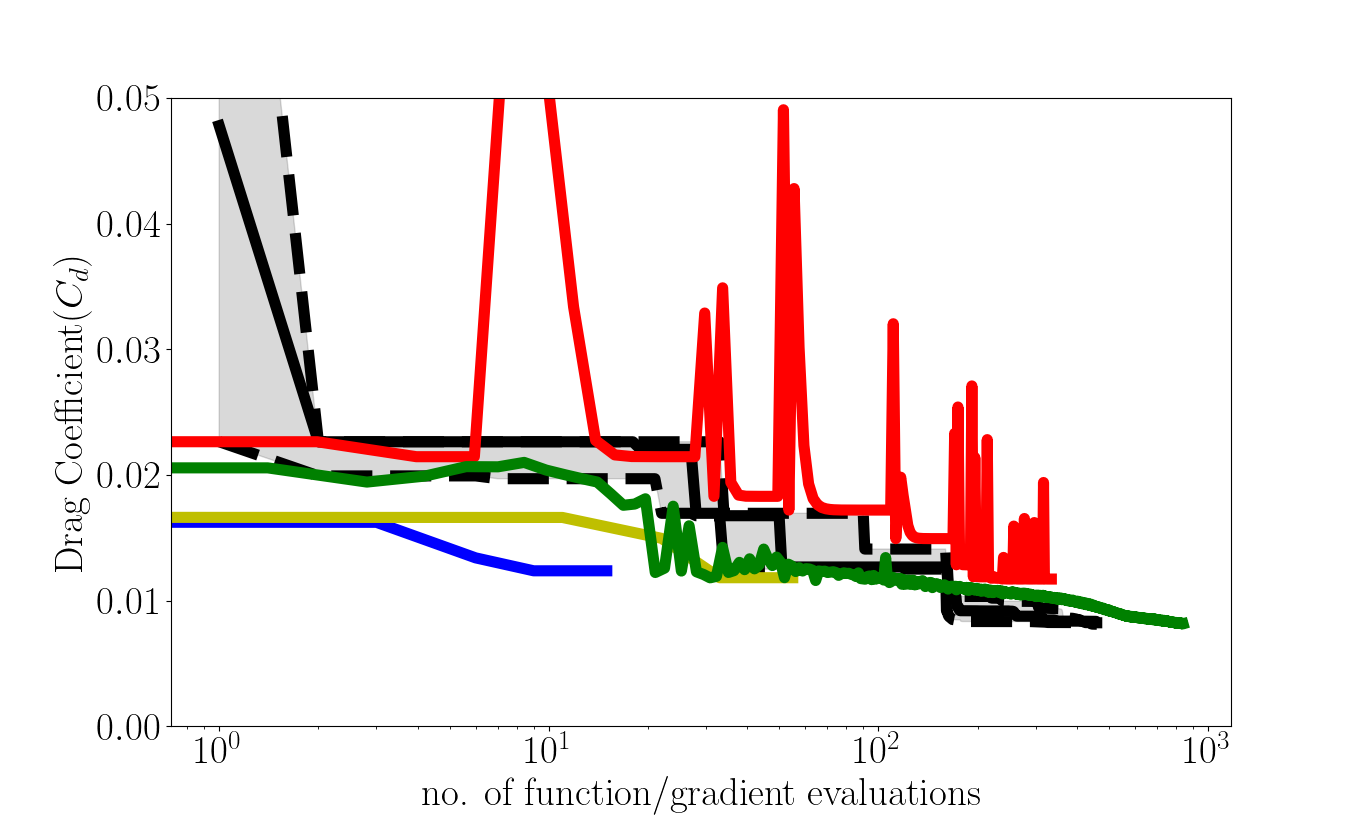} 
    \subcaption{16 control points}\label{feval-16pts-RAE2822}
    \end{minipage}
    \hspace{\fill} 
    \begin{minipage}{0.47\textwidth}
    \includegraphics[width=\linewidth]{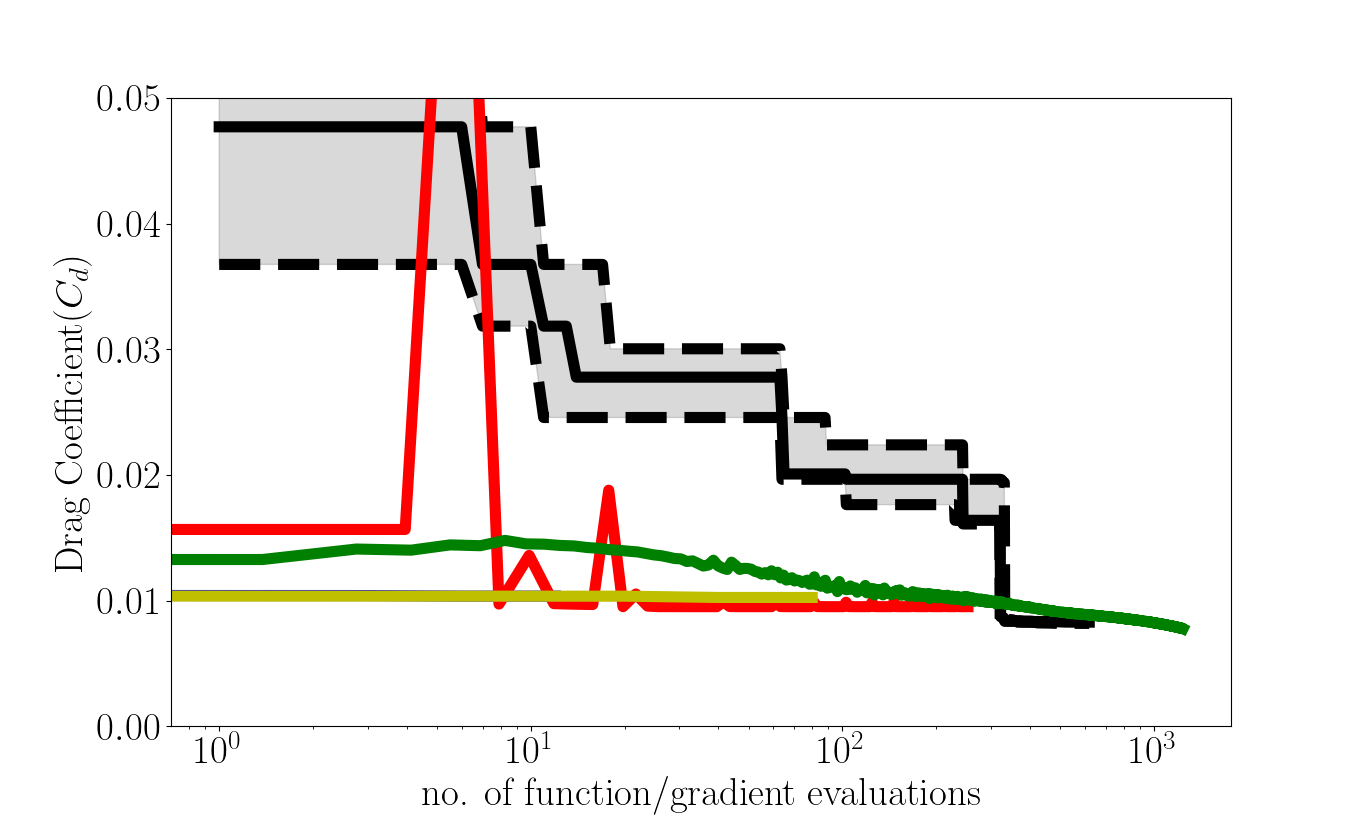} 
    \subcaption{32 control points}\label{feval-32pts-RAE2822}
    \end{minipage}
    \begin{center}
     \begin{minipage}{\textwidth}
     \includegraphics[width=0.9\linewidth]{Legends_RAE_Unconstrained.png} 
     \end{minipage}
     \end{center}
\caption{{\bf Unconstrained RAE2822.} Convergence histories versus number of function/gradient evaluations.} \label{feval-RAE2822}
\end{figure}
\Cref{feval-RAE2822} shows the convergence history versus the number of function evaluations. Again, notice that the benefit of derivative-free methods is best demonstrated in this context.
The derivative-based methods prematurely terminate with higher objective values  despite utilizing expensive gradient and (quasi) $2nd$ order information. Finally, the gradient norm plots in \Cref{norm-RAE2822} show that strict optimality is not quite achieved in any of the algorithms, suggesting that the domain might not contain a stationary point or, potentially, the objective function is non-smooth. To verify this, we tested the same experiment with a different parametrization, namely the Hicks-Henne bump functions \cite{Hicks1978}, shown in \Cref{hicks-henne-RAE2822}. We observe that the gradient norm still plateaus at $10^{-1}$, substantiating our claim.

\changethree{In the unconstrained RAE2822 (RANS) case, SciPy's SLSQP frequently terminated due to line-search failure or reaching iteration limits before satisfying the stopping criteria used in this study; therefore, we omit its trajectory from the main plots for that case and report its termination status separately. We interpret this as an empirical observation for this solver configuration and the gradient fidelity available in our pipeline, rather than as a general statement about SQP methods.}

\begin{figure}[htb!]
\begin{minipage}{0.47\textwidth}
    \includegraphics[width=\linewidth]{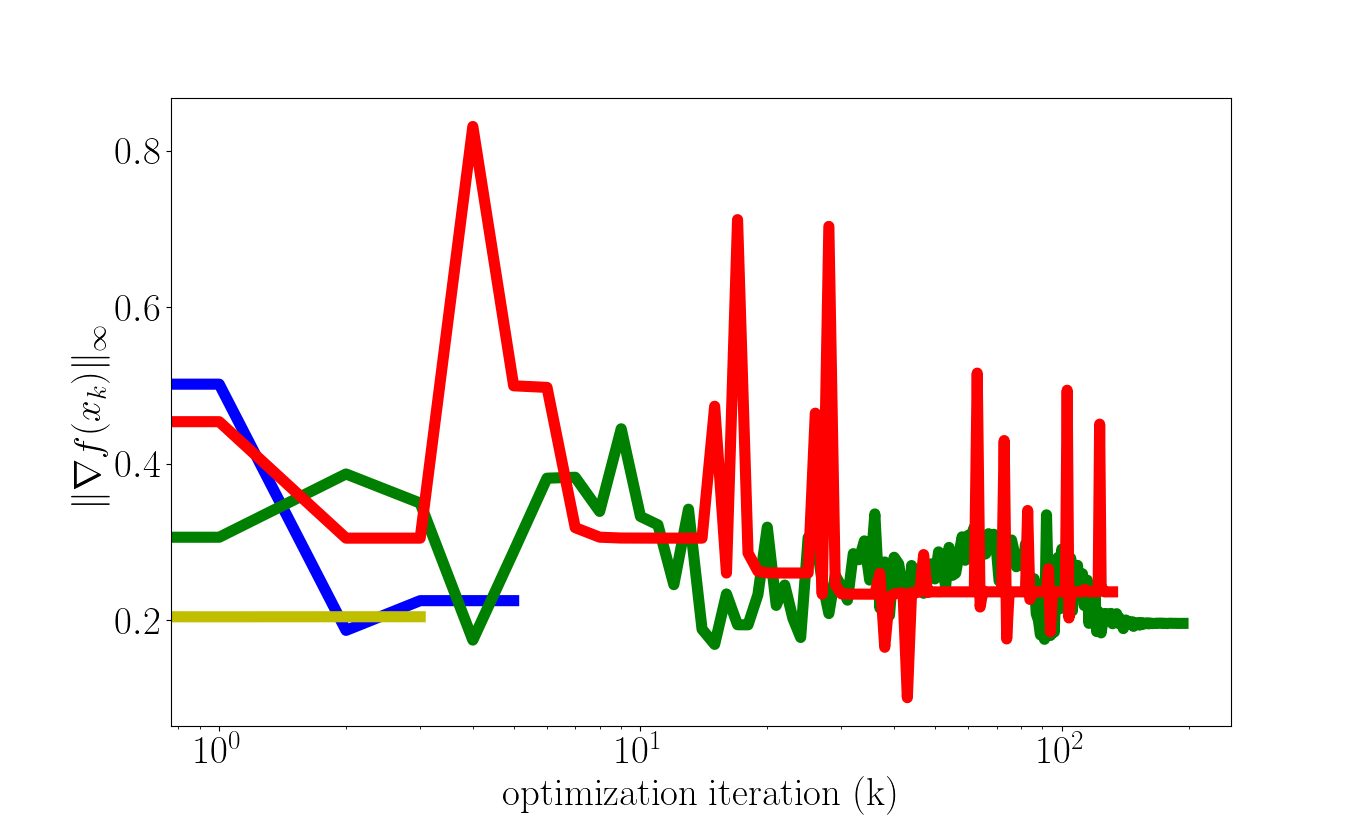} 
    \subcaption{4 control points}\label{norm-4pts-RAE2822}
    \end{minipage}
    \hspace{\fill} 
    \begin{minipage}{0.47\textwidth}
    \includegraphics[width=\linewidth]{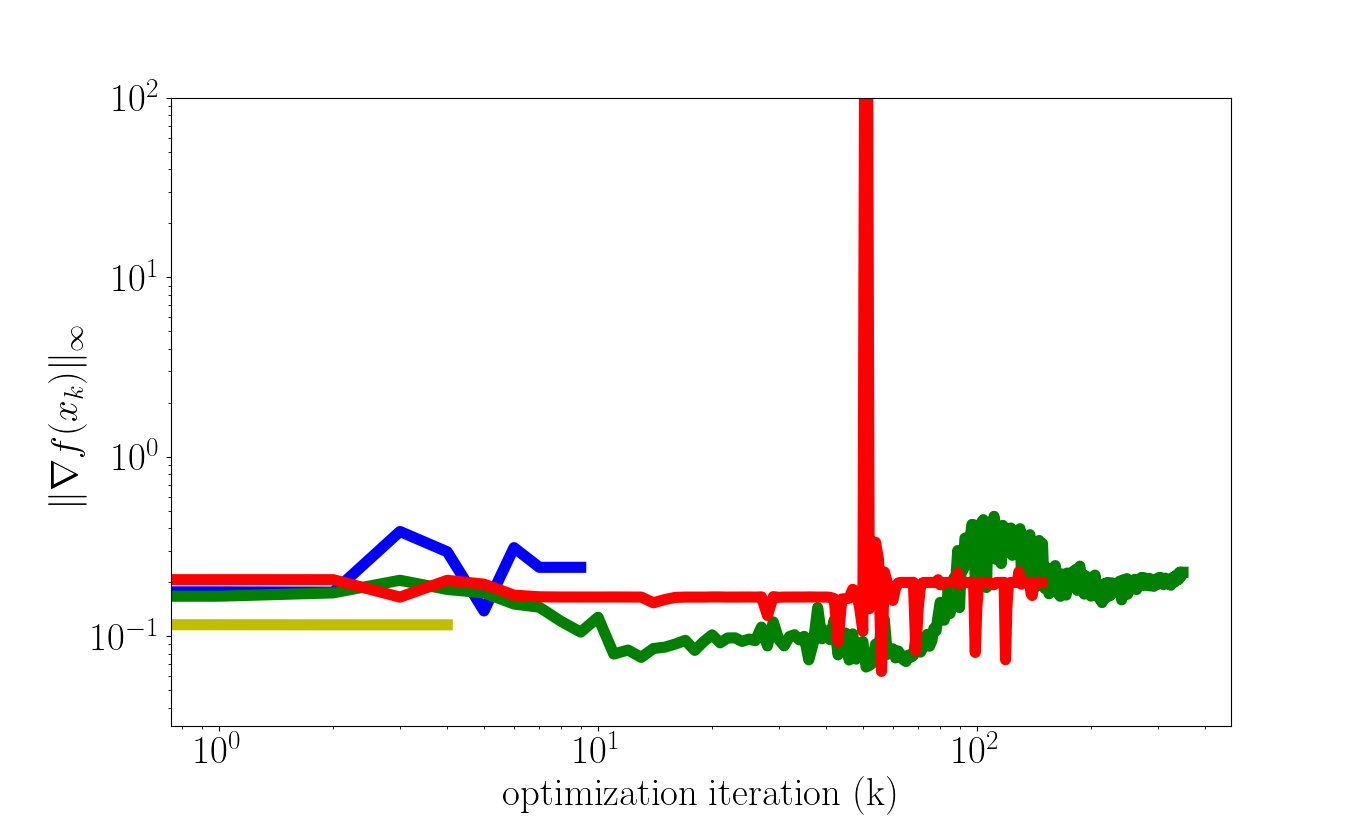} 
    \subcaption{8 control points}\label{norm-8pts-RAE2822}
    \end{minipage}


    \begin{minipage}{0.47\textwidth}
    \includegraphics[width=\linewidth]{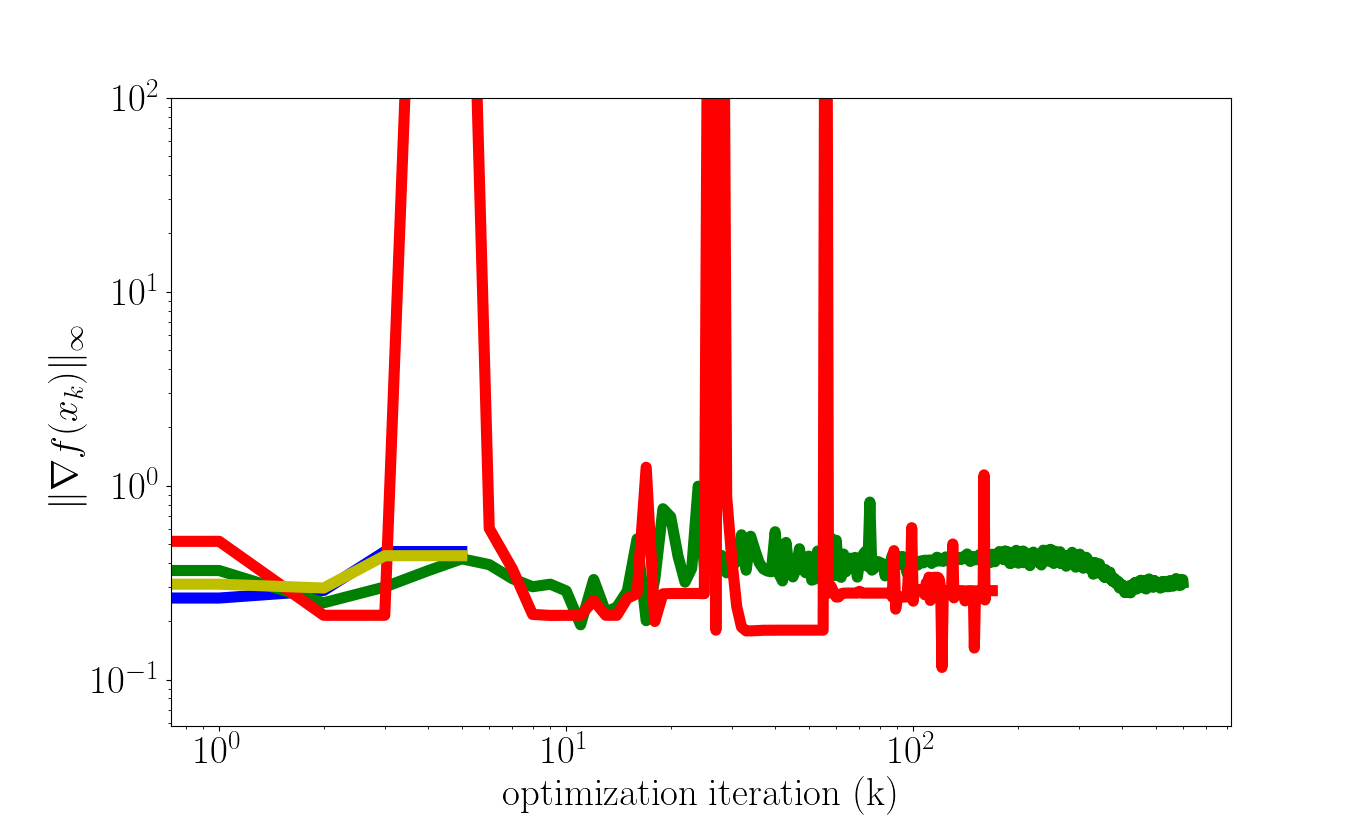} 
    \subcaption{16 control points}\label{norm-16pts-RAE2822}
    \end{minipage}
    \hspace{\fill} 
    \begin{minipage}{0.47\textwidth}
    \includegraphics[width=\linewidth]{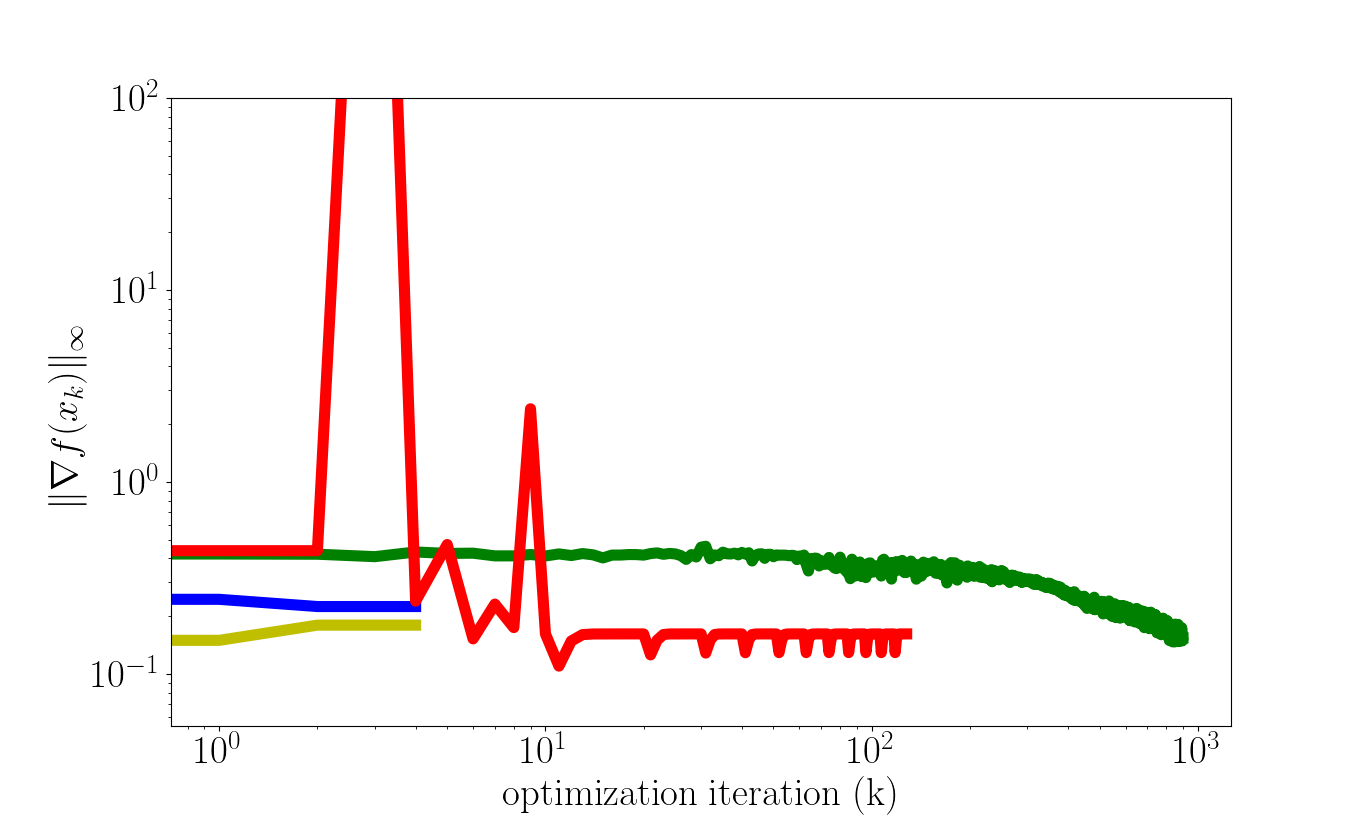} 
    \subcaption{32 control points}\label{norm-32pts-RAE2822}
    \end{minipage}
    \begin{center}
     \begin{minipage}{\textwidth}
     \includegraphics[width=0.9\linewidth]{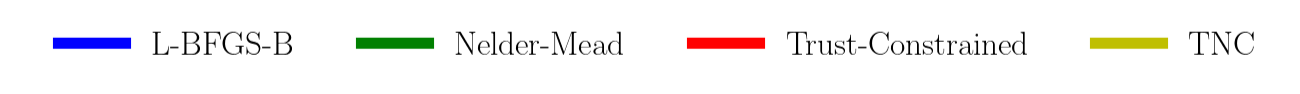} 
     \end{minipage}
     \end{center}
\caption{{\bf Unconstrained RAE2822.} Gradient infinity norm history.} \label{norm-RAE2822}
\end{figure}

\begin{figure}[htb!]
\begin{minipage}{0.47\textwidth}
    \includegraphics[width=\linewidth]{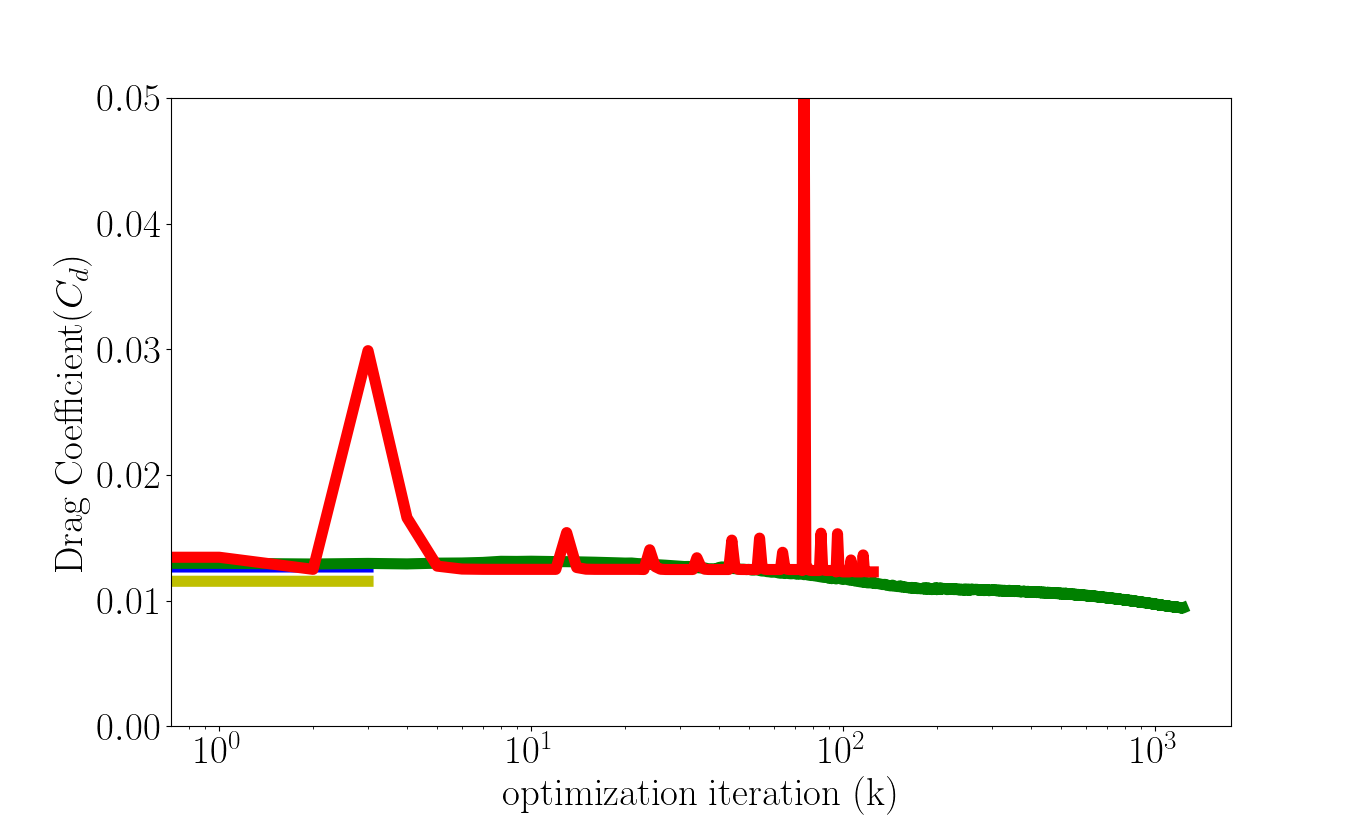} 
    \subcaption{$C_d$ vs iteration}\label{Cd_hh-38pts-RAE2822}
    \end{minipage}
    \hspace{\fill} 
    \begin{minipage}{0.47\textwidth}
    \includegraphics[width=\linewidth]{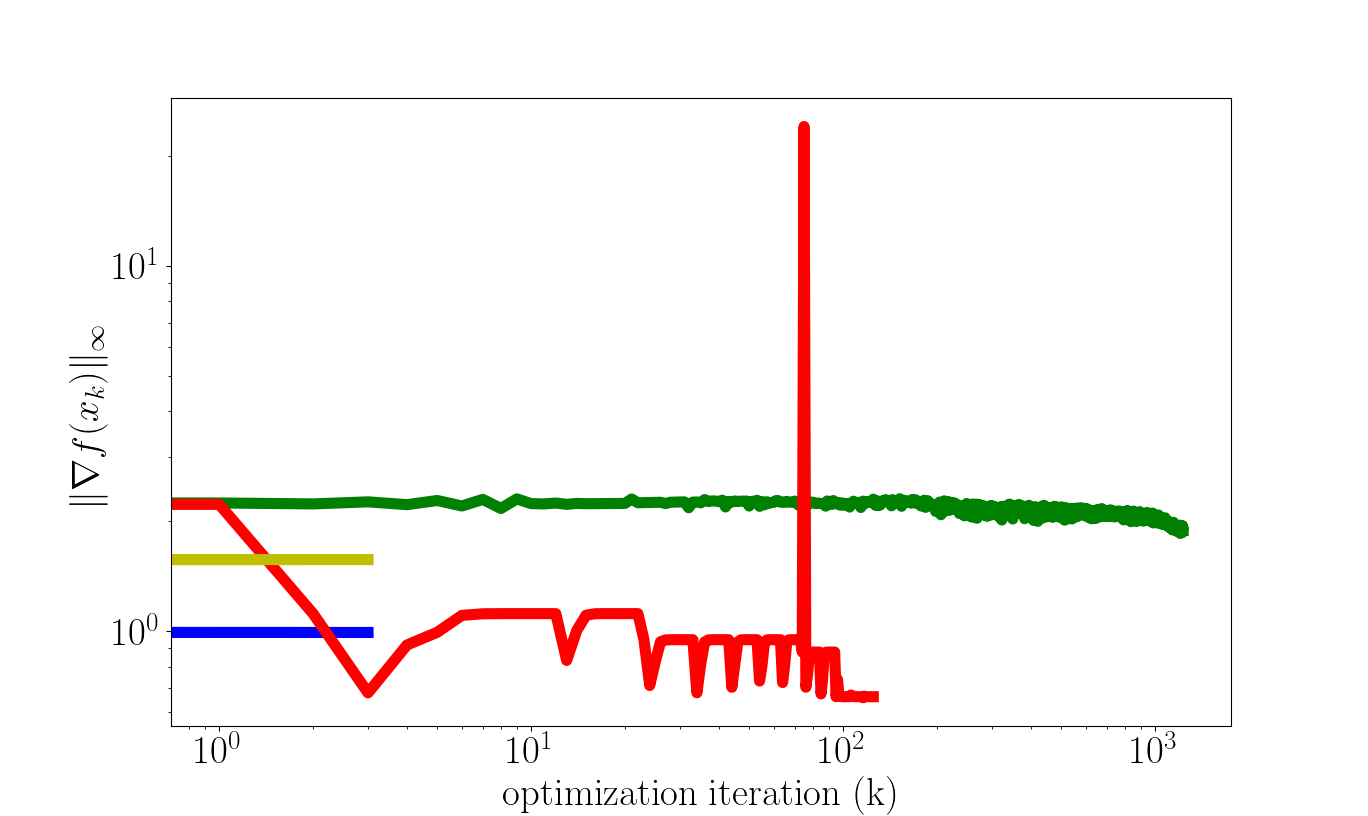} 
    \subcaption{infinity norm of gradient of $C_d$ vs. iteration}\label{norm_hh-38pts-RAE2822}
    \end{minipage}
    \begin{center}
     \begin{minipage}{\textwidth}
     \includegraphics[width=0.9\linewidth]{legends_hicks_henne.png} 
     \end{minipage}
     \end{center}
\caption{{\bf Unconstrained RAE2822.} Convergence histories Hicks and Henne parametrization.} \label{hicks-henne-RAE2822}
\end{figure}

\subsection{Unconstrained ONERAM6}
\begin{figure}[htb!]
    \centering
    \begin{subfigure}{0.33\textwidth}
\includegraphics[width=\linewidth]{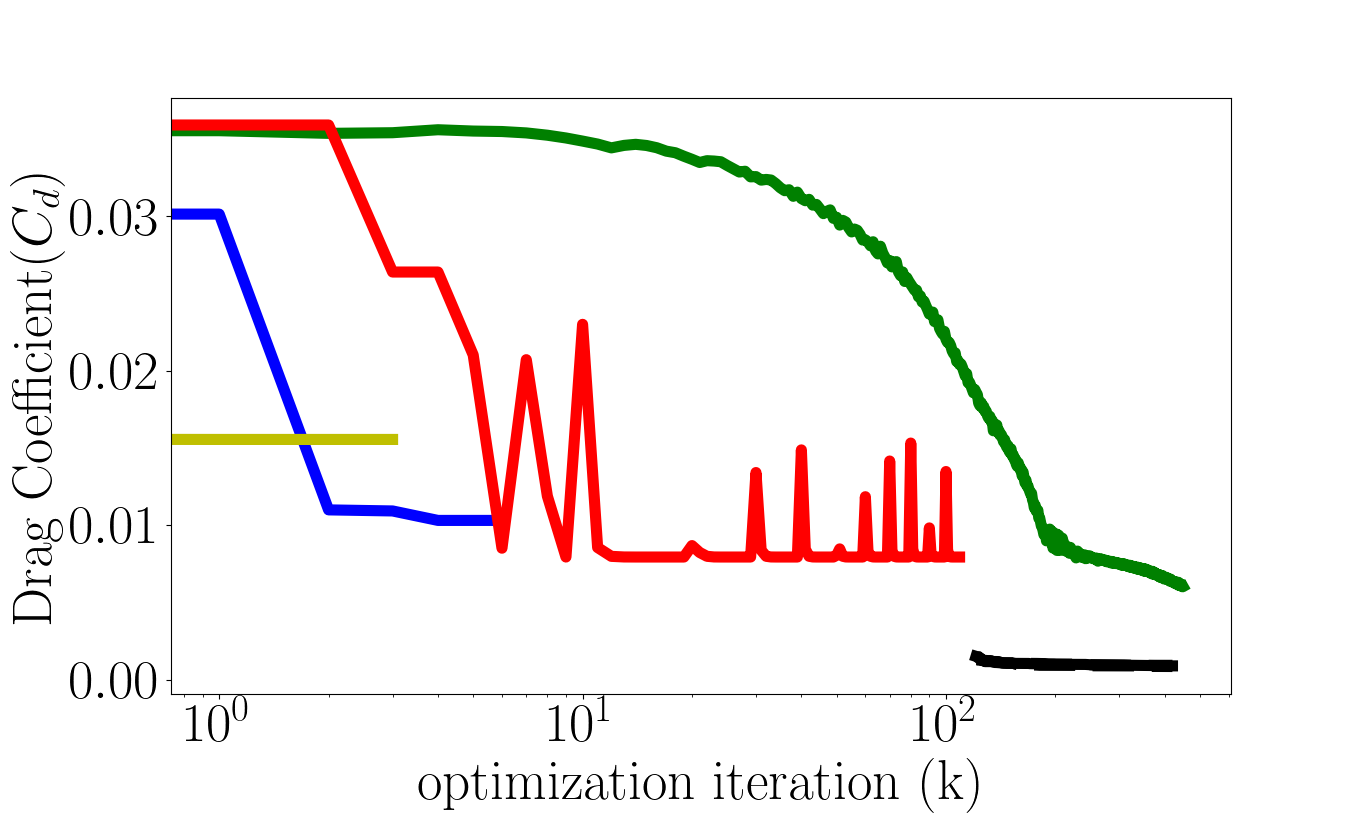}
\caption{$C_D$ vs iteration}\label{Cd-12pts-ONERA}
    \end{subfigure}%
\begin{subfigure}{0.33\textwidth}
\includegraphics[width=\linewidth]{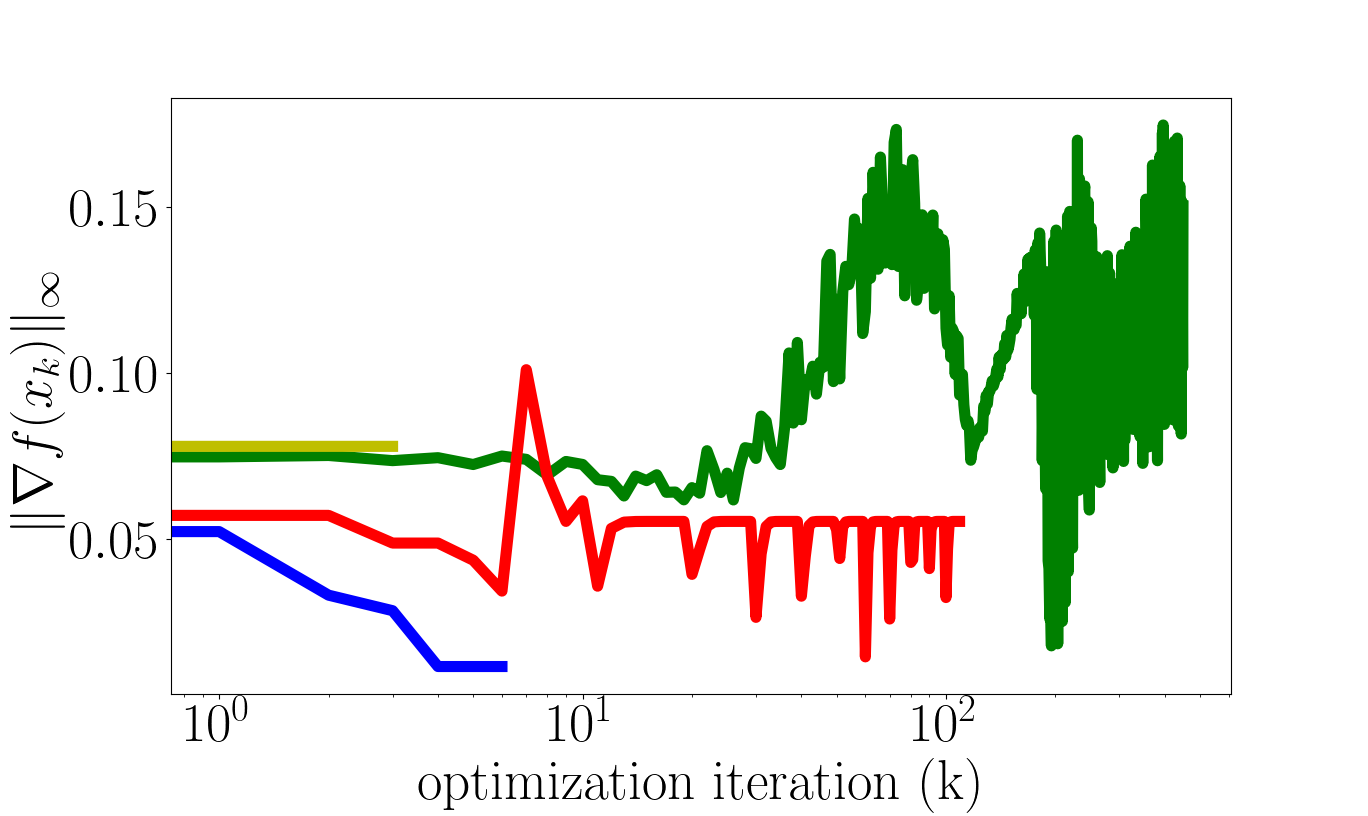} \caption{$\|\nabla f\|_\infty$ vs. iteration}\label{norm-12pts-ONERA}
    \end{subfigure}%
    \begin{subfigure}{0.33\textwidth}
\includegraphics[width=\linewidth]{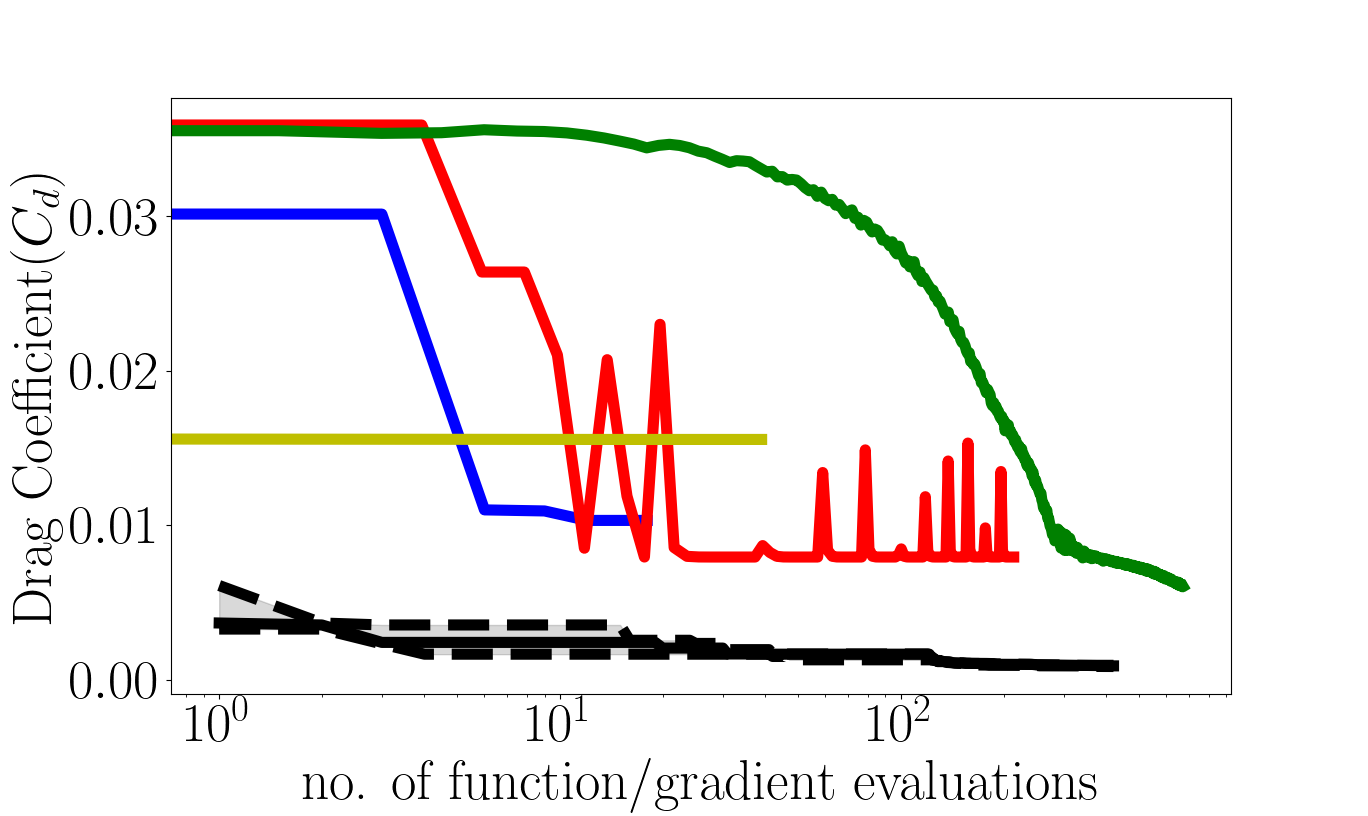} 
\caption{$C_D$ vs no. of function evals}\label{feval-12pts-ONERA}
    \end{subfigure}
    \begin{subfigure}{\textwidth}
        \centering
        \includegraphics[width=0.9\linewidth]{Legends_RAE_Unconstrained.png} 
    \end{subfigure}
    \caption{{\bf Unconstrained ONERAM6.} Convergence histories versus number of function/gradient evaluations.}
    \label{12pts-ONERA}
\end{figure}

For the ONERAM6 test case, the approach is repeated, but considering the computational complexity and simulation time, we demonstrate the results for only one case considering 12 control points. The domain is set as $\mcl{X} = [-0.00025, 0.00025]^{12}$. \Cref{12pts-ONERA} shows results for $M=0.8395$ and $\alpha=3.06$ degrees under inviscid conditions. The control points' placement and bounds are designed to ensure smooth wing deformation, avoiding any issues in the direct and adjoint simulations.   

The convergence history of $C_D$ with iteration is presented in \Cref{Cd-12pts-ONERA}. The trends observed for the derivative-based methods \bfgs~, \tc~, and \tnc~, are consistent with the previous two test cases discussed in \Cref{NACA0012} and \Cref{RAE2822}. \changetwo{In this constrained case as well, \slsqp~did not satisfy our KKT-based convergence tolerance within the evaluation budget, terminating instead due to linesearch failure or iteration limits despite variations of its hyperparameters; consequently, we do not include its trajectories in the main convergence plots.} However, notice that \bo~finds the lowest objective value compared to all other methods. In terms of the number of function evaluations (right of \Cref{feval-12pts-ONERA}), \bo~stands out by demonstrating the fastest convergence.

\subsection{Constrained RAE2822}{\label{RAE-constrained}}
The optimization problem we considered for this experiment can be formulated as described in Equation \ref{eq:opt_problem}.
\begin{equation} \label{eq:opt_problem}
\begin{split}
\min_{\mathbf{x} \in \mcl{X}} \quad & C_d \\
\text{s.t.} \quad & C_{mz} \leq 0.092, \\
& C_{l} = 0.75, \\
& t_{\text{max}} \geq 0.12,
\end{split}
\end{equation}
where $C_{mz}$ denotes the pitching moment coefficient of the airfoil about the quarter-chord, $C_l$ denotes the lift coefficient, and $t_{\text{max}}$ denotes the maximum airfoil thickness-to-chord ratio.   

\begin{figure}[htb!]
    \centering
    \begin{subfigure}{0.33\textwidth}
\includegraphics[width=\linewidth]{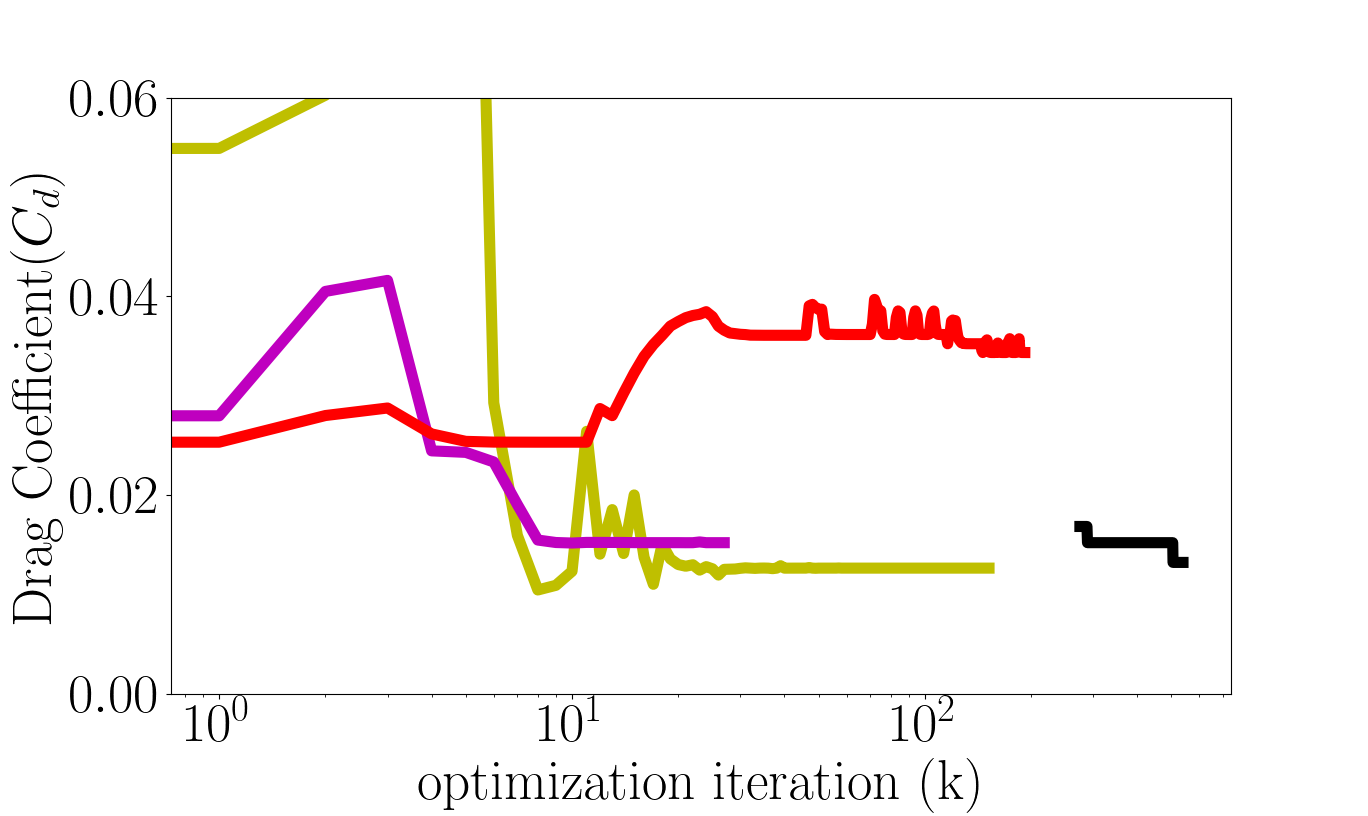}
\caption{4 control points}\label{Cd-4pts-RAE-Constrained}
    \end{subfigure}%
\begin{subfigure}{0.33\textwidth}
\includegraphics[width=\linewidth]{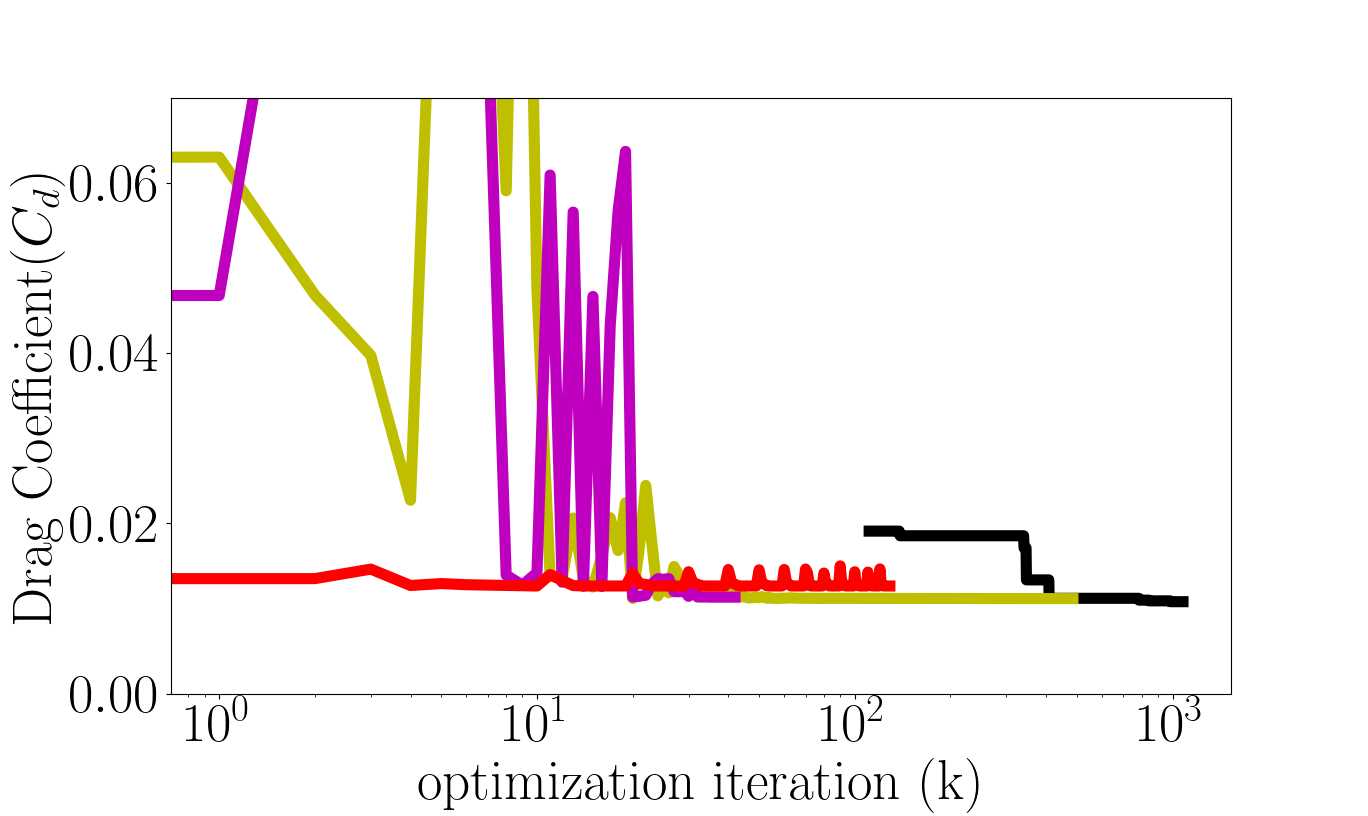} \caption{8 control points}\label{Cd-8pts-RAE-Constrained}
    \end{subfigure}%
    \begin{subfigure}{0.33\textwidth}
\includegraphics[width=\linewidth]{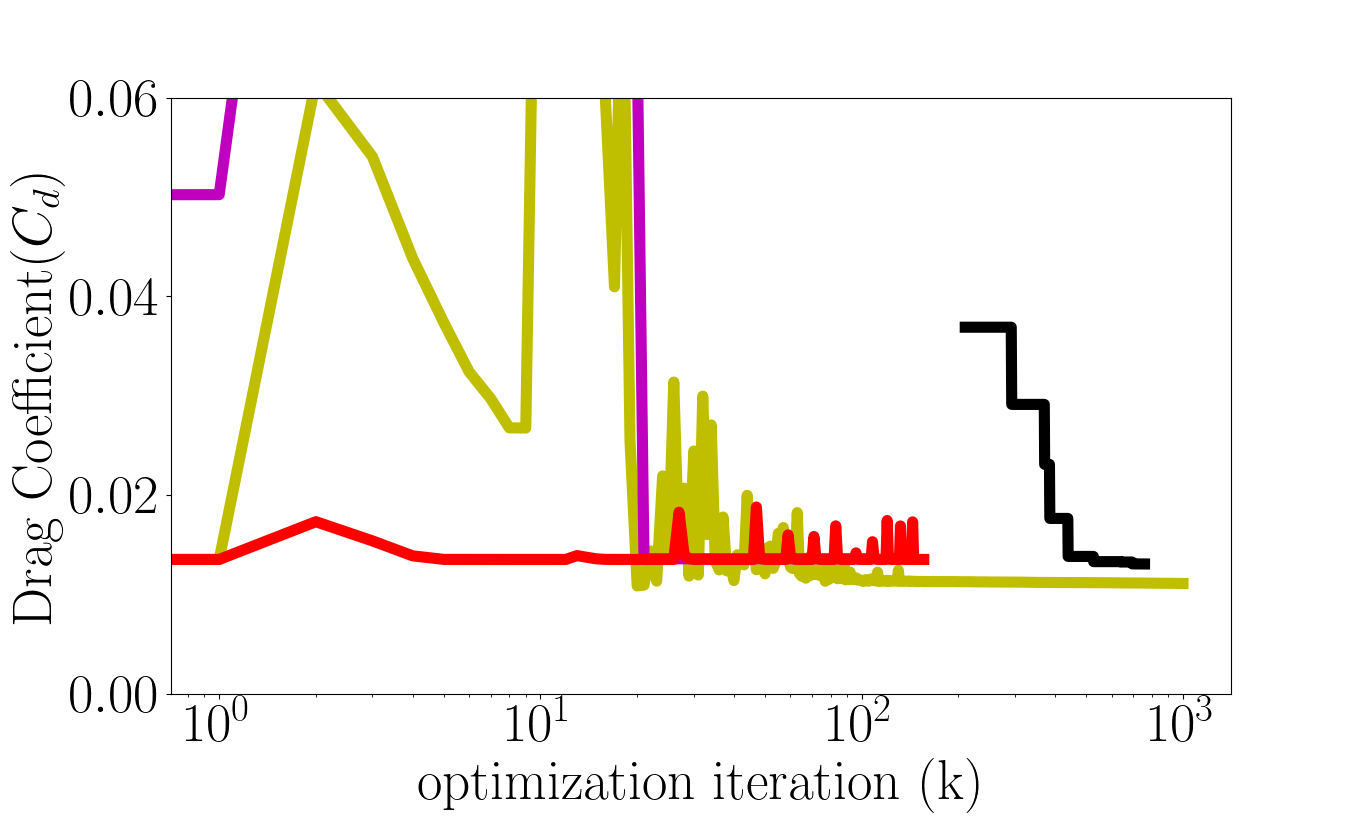} 
\caption{16 control points}\label{Cd-16pts-RAE-Constrained}
    \end{subfigure}
    \begin{subfigure}{\textwidth}
        \centering
        \includegraphics[width=0.9\linewidth]{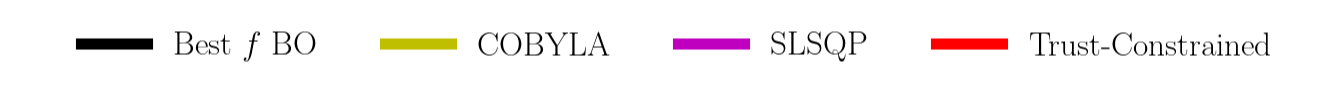} 
    \end{subfigure}
    \caption{{\bf Constrained RAE2822.} Convergence histories versus optimization iterations. \change{Black are best feasible objective values.}}
    \label{Cd-RAE-constrained}
\end{figure}

\begin{figure}[htb!]
    \centering
    \begin{subfigure}{0.33\textwidth}
\includegraphics[width=\linewidth]{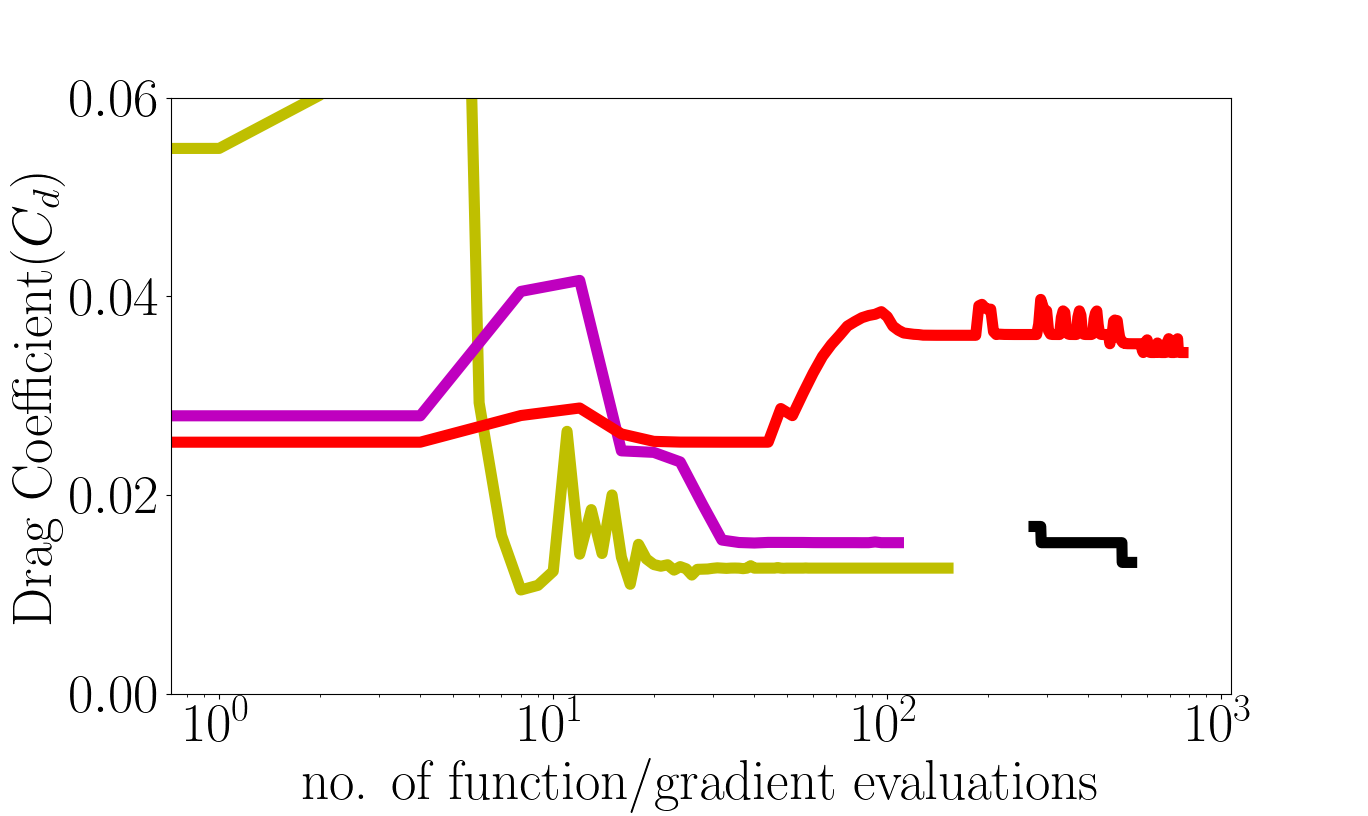}
\caption{4 control points}\label{feval-4pts-RAE-Constrained}
    \end{subfigure}%
\begin{subfigure}{0.33\textwidth}
\includegraphics[width=\linewidth]{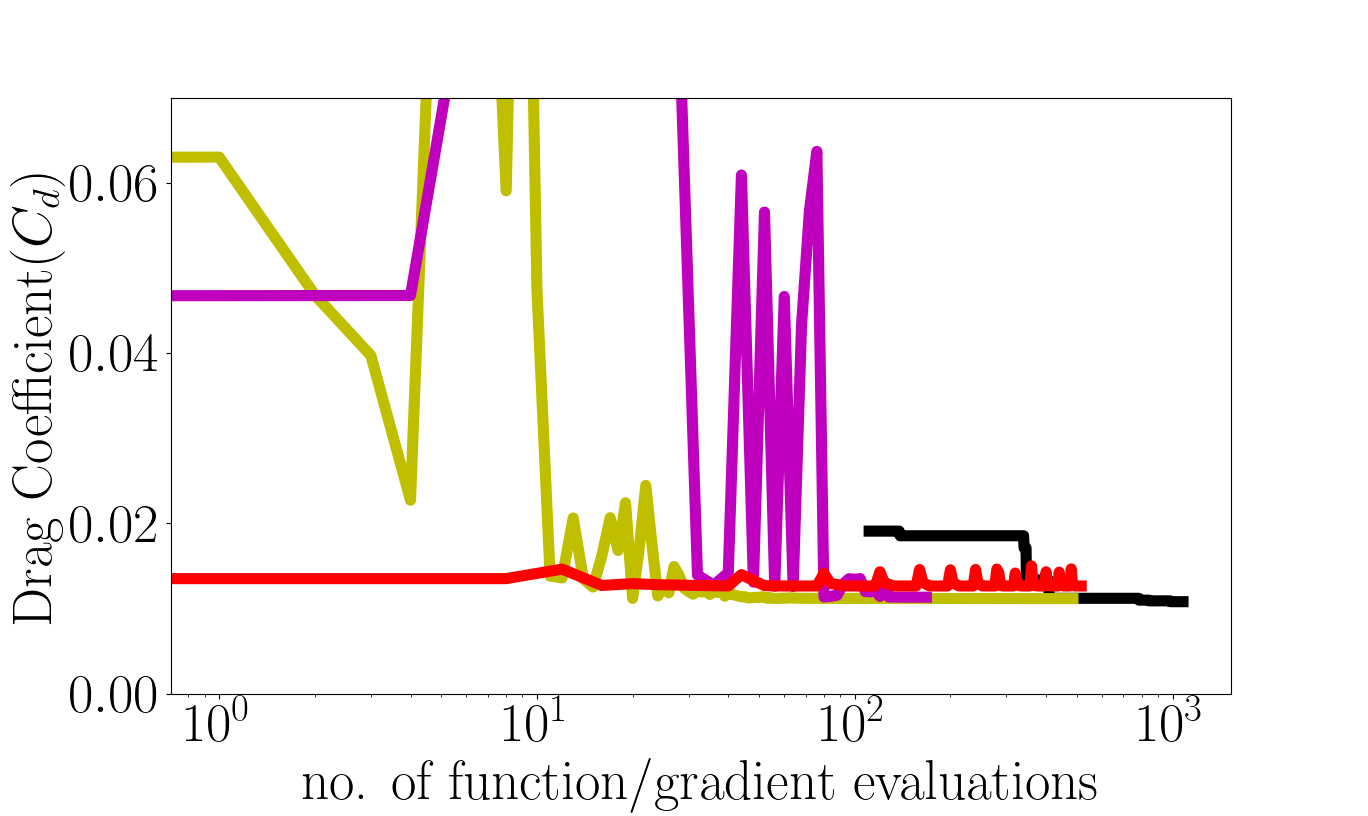} \caption{8 control points}\label{feval-8pts-RAE-Constrained}
    \end{subfigure}%
    \begin{subfigure}{0.33\textwidth}
\includegraphics[width=\linewidth]{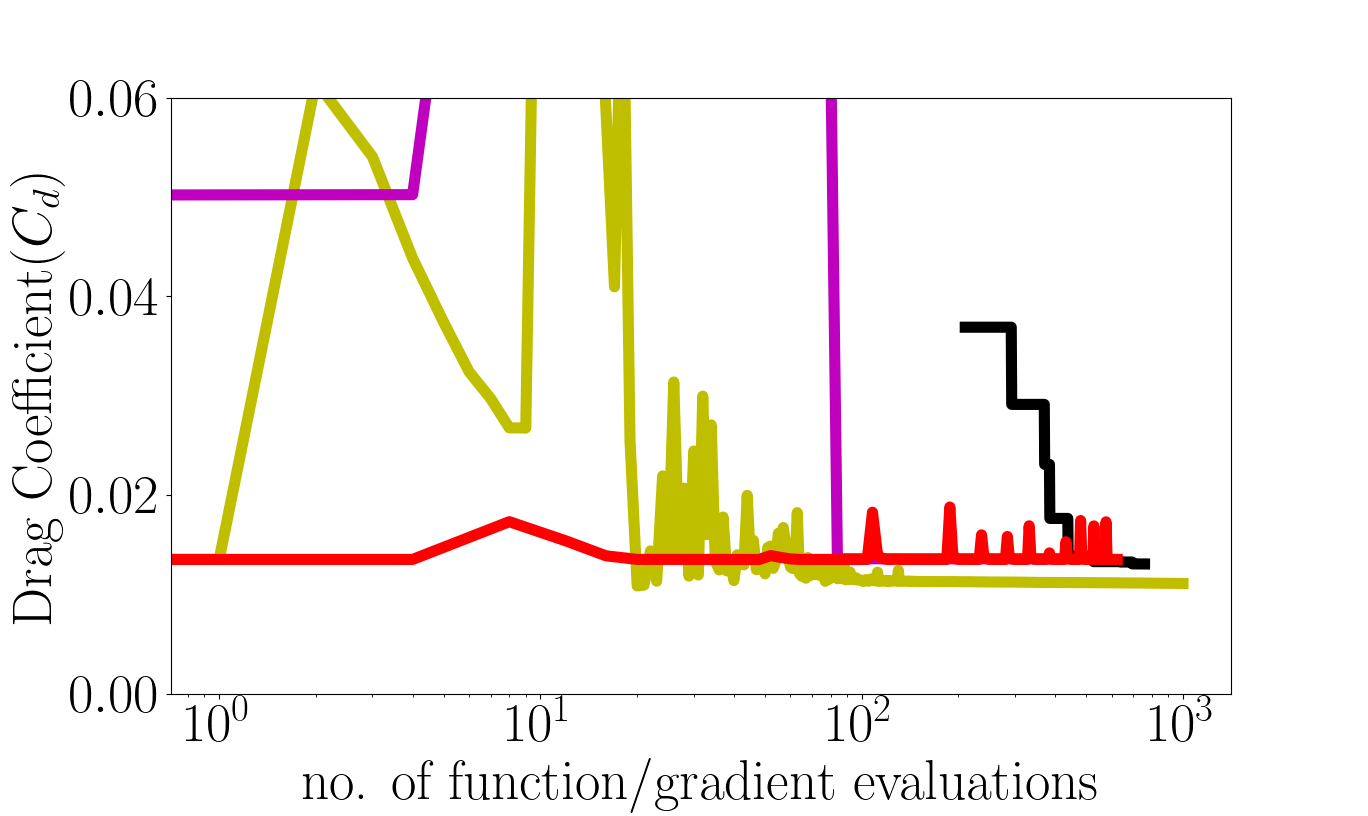} 
\caption{16 control points}\label{feval-16pts-RAE-Constrained}
    \end{subfigure}
    \begin{subfigure}{\textwidth}
        \centering
        \includegraphics[width=0.9\linewidth]{legends_con_RAE.png} 
    \end{subfigure}
    \caption{{\bf Constrained RAE2822.} Convergence histories versus number of function/gradient evaluations. \change{Black are best feasible objective values.}}
    \label{feval-RAE-constrained}
\end{figure}

\begin{figure}[htb!]
    \centering
    \begin{subfigure}{0.33\textwidth}
\includegraphics[width=\linewidth]{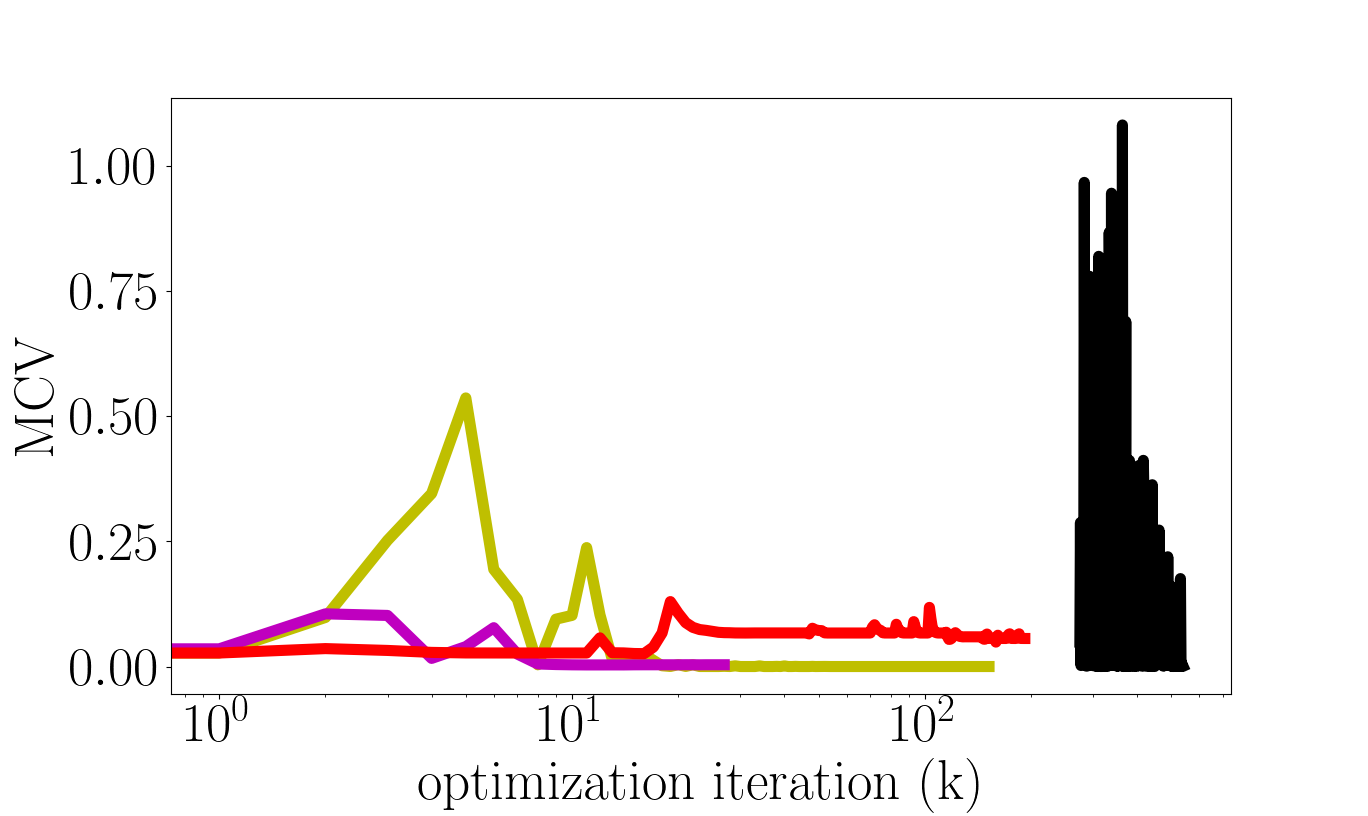}
\caption{4 control points}\label{MCV-4pts-RAE-Constrained}
    \end{subfigure}%
\begin{subfigure}{0.33\textwidth}
\includegraphics[width=\linewidth]{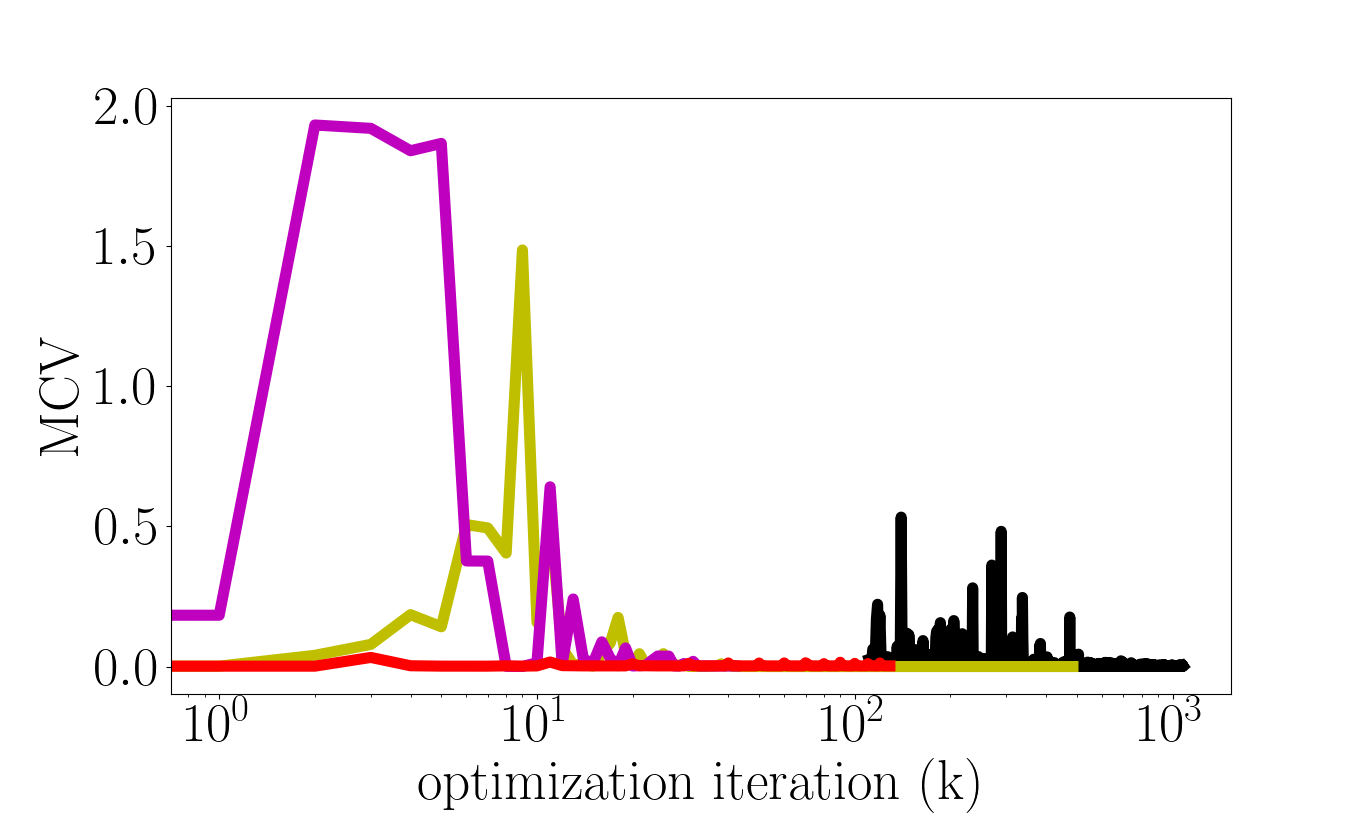} \caption{8 control points}\label{MCV-8pts-RAE-Constrained}
    \end{subfigure}%
    \begin{subfigure}{0.33\textwidth}
\includegraphics[width=\linewidth]{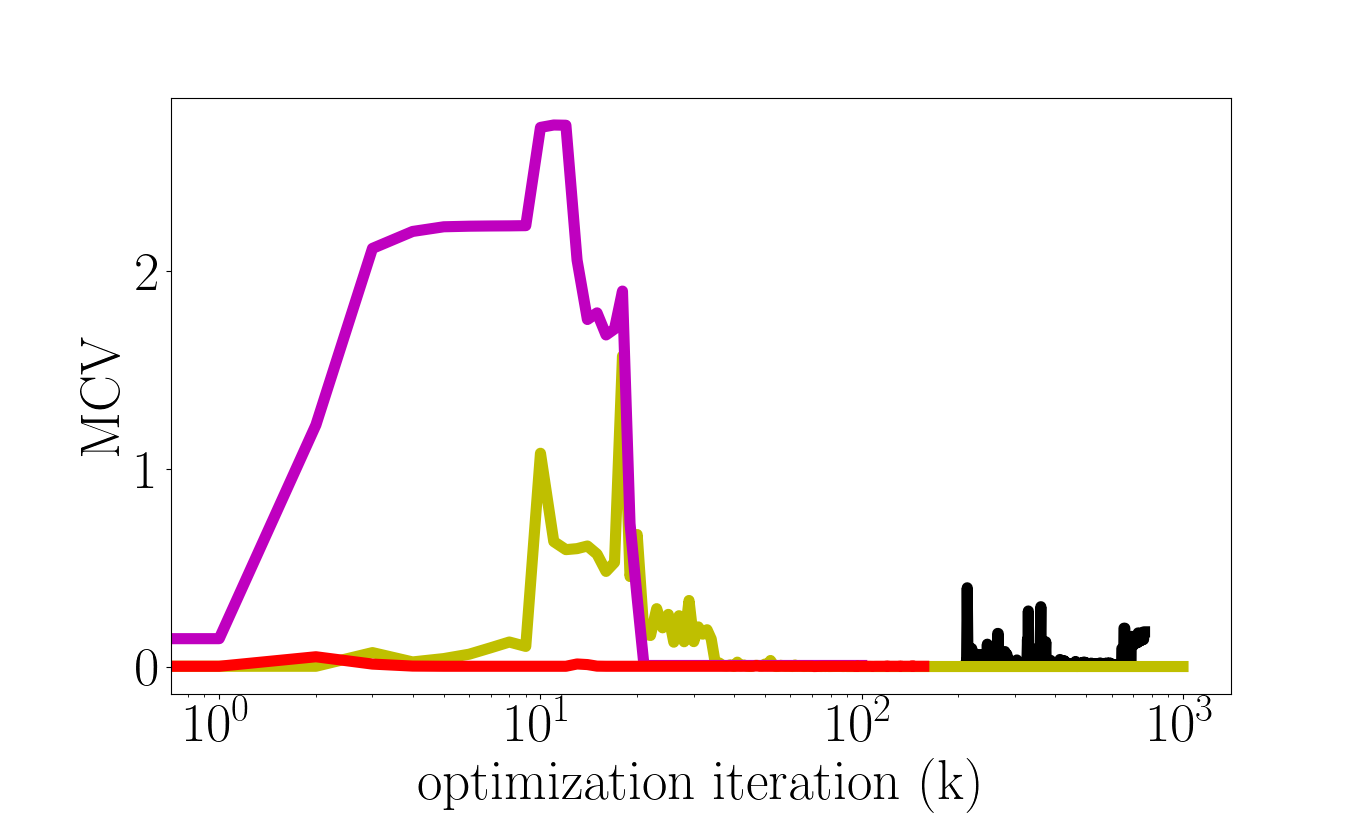} 
\caption{16 control points}\label{MCV-16pts-RAE-Constrained}
    \end{subfigure}
    \begin{subfigure}{\textwidth}
        \centering
        \includegraphics[width=0.9\linewidth]{legends_con_RAE.png} 
    \end{subfigure}
    \caption{{\bf Constrained RAE2822.} Convergence histories of maximum constraint violation (MCV).}
    \label{MCV-RAE-constrained}
\end{figure}

The optimization results, \change{objectives and MCV histories}, are shown in figures \ref{Cd-RAE-constrained}, \ref{feval-RAE-constrained}, and \ref{MCV-RAE-constrained}. In terms of the objective function \change{(figures \ref{Cd-RAE-constrained} and \ref{feval-RAE-constrained})}, notice that derivative-free methods, \cobyla~and \bo, give the lowest objective values among all the methods compared; particularly, \cobyla~and \bo~are quite similar in performance. The maximum constraint violation (MCV) is shown in \Cref{MCV-RAE-constrained}. Notice that \tc~can sometimes (4 control points) fail to find a feasible solution in the constrained setting. In contrast, \bo~constantly attempts to balance exploration and exploitation--as a result, infeasible iterates are still acceptable as they help improve the GP surrogate models. However, for the objective function plots (\Cref{Cd-4pts-RAE-Constrained}, \ref{Cd-8pts-RAE-Constrained} and \ref{Cd-16pts-RAE-Constrained}), we only show the best feasible objective value.

Notice that the history of MCV directly correlates with the convergence behavior of the objective derivative-based algorithms. \slsqp~ experiences extreme constraint violations in early iterations, particularly for high-dimensional cases where MCV spikes coincide with $C_d$ oscillations, highlighting its struggle to maintain feasibility while reducing drag. However, it rapidly improves feasibility, requiring fewer iterations to converge. \tc~ maintains the lowest MCV across all cases, demonstrating better constraint handling, but its conservative approach limits $C_d$ reduction. \cobyla~ shows a moderate decline in constraint violations, resulting in slower convergence. \bo~, on the other hand, achieves the best balance between constraint satisfaction and objective function improvement, avoiding extreme MCV spikes seen in other methods, especially in high-dimensional cases, though occasional constraint violations slow its overall convergence.

\begin{itemize}

\item \textbf{SLSQP / trust-constr (local enforcement):} feasibility is driven by local models (linear/quadratic constraint approximations) and merit/penalty or filter logic with line-search or trust-region acceptance, which can be conservative or brittle when constraints are tight and local models are inaccurate.
\item \textbf{COBYLA (local model, derivative-free):} feasibility is managed via linear interpolation models inside a trust region; it can be robust to gradient issues but still relies on local approximations.
\item \textbf{Constrained BO (global, probabilistic):} feasibility is modeled globally through surrogates for constraints, and candidate selection is based on expected utility under uncertainty (e.g., incorporating probability of feasibility or analogous criteria). This can reduce the frequency of large infeasible steps and can better explore narrow feasible corridors within a fixed evaluation budget.
\end{itemize}

\subsection{Constrained ONERA M6}

Finally, we present results on the constrained optimization of ONERAM6 under inviscid conditions. 
The optimization problem is described in Equation \ref{eq:opt_problem_2}.
\begin{equation} \label{eq:opt_problem_2}
\begin{split}
\min_{\mathbf{x} \in \mcl{X}} \quad & C_D \\
\text{s.t.} \quad & t_1 \geq 0.078, \\
& t_2 \geq 0.072, \\
& t_3 \geq 0.066, \\
& t_4 \geq 0.061, \\
& t_5 \geq 0.055, \\
& C_{L} = 0.292,
\end{split}
\end{equation}
where $t_i$ represents the thickness-to-chord ratio at the $i$th spanwise location and $C_L$ represents the lift coefficient.  

\begin{figure}[htb!]
    \centering
    \begin{subfigure}{0.33\textwidth}
\includegraphics[width=\linewidth]{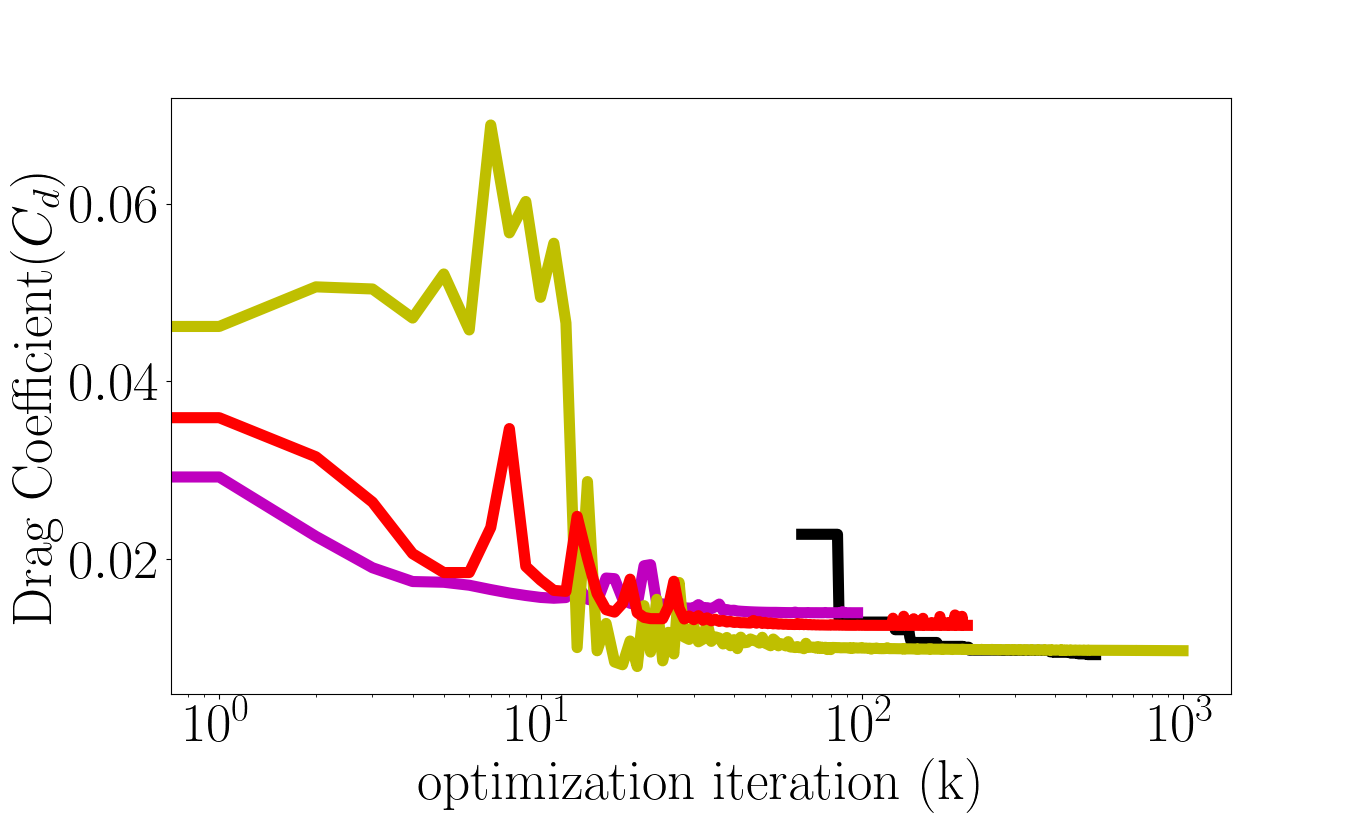}
\caption{$C_D$ vs iteration}\label{Cd-ONERA-12pts-Constrained}
    \end{subfigure}%
\begin{subfigure}{0.33\textwidth}
\includegraphics[width=\linewidth]{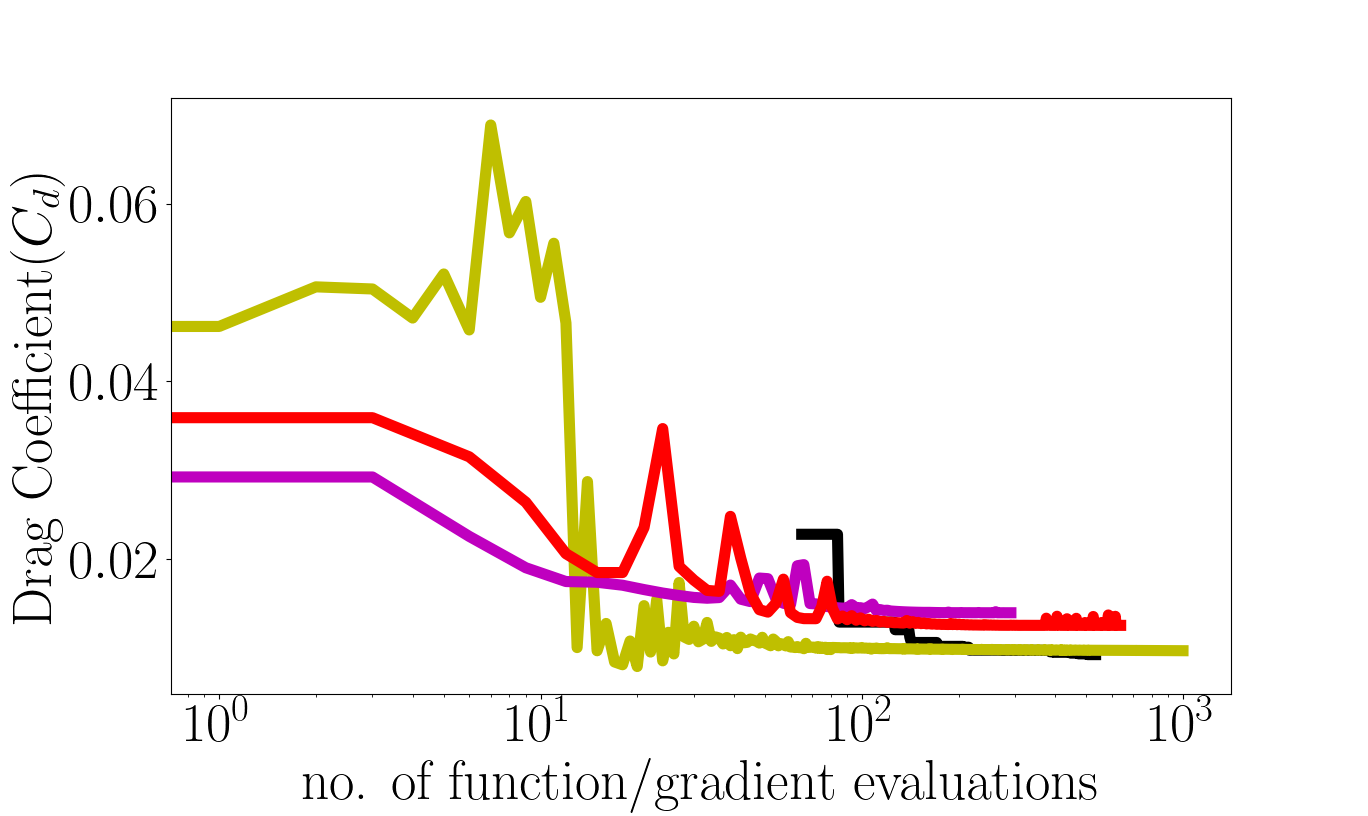} \caption{$C_D$ vs no. of function evals}\label{feval-12pts-ONERA-Constrained}
    \end{subfigure}%
    \begin{subfigure}{0.33\textwidth}
\includegraphics[width=\linewidth]{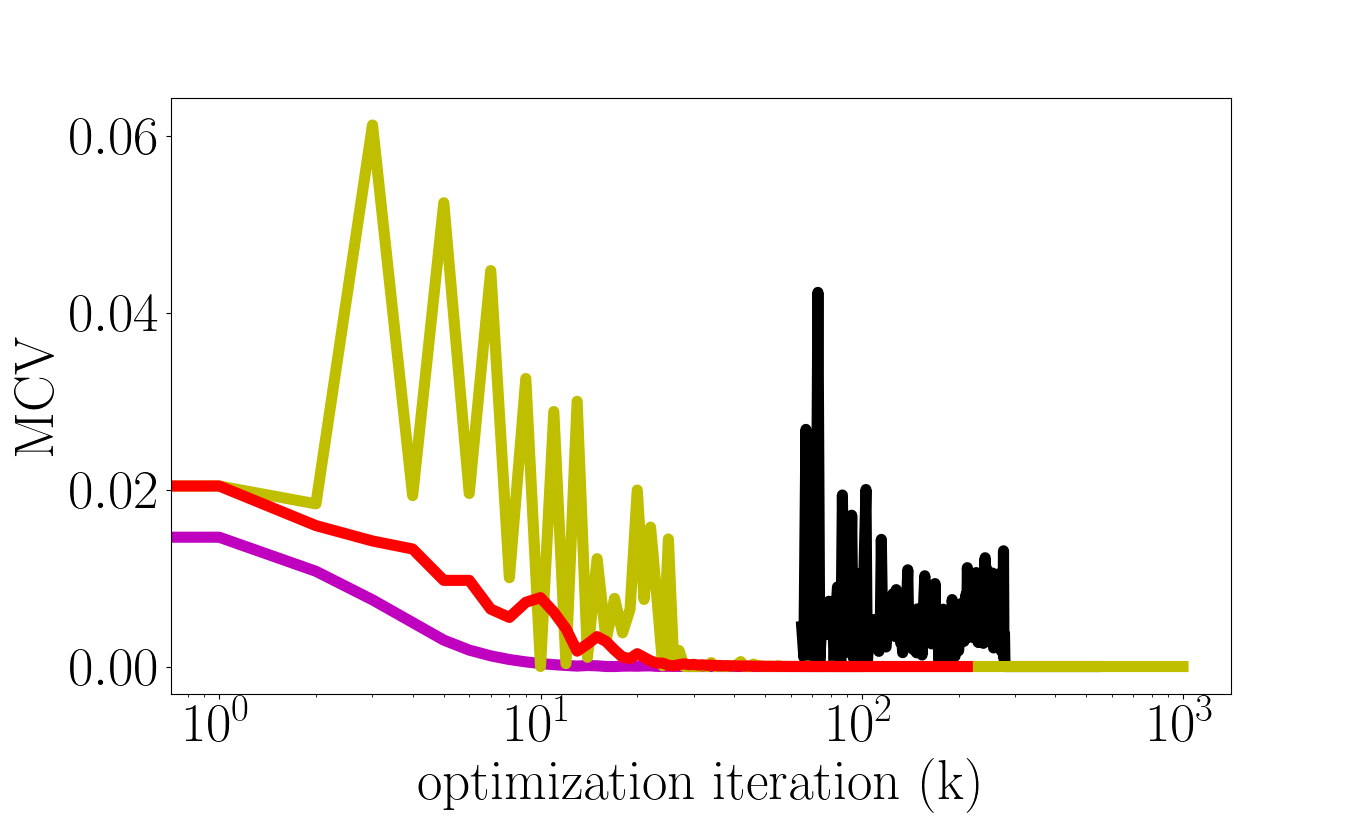} 
\caption{MCV vs iteration}\label{MCV-12pts-ONERA-Constrained}
    \end{subfigure}
    \begin{subfigure}{\textwidth}
        \centering
        \includegraphics[width=0.9\linewidth]{legends_con_RAE.png} 
    \end{subfigure}
    \caption{{\bf Constrained ONERAM6.} Convergence histories.}
    \label{ONERA-constrained}
\end{figure}

\Cref{Cd-ONERA-12pts-Constrained} illustrates the convergence history of the objective ($C_D$). Again, notice that derivative-free methods \change{(\cobyla~and \bo)} outperform derivative-based methods by achieving lower objective values, with \bo~ \change{having and edge over} \cobyla. The MCV plot shows that all methods find a feasible solution relatively quickly (within $\sim 20$ iterations), including \cobyla. This experiment further demonstrates that derivative-free methods can outperform derivative-based methods for three-dimensional problems with higher dimensional parametrization.
\changetwo{It should be noted that the wall times for all optimization algorithms mirrors the trends observed in terms of the number of function and gradient evals. This is because, the cost of the CFD calls (primal and adjoint) were substantially larger than the cost of computing an optimization search step.}
\changetwo{Finally, we summarize the number of function and gradient evaluations of all experiments in Table 3.}

\changethree{One thing worth mentioning is the ability of all derivative-free optimization methods tested in this work to find feasible solutions consistently. For gradient-based local optimization methods, feasibility is driven by local models (linear/quadratic constraint approximations) and merit/penalty or filter logic with line-search or trust-region acceptance, which can be conservative or brittle when constraints are tight and local models are inaccurate. In derivative-free methods like \cobyla~feasibility is managed via linear interpolation models inside a trust region; it can be robust to gradient issues but still relies on local approximations. However, in \bo, feasibility is modeled globally through surrogates for constraints, and candidate selection is based on expected utility under uncertainty (e.g., incorporating probability of feasibility or analogous criteria). This can reduce the frequency of large infeasible steps and can better explore narrow feasible corridors within a fixed evaluation budget.}

\section{Conclusions and future work}
\label{sec:conclusions}
We present the first of a sequence of studies aimed at evaluating the feasibility of derivative-free optimization methods for practical aerodynamic design optimization problems. To our knowledge, a systematic benchmarking of derivative-free approaches against derivative-based methods doesn't exist -- this work fills that gap. Based on two airfoil and one wing test case, up to $32$ dimensions, we benchmark three derivative-free algorithms (\nm, \cobyla, and \bo) against four derivative-based approaches (\bfgs, \slsqp, \tc, \tnc), with and without constraints. Our experiments reveal that the additional cost invested in computing gradients does not necessarily translate into efficient convergence to a stationary point; indeed, in the majority of our experiments, we observe that derivative-based methods fail to achieve first-order optimality. Furthermore, in the constrained setting, derivative-based methods do not always find feasible designs. \change{We measure algorithm performance in terms of identifying the lowest feasible objective value with the fewest evaluations of the objective and constraint functions. In this regard, DFO methods \nm~and \bo~outperform gradient-based methods in the unconstrained RAE2822 and ONERAM6 experiments; DFO methods \cobyla~and \bo~outperform gradient based methods in the constrained RAE2822 and ONERAM6 experiments -- this amounts to four out of the five experiments considered. This establishes the competitiveness of DFO methods in aerodynamic design optimization problems in moderate dimensions.}

We note that, for some of the derivative-based algorithms, we are conservative in estimating the number of function and gradient evaluations to be two per iteration; in practice, computing line search step lengths can cost much more than that~\cite{nocedal1999numerical}. Despite this, a competitive performance demonstrated by derivative-free methods provides substantial evidence that they are a realistic option for practical aerodynamic design optimization. Crucially, surrogate-based methods such as \bo~are particularly attractive due to their automatic ability to balance exploration and exploitation, despite depending only on zeroth-order information. It is worth noting that \bo~offers further flexibility in terms of the surrogate model and acquisition policy choices, which provide scope for further improvement in performance. \changethree{Within the prescribed budgets and solver tolerances, the tested gradient-based solvers
did not consistently drive the reported gradient norm/KKT residual to small values,
which is consistent with sensitivity to gradient/constraint inconsistencies and/or local
non-smoothness in the CFD-based objectives/constraints.}

Naturally, one of the primary directions for future work is to perform the same study in higher dimensions to observe whether the conclusions still hold. While it is known that global surrogate models, such as GPs, scale poorly with dimensions, several measures already exist to scale them to high dimensions~\cite{eriksson2021high,eriksson2019scalable,li2016high,wang2016bayesian,nayebi2019framework,lizotte2008practical,renganathan2025qpots}. We plan to leverage existing methods to benchmark \bo~on high-dimensional aerodynamic design optimization problems. Another direction of future work is to 
evaluate the use of non-stationary surrogates (e.g., deep Gaussian processes~\cite{booth2024actively, booth2025contour, rajaram2020deep,rajaram2021empirical}) for their potential to improve \bo~performance. 

\section*{Appendix}

\begin{longtable}{ll l r r r}
\caption{Summary of total function evaluations ($n_f$) and gradient evaluations ($n_g$) for each test case and optimizer. For gradient-based methods, we assume $n_f = n_g$.}\\
\toprule
Case & $d$ & Optimizer & $n_f$ & $n_g$ & $n_f{+}n_g$\\
\midrule
\endfirsthead
\toprule
Case & $d$ & Optimizer & $n_f$ & $n_g$ & $n_f{+}n_g$\\
\midrule
\endhead
\midrule \multicolumn{6}{r}{\emph{Continued on next page}}\\
\endfoot
\bottomrule
\endlastfoot
\multicolumn{6}{l}{\textbf{NACA0012 (unconstrained) (Fig.~3)}}\\
\midrule
\multirow{6}{*}{NACA0012} & \multirow{6}{*}{4} & L-BFGS-B & 27 & 27 & 54\\
 &  & SLSQP & 20 & 20 & 40\\
 &  & TNC & 67 & 67 & 134\\
 &  & Trust-Constrained & 177 & 177 & 354\\
 &  & Nelder-Mead & 132 & 0 & 132\\
 &  & BO & 59 & 0 & 59\\
\addlinespace[0.5ex]
\multirow{6}{*}{NACA0012} & \multirow{6}{*}{8} & L-BFGS-B & 9 & 10 & 19\\
 &  & SLSQP & 39 & 39 & 78\\
 &  & TNC & 47 & 48 & 95\\
 &  & Trust-Constrained & 161 & 161 & 322\\
 &  & Nelder-Mead & 80 & 0 & 80\\
 &  & BO & 30 & 0 & 30\\
\addlinespace[0.5ex]
\multirow{6}{*}{NACA0012} & \multirow{6}{*}{16} & L-BFGS-B & 16 & 16 & 32\\
 &  & SLSQP & 34 & 34 & 68\\
 &  & TNC & 94 & 94 & 188\\
 &  & Trust-Constrained & 189 & 190 & 379\\
 &  & Nelder-Mead & 1311 & 0 & 1311\\
 &  & BO & 441 & 0 & 441\\
\addlinespace[0.5ex]
\midrule
\multicolumn{6}{l}{\textbf{RAE2822 (unconstrained) (Fig.~6)}}\\
\midrule
\multirow{5}{*}{RAE2822} & \multirow{5}{*}{4} & L-BFGS-B & 8 & 8 & 16\\
 &  & TNC & 25 & 26 & 51\\
 &  & Trust-Constrained & 131 & 131 & 262\\
 &  & Nelder-Mead & 326 & 0 & 326\\
 &  & BO & 137 & 0 & 137\\
\addlinespace[0.5ex]
\multirow{5}{*}{RAE2822} & \multirow{5}{*}{8} & L-BFGS-B & 14 & 14 & 28\\
 &  & TNC & 73 & 73 & 146\\
 &  & Trust-Constrained & 148 & 149 & 297\\
 &  & Nelder-Mead & 515 & 0 & 515\\
 &  & BO & 285 & 0 & 285\\
\addlinespace[0.5ex]
\multirow{5}{*}{RAE2822} & \multirow{5}{*}{16} & L-BFGS-B & 7 & 8 & 15\\
 &  & TNC & 140 & 140 & 280\\
 &  & Trust-Constrained & 173 & 173 & 346\\
 &  & Nelder-Mead & 866 & 0 & 866\\
 &  & BO & 472 & 0 & 472\\
\addlinespace[0.5ex]
\multirow{5}{*}{RAE2822} & \multirow{5}{*}{32} & L-BFGS-B & -- & -- & --\\
 &  & TNC & 131 & 131 & 262\\
 &  & Trust-Constrained & 130 & 130 & 260\\
 &  & Nelder-Mead & 1262 & 0 & 1262\\
 &  & BO & 625 & 0 & 625\\
\addlinespace[0.5ex]
\midrule
\multicolumn{6}{l}{\textbf{ONERAM6 (unconstrained) (Fig.~9c)}}\\
\midrule
\multirow{5}{*}{ONERAM6} & \multirow{5}{*}{12} & L-BFGS-B & 19 & 19 & 38\\
 &  & TNC & 146 & 147 & 293\\
 &  & Trust-Constrained & 279 & 279 & 558\\
 &  & Nelder-Mead & 705 & 0 & 705\\
 &  & BO & 434 & 0 & 434\\
\addlinespace[0.5ex]
\midrule
\multicolumn{6}{l}{\textbf{RAE2822 (constrained) (Fig.~11)}}\\
\midrule
\multirow{4}{*}{RAE2822} & \multirow{4}{*}{4} & COBYLA & 1030 & 0 & 1030\\
 &  & SLSQP & 186 & 186 & 372\\
 &  & Trust-Constrained & 423 & 423 & 846\\
 &  & BO & 1030 & 0 & 1030\\
\addlinespace[0.5ex]
\multirow{4}{*}{RAE2822} & \multirow{4}{*}{8} & COBYLA & 1091 & 0 & 1091\\
 &  & SLSQP & 161 & 161 & 322\\
 &  & Trust-Constrained & 399 & 399 & 798\\
 &  & BO & 1142 & 0 & 1142\\
\addlinespace[0.5ex]
\multirow{4}{*}{RAE2822} & \multirow{4}{*}{16} & COBYLA & 491 & 0 & 491\\
 &  & SLSQP & 125 & 125 & 250\\
 &  & Trust-Constrained & 342 & 342 & 684\\
 &  & BO & 788 & 0 & 788\\
\addlinespace[0.5ex]
\midrule
\multicolumn{6}{l}{\textbf{ONERAM6 (constrained) (Fig.~13b)}}\\
\midrule
\multirow{4}{*}{ONERAM6} & \multirow{4}{*}{12} & COBYLA & 1046 & 0 & 1046\\
 &  & SLSQP & 153 & 154 & 307\\
 &  & Trust-Constrained & 336 & 336 & 672\\
 &  & BO & 543 & 0 & 543\\
\addlinespace[0.5ex]
\label{tab:summary-evals}
\end{longtable}

\clearpage
\section{Statements and Declarations}

\subsection{Funding} N/A
\subsection{Conflict of Interest}The authors have no conflicts of interest to disclose.
\subsection{Author Contributions}P. Plaban -- CFD preparation and analysis, derivative-based algorithms, manuscript drafting. P. Bachman -- CFD analysis and derivative-free algorithms. A. Renganathan -- design of experiments, manuscript drafting, and research supervision.
\subsection{Data Availability} N/A
\subsection{Ethics Approval and Consent to Participate} Yes
\subsection{Replication of Results}Code publicly available at \url{https://github.com/csdlpsu/aso}.

\section{Replication of results}
The code for running our experiments and replicating the results is publicly available at \url{https://github.com/csdlpsu/aso}.

\clearpage
\bibliographystyle{plainnat}
\bibliography{references}

@inproceedings{gardner2014bayesian,
  title={Bayesian optimization with inequality constraints.},
  author={Gardner, Jacob R and Kusner, Matt J and Xu, Zhixiang Eddie and Weinberger, Kilian Q and Cunningham, John P},
  booktitle={ICML},
  volume={2014},
  pages={937--945},
  year={2014}
}

@article{gill2005snopt,
  title={SNOPT: An SQP algorithm for large-scale constrained optimization},
  author={Gill, Philip E and Murray, Walter and Saunders, Michael A},
  journal={SIAM review},
  volume={47},
  number={1},
  pages={99--131},
  year={2005},
  publisher={SIAM}
}

@book{omojokun1989trust,
  title={Trust region algorithms for optimization with nonlinear equality and inequality constraints},
  author={Omojokun, Emmanuel Omotayo},
  year={1989},
  publisher={University of Colorado at Boulder}
}

@article{virtanen2020scipy,
  title={SciPy 1.0: fundamental algorithms for scientific computing in Python},
  author={Virtanen, Pauli and Gommers, Ralf and Oliphant, Travis E and Haberland, Matt and Reddy, Tyler and Cournapeau, David and Burovski, Evgeni and Peterson, Pearu and Weckesser, Warren and Bright, Jonathan and others},
  journal={Nature methods},
  volume={17},
  number={3},
  pages={261--272},
  year={2020},
  publisher={Nature Publishing Group US New York}
}

@article{letham2019constrained,
  title={Constrained Bayesian optimization with noisy experiments},
  author={Letham, Benjamin and Karrer, Brian and Ottoni, Guilherme and Bakshy, Eytan},
  year={2019}
}

@article{renganathan2021enhanced,
  title={Enhanced data efficiency using deep neural networks and Gaussian processes for aerodynamic design optimization},
  author={Renganathan, S Ashwin and Maulik, Romit and Ahuja, Jai},
  journal={Aerospace Science and Technology},
  volume={111},
  pages={106522},
  year={2021},
  publisher={Elsevier},
  doi={10.1016/j.ast.2021.106522}
}

@book{nocedal1999numerical,
  title={Numerical optimization},
  author={Nocedal, Jorge and Wright, Stephen J},
  year={1999},
  publisher={Springer}
}

@inproceedings{eriksson2021high,
  title={High-dimensional Bayesian optimization with sparse axis-aligned subspaces},
  author={Eriksson, David and Jankowiak, Martin},
  booktitle={Uncertainty in Artificial Intelligence},
  pages={493--503},
  year={2021},
  organization={PMLR}
}

@article{eriksson2019scalable,
  title={Scalable global optimization via local Bayesian optimization},
  author={Eriksson, David and Pearce, Michael and Gardner, Jacob and Turner, Ryan D and Poloczek, Matthias},
  journal={Advances in neural information processing systems},
  volume={32},
  year={2019}
}

@inproceedings{li2016high,
  title={High dimensional Bayesian optimization via restricted projection pursuit models},
  author={Li, Chun-Liang and Kandasamy, Kirthevasan and P{\'o}czos, Barnab{\'a}s and Schneider, Jeff},
  booktitle={Artificial Intelligence and Statistics},
  pages={884--892},
  year={2016},
  organization={PMLR}
}

@article{wang2016bayesian,
  title={Bayesian optimization in a billion dimensions via random embeddings},
  author={Wang, Ziyu and Hutter, Frank and Zoghi, Masrour and Matheson, David and De Feitas, Nando},
  journal={Journal of Artificial Intelligence Research},
  volume={55},
  pages={361--387},
  year={2016},
  doi={10.1613/jair.4806}
}

@inproceedings{nayebi2019framework,
  title={A framework for Bayesian optimization in embedded subspaces},
  author={Nayebi, Amin and Munteanu, Alexander and Poloczek, Matthias},
  booktitle={International Conference on Machine Learning},
  pages={4752--4761},
  year={2019},
  organization={PMLR}
}

@inproceedings{renganathan2025qpots,
  title={qPOTS: Efficient Batch Multiobjective Bayesian Optimization via Pareto Optimal Thompson Sampling},
  author={Renganathan, Ashwin and Carlson, Kade},
  booktitle={International Conference on Artificial Intelligence and Statistics},
  pages={4051--4059},
  year={2025},
  organization={PMLR}
}

@article{lizotte2008practical,
  title={Practical bayesian optimization},
  author={Lizotte, Daniel James},
  year={2008},
  isbn      = { 978-0-494-46365-9}
}

@inproceedings{marchildon2024,
  author    = {André L. Marchildon and David W. Zingg},
  title     = {Gradient-Enhanced Bayesian Optimization with Application to Aerodynamic Shape Optimization},
  booktitle = {AIAA Aviation Forum and ASCEND 2024},
  year      = {2024},
  doi       = {10.2514/6.2024-4405}
}

@article{booth2025contour,
  title={Contour location for reliability in airfoil simulation experiments using deep gaussian processes},
  author={Booth, Annie S and Renganathan, S Ashwin and Gramacy, Robert B},
  journal={The Annals of Applied Statistics},
  volume={19},
  number={1},
  pages={191--211},
  year={2025},
  publisher={Institute of Mathematical Statistics},
  doi={10.1214/24-AOAS1951}
}

@inproceedings{booth2024actively,
  title={Actively learning deep Gaussian process models for failure contour and probability estimation.},
  author={Booth, Annie S and Gramacy, Robert and Renganathan, Ashwin},
  booktitle={AIAA SCITECH 2024 Forum},
  pages={0577},
  year={2024},
  doi={10.2514/6.2024-0577}
}

@article{rajaram2021empirical,
  title={Empirical assessment of deep gaussian process surrogate models for engineering problems},
  author={Rajaram, Dushhyanth and Puranik, Tejas G and Ashwin Renganathan, S and Sung, WoongJe and Fischer, Olivia Pinon and Mavris, Dimitri N and Ramamurthy, Arun},
  journal={Journal of Aircraft},
  volume={58},
  number={1},
  pages={182--196},
  year={2021},
  publisher={American Institute of Aeronautics and Astronautics},
  doi={10.2514/1.C036026}
}

@inproceedings{rajaram2020deep,
  title={Deep Gaussian process enabled surrogate models for aerodynamic flows},
  author={Rajaram, Dushhyanth and Puranik, Tejas G and Renganathan, Ashwin and Sung, WoongJe and Pinon-Fischer, Olivia J and Mavris, Dimitri N and Ramamurthy, Arun},
  booktitle={AIAA scitech 2020 forum},
  pages={1640},
  year={2020},
  doi={10.2514/6.2020-1640}
}

@article{shahriari2015taking,
  title={Taking the human out of the loop: A review of Bayesian optimization},
  author={Shahriari, Bobak and Swersky, Kevin and Wang, Ziyu and Adams, Ryan P and De Freitas, Nando},
  journal={Proceedings of the IEEE},
  volume={104},
  number={1},
  pages={148--175},
  year={2015},
  publisher={IEEE},
  doi={10.1109/JPROC.2015.2494218}
}

@article{greenhill2020bayesian,
  title={Bayesian optimization for adaptive experimental design: A review},
  author={Greenhill, Stewart and Rana, Santu and Gupta, Sunil and Vellanki, Pratibha and Venkatesh, Svetha},
  journal={IEEE access},
  volume={8},
  pages={13937--13948},
  year={2020},
  publisher={IEEE},
  doi={ 10.1109/ACCESS.2020.2966228}
}

@article{thompson1933likelihood,
  title={On the likelihood that one unknown probability exceeds another in view of the evidence of two samples},
  author={Thompson, William R},
  journal={Biometrika},
  volume={25},
  number={3-4},
  pages={285--294},
  year={1933},
  publisher={Oxford University Press},
  doi={10.2307/2332286}
}

@article{jones1998efficient,
  doi          = {10.1023/A:1008306431147},
  title={Efficient global optimization of expensive black-box functions},
  author={Jones, Donald R. and Schonlau, Matthias and Welch, William J.},
  journal={Journal of Global Optimization},
  volume={13},
  number={4},
  pages={455--492},
  year={1998},
  publisher={Springer}
}

@inproceedings{hernandez2014predictive,
  title={Predictive entropy search for efficient global optimization of black-box functions},
  author={Hern{\'a}ndez-Lobato, Jos{\'e} Miguel and Hoffman, Matthew W. and Ghahramani, Zoubin},
  booktitle={Advances in Neural Information Processing Systems},
  pages={918--926},
  year={2014}
}

@inproceedings{wang2017max,
  title={Max-value entropy search for efficient {Bayesian} optimization},
  author={Wang, Zi and Jegelka, Stefanie},
  booktitle={Proceedings of the 34th International Conference on Machine Learning},
  volume={70},
  pages={3627--3635},
  year={2017},
  organization={PMLR}
}

@inproceedings{hernandez2016predictive,
  title={Predictive entropy search for multi-objective bayesian optimization},
  author={Hern{\'a}ndez-Lobato, Daniel and Hernandez-Lobato, Jose and Shah, Amar and Adams, Ryan},
  booktitle={International conference on machine learning},
  pages={1492--1501},
  year={2016},
  organization={PMLR}
}

@article{renganathan2023camera,
  title={CAMERA: A method for cost-aware, adaptive, multifidelity, efficient reliability analysis},
  author={Renganathan, S Ashwin and Rao, Vishwas and Navon, Ionel M},
  journal={Journal of Computational Physics},
  volume={472},
  pages={111698},
  year={2023},
  publisher={Elsevier},
  doi={10.1016/j.jcp.2022.111698}
}

@inproceedings{renganathan2024efficient,
  title={Efficient reliability analysis with multifidelity Gaussian processes and normalizing flows},
  author={Renganathan, Ashwin},
  booktitle={AIAA SCITECH 2024 Forum},
  pages={0576},
  year={2024},
  doi={10.2514/6.2024-0576}
}

@article{balandat2020botorch,
  title={BoTorch: A framework for efficient Monte-Carlo Bayesian optimization},
  author={Balandat, Maximilian and Karrer, Brian and Jiang, Daniel and Daulton, Samuel and Letham, Ben and Wilson, Andrew G and Bakshy, Eytan},
  journal={Advances in neural information processing systems},
  volume={33},
  pages={21524--21538},
  year={2020}, 
  doi={}

}

@article{lewis1999pattern,
  title={Pattern search algorithms for bound constrained minimization},
  author={Lewis, Robert Michael and Torczon, Virginia},
  journal={SIAM Journal on optimization},
  volume={9},
  number={4},
  pages={1082--1099},
  year={1999},
  publisher={SIAM},
  doi={10.1137/S1052623496300507}
}

@article{olsson1975nelder,
  title={The Nelder-Mead simplex procedure for function minimization},
  author={Olsson, Donald M and Nelson, Lloyd S},
  journal={Technometrics},
  volume={17},
  number={1},
  pages={45--51},
  year={1975},
  publisher={Taylor \& Francis}
}

@inproceedings{meo2024,
  author    = {Azzurra Meo and Alessandro Aponte and Alessandro Munafò and Giulio Gori and Marco Panesi},
  title     = {A Bayesian Framework for the Shape Optimization of Atmospheric Reentry Vehicles},
  booktitle = {AIAA Aviation Forum and ASCEND, 2024},
  year      = {2024},
  publisher = {American Institute of Aeronautics and Astronautics Inc, AIAA},
  doi       = {10.2514/6.2024-3563}
}

@article{tfaily2024,
  author      = {Ali Tfaily and Youssef Diouane and Nathalie Bartoli and Michael Kokkolaras},
  title       = {Bayesian Optimization with Hidden Constraints for Aircraft Design},
  institution = {Les Cahiers du GERAD},
  number      = {G–2024–10},
  year        = {2024},
  journal={Structural and Multidisciplinary Optimization},
  volume={67},
  pages={123},
  publisher={Springer},
  doi         = {10.1007/s00158-024-03833-8}
}

@article{reist2020,
  author    = {Thomas A. Reist and David Koo and David W. Zingg and Pascal Bochud and Philippe Castonguay and David Leblond},
  title     = {Cross-Validation of High-Fidelity Aerodynamic Shape Optimization Methodologies for Aircraft Wing-Body Optimization},
  journal   = {AIAA Journal},
  volume    = {58},
  number    = {6},
  pages     = {2581--2595},
  year      = {2020},
  doi       = {10.2514/1.J059091}
}

@article{jim2021,
  author    = {Timothy M. S. Jim and Ghifari A. Faza and Pramudita S. Palar and Koji Shimoyama},
  title     = {Bayesian Optimization of a Low-Boom Supersonic Wing Planform},
  journal   = {AIAA Journal},
  volume    = {59},
  number    = {11},
  pages     = {4514--4529},
  year      = {2021},
  doi       = {10.2514/1.J060225}
}

@article{liu2022,
  author    = {Ruo-Lin Liu and Qiang Zhao and Xian-Jun He and Xin-Yi Yuan and Wei-Tao Wu and Ming-Yu Wu},
  title     = {Airfoil Optimization Based on Multi-Objective Bayesian Optimization},
  journal   = {Journal of Mechanical Science and Technology},
  volume    = {36},
  number    = {11},
  pages     = {5561--5573},
  year      = {2022},
  doi       = {10.1007/s12206-022-1020-y}
}

@article{xu2024,
  author    = {D. Y. Xu and Y. Shen and W. Huang and Z. Y. Guo and H. Zhang and D. F. Xu},
  title     = {Parametric Modelling and Variable-Fidelity Bayesian Optimization of Aerodynamics for a Reusable Flight Vehicle},
  journal   = {Fluid Dynamics},
  volume    = {59},
  number    = {6},
  pages     = {123--135},
  year      = {2024},
  doi       = {10.1134/S0015462824603814}
}

@inproceedings{erhard2024,
  author    = {Racheal M. Erhard and Juan J. Alonso},
  title     = {Multi-Fidelity Bayesian Optimization of a Coaxial Rotor for eVTOL Aircraft},
  booktitle = {AIAA SciTech Forum},
  year      = {2024},
  doi       = {10.2514/6.2024-2506}
}

@article{kundu2018,
  author    = {Abhishek Kundu and H. G. Matthies and M. I. Friswell},
  title     = {Probabilistic Optimization of Engineering Systems with Prescribed Target Design in a Reduced Parameter Space},
  journal   = {Computational Methods in Applied Mechanics and Engineering},
  volume    = {339},
  pages     = {640--660},
  year      = {2018},
  doi       = {10.1016/j.cma.2018.03.041}
}

@article{jameson1988control,
  author    = {Antony Jameson},
  title     = {Aerodynamic Design via Control Theory},
  journal   = {Journal of Scientific Computing},
  volume    = {3},
  number    = {3},
  pages     = {233--260},
  year      = {1988},
  doi       = {10.1007/BF01061285}
}

@article{nelder1965simplex,
  author    = {J. A. Nelder and R. Mead},
  title     = {A Simplex Method for Function Minimization},
  journal   = {The Computer Journal},
  volume    = {7},
  number    = {4},
  pages     = {308--313},
  year      = {1965},
  publisher = {Oxford University Press},
  doi       = {10.1093/comjnl/7.4.308}
}

@incollection{powell1994direct,
  author    = {M. J. D. Powell},
  title     = {A Direct Search Optimization Method That Models the Objective and Constraint Functions by Linear Interpolation},
  booktitle = {Advances in Optimization and Numerical Analysis},
  editor    = {S. Gomez and J.-P. Hennart},
  series    = {Mathematics and Its Applications},
  volume    = {275},
  pages     = {51--67},
  year      = {1994},
  publisher = {Springer},
  doi       = {10.1007/978-94-015-8330-5_4}
}

@article{jones1993direct,
  author    = {Jones, Donald R and Perttunen, Cary D and Stuckman, Bruce E},
  title     = {Lipschitzian Optimization Without the Lipschitz Constant},
  journal   = {Journal of Optimization Theory and Applications},
  volume    = {79},
  number    = {1},
  pages     = {157--181},
  year      = {1993},
  doi       = {10.1007/BF00941892}
}

@article{li2020efficient,
  title     = {Efficient Aerodynamic Shape Optimization with Deep-Learning-Based Geometric Filtering},
  author    = {Li, Jichao and Zhang, Mengqi and Martins, Joaquim R. R. A. and Shu, Chang},
  journal   = {AIAA Journal},
  volume    = {58},
  number    = {10},
  pages     = {4243--4259},
  year      = {2020},
  publisher = {American Institute of Aeronautics and Astronautics},
  doi       = {10.2514/1.J059254}
}

@inproceedings{chen2019aerodynamic,
  author    = {Wei Chen and Kevin Chiu and Mark Fuge},
  title     = {Aerodynamic Design Optimization and Shape Exploration Using Generative Adversarial Networks},
  booktitle = {AIAA SciTech Forum},
  year      = {2019},
  doi       = {10.2514/6.2019-2351}
}

@article{sheikh2022optimization,
  title     = {Optimization of the Shape of a Hydrokinetic Turbine's Draft Tube and Hub Assembly Using Design-by-Morphing with Bayesian Optimization},
  author    = {Sheikh, Haris Moazam and Callan, Tess A. and Hennessy, Kealan J. and Marcus, Philip S.},
  journal   = {Computer Methods in Applied Mechanics and Engineering},
  volume    = {401},
  pages     = {115654},
  year      = {2022},
  publisher = {Elsevier},
  doi       = {10.1016/j.cma.2022.115654}
}

@article{li2022machine,
  title     = {Machine Learning in Aerodynamic Shape Optimization},
  author    = {Li, Jichao and Du, Xiaosong and Martins, Joaquim R. R. A.},
  journal   = {Progress in Aerospace Sciences},
  volume    = {134},
  pages     = {100849},
  year      = {2022},
  publisher = {Elsevier},
  doi       = {10.1016/j.paerosci.2022.100849}
}

@article{bouhlel2019python,
  title={A Python surrogate modeling framework with derivatives},
  author={Bouhlel, Mohamed Amine and Hwang, John T and Bartoli, Nathalie and Lafage, R{\'e}mi and Morlier, Joseph and Martins, Joaquim RRA},
  journal={Advances in Engineering Software},
  volume={135},
  pages={102662},
  year={2019},
  publisher={Elsevier},
  doi={10.1016/j.advengsoft.2019.03.005}
}

@article{queipo2005surrogate,
  author  = {Queipo, Nestor V. and Haftka, Raphael T. and Shyy, Wei and Goel, Tapabrata and Vaidyanathan, Rajiv and Tucker, Paul K.},
  title   = {Surrogate-Based Analysis and Optimization},
  journal = {Progress in Aerospace Sciences},
  volume  = {41},
  number  = {1},
  pages   = {1--28},
  year    = {2005},
  doi     = {10.1016/j.paerosci.2005.02.001}
}

@article{sederberg1986ffd,
  author    = {Sederberg, Thomas W. and Parry, Scott R.},
  title     = {Free-form deformation of solid geometric models},
  booktitle={Proceedings of the 13th annual conference on Computer graphics and interactive techniques},
  journal = {SIGGRAPH Comput. Graph.},
  volume    = {20},
  number    = {4},
  pages     = {151--160},
  year      = {1986},
  doi       = {10.1145/15922.15903}
}

@article{hsu1992direct,
  author    = {Hsu, W. M. and Hughes, J. F. and Kaufman, H.},
  title     = {Direct manipulation of free-form deformations},
  journal   = {ACM SIGGRAPH Computer Graphics},
  volume    = {26},
  number    = {2},
  pages     = {177--184},
  year      = {1992},
  doi       = {10.1145/142920.134036}
}

@article{frazier2018tutorial,
  author  = {Frazier, Peter I.},
  title   = {A Tutorial on Bayesian Optimization},
  journal = {arXiv preprint arXiv:1807.02811},
  year    = {2018},
  doi     = {10.48550/arXiv.1807.02811}
}

@book{rasmussen2006gpml,
  author    = {Rasmussen, Carl Edward and Williams, Christopher K. I.},
  title     = {Gaussian Processes for Machine Learning},
  publisher = {The MIT Press},
  year      = {2006},
  isbn      = {9780262182539},
  doi       = {10.7551/mitpress/3206.001.0001}
}

@book{garnett2023bo,
  author    = {Garnett, Roman},
  title     = {Bayesian Optimization},
  publisher = {Cambridge University Press},
  year      = {2023},
  isbn      = {9781108425780},
  doi       = {10.1017/9781108348973}
}

@incollection{palacios_stanford_2013,
	title = {Stanford {University} {Unstructured} ({SU}{\textasciicircum}2): {An} open-source integrated computational environment for multi-physics simulation and design},
	booktitle = {51st {AIAA} {Aerospace} {Sciences} {Meeting} including the {New} {Horizons} {Forum} and {Aerospace} {Exposition}},
	publisher = {American Institute of Aeronautics and Astronautics},
	author = {Palacios, Francisco and Alonso, Juan and Duraisamy, Karthikeyan and Colonno, Michael and Hicken, Jason and Aranake, Aniket and Campos, Alejandro and Copeland, Sean and Economon, Thomas and Lonkar, Amrita and Lukaczyk, Trent and Taylor, Thomas},
	year = {2013},
	doi = {10.2514/6.2013-287}
}

@article{jameson2017origins,
  title={Origins and further development of the Jameson--Schmidt--Turkel scheme},
  author={Jameson, Antony},
  journal={AIAA Journal},
  volume={55},
  number={5},
  pages={1487--1510},
  year={2017},
  publisher={American Institute of Aeronautics and Astronautics},
  doi={doi.org/10.2514/1.J055493}
}

@article{Spalart1992,
  author    = {Philippe R. Spalart and Steven R. Allmaras},
  title     = {A One-Equation Turbulence Model for Aerodynamic Flows},
  journal   = {AIAA Journal},
  year      = {1992},
  volume    = {30},
  number    = {4},
  pages     = {439--446},
  doi       = {10.2514/6.1992-439}
}

@article{Zhu1997,
  author  = {Ciyou Zhu and Richard H. Byrd and Peihuang Lu and Jorge Nocedal},
  title   = {Algorithm 778: L-BFGS-B: Fortran Subroutines for Large-Scale Bound-Constrained Optimization},
  journal = {ACM Transactions on Mathematical Software},
  year    = {1997},
  volume  = {23},
  number  = {4},
  pages   = {550--560},
  doi     = {10.1145/279232.279236}
}

@book{Gill1981,
  author    = {Philip E. Gill and Walter Murray and Margaret H. Wright},
  title     = {Practical Optimization},
  publisher = {SIAM},
  year      = {2019},
  doi       = {10.1137/1.9781611975604}
}

@article{Nash1984,
  author  = {Stephen G. Nash},
  title   = {Newton-Type Minimization Via the Lanczos Method},
  journal = {SIAM Journal on Numerical Analysis},
  year    = {1984},
  volume  = {21},
  number  = {4},
  pages   = {770--788},
  doi     = {10.1137/0721052}
}

@book{Conn2000,
  author    = {Andrew R. Conn and Nicholas I. M. Gould and Philippe L. Toint},
  title     = {Trust Region Methods},
  publisher = {Society for Industrial and Applied Mathematics (SIAM)},
  year      = {2000},
  address   = {Philadelphia, PA},
  isbn      = {978-0-89871-460-9},
  doi       = {10.1137/1.9780898719857}
}

@article{Hicks1978,
  author    = {Raymond M. Hicks and Peter A. Henne},
  title     = {Wing Design by Numerical Optimization},
  journal   = {Journal of Aircraft},
  volume    = {15},
  number    = {7},
  pages     = {407--412},
  year      = {1978},
  doi       = {10.2514/3.58379}
}
\end{document}